\documentclass[12pt,a4paper,reqno]{amsart}
\usepackage[T1]{fontenc}
\usepackage{lmodern}
\usepackage{microtype}
\usepackage{setspace}
\usepackage{amsmath,amssymb,amsthm,mathtools,mathrsfs}
\usepackage[margin=1.3in]{geometry}
\usepackage[colorlinks=true,allcolors=black,bookmarksnumbered=true]{hyperref}
\hypersetup{pdftitle={Global Weak Solutions to the Three-Dimensional Compressible Navier--Stokes System for gamma > 4/3}, pdfauthor={Guangyi Hong, Hongyun Peng, Changjiang Zhu}}
\makeatletter
\let\polish@tocwrite\@tocwrite
\def\@tocwrite#1#2{\polish@tocwrite{#1}{\texorpdfstring{#2}{\@currentlabelname}}}
\makeatother
\newtheoremstyle{tthreeplain}{6.5pt}{6.5pt}{\itshape}{}%
  {\bfseries}{.}{.5em}{}
\newtheoremstyle{tthreedefinition}{6.5pt}{6.5pt}{\normalfont}{}%
  {\bfseries}{.}{.5em}{}
\theoremstyle{tthreeplain}
\newtheorem{theorem}{Theorem}[section]
\newtheorem{proposition}[theorem]{Proposition}
\newtheorem{lemma}[theorem]{Lemma}

\theoremstyle{tthreedefinition}
\newtheorem{definition}[theorem]{Definition}
\newtheorem{remark}[theorem]{Remark}
\providecommand{\T}{}
\renewcommand{\T}{{\mathbb T^3}}
\newcommand{\R}{\mathbb R}
\newcommand{\dd}{\,\mathrm d}
\newcommand{\Div}{\operatorname{div}}
\newcommand{\tr}{\operatorname{tr}}
\newcommand{\supp}{\operatorname{supp}}
\newcommand{\Id}{\mathrm{Id}}
\newcommand{\av}[1]{\langle #1\rangle}
\newcommand{\esssup}{\operatorname*{ess\,sup}}
\newcommand{\ind}{\mathbf1}
\numberwithin{equation}{section}
\allowdisplaybreaks[2]
\makeatletter
\renewenvironment{abstract}{%
  \global\setbox\abstractbox=\vbox\bgroup
    \normalfont\small\noindent\strut
    \textbf{Abstract.}\enspace\ignorespaces
}{%
    \unskip\strut\par
  \egroup
}
\def\@setkeywords{{\itshape Keywords:}\enspace\@keywords\@addpunct.}
\def\ps@headings{%
  \let\@mkboth\@gobbletwo
  \def\@evenhead{\normalfont\scriptsize\thepage\hfil}%
  \def\@oddhead{\normalfont\scriptsize\hfil\thepage}%
  \let\@evenfoot\@empty
  \let\@oddfoot\@empty
}
\AtBeginDocument{\pagestyle{headings}}%
\def\@setaddresses{%
  \par\addvspace{19.5pt}%
  \begingroup
    \small\raggedright
    \count@=0\relax
    \def\author##1{}%
    \def\address##1##2{%
      \advance\count@\@ne
      \par\addvspace{13pt}%
      \noindent
      \ifcase\count@
       \or \textbf{Guangyi Hong}%
       \or \textbf{Hongyun Peng}%
       \or \textbf{Changjiang Zhu}%
      \fi
      \par\nobreak
      \noindent ##2\par\nobreak
    }%
    \def\email##1##2{\noindent\href{mailto:##2}{##2}\par}%
    \def\\{\par\nobreak\noindent}%
    \addresses
  \endgroup
}
\def\@maketitle{%
  \normalfont\normalsize
  \@adminfootnotes
  \@mkboth{\@nx\shortauthors}{\@nx\shorttitle}%
  \thispagestyle{empty}%
  \vspace*{26pt}%
  \begingroup
    \centering
    {\Large\bfseries\boldmath \@title\par}%
    \vspace{8mm}%
    {\large\def\and{,\space\ignorespaces}\authors\par}%
    \vspace{5mm}%
  \endgroup
  \ifvoid\abstractbox\else
    \noindent\unvbox\abstractbox\par
  \fi
  \vspace{26pt}%
}
\makeatother

\title[Compressible Navier-Stokes system]{Global renormalized finite-energy weak solutions to the three-dimensional compressible Navier-Stokes system on tori for $\gamma>4/3$}

\author{Guangyi Hong}
\address{School of Mathematics,
South China University of Technology\\
Guangzhou 510641, P. R. China}
\email{magyhong@scut.edu.cn}

\author{Hongyun Peng}
\address{School of Mathematics and Statistics,
Guangdong University of Technology\\
Guangzhou 510006, P. R. China}
\email{hypengmath@gdut.edu.cn}

\author{Changjiang Zhu}
\address{School of Mathematics,
South China University of Technology\\
Guangzhou 510641, P. R. China}
\email{machjzhu@scut.edu.cn}

\thanks{Corresponding author: Changjiang Zhu.}

\date{}

\begin{document}
\setstretch{1.08}
\setlength{\jot}{3pt}
\begin{abstract}
Building on the foundational theory of finite-energy weak solutions initiated by P.-L. Lions and further developed by E. Feireisl and collaborators, we establish the global existence of finite-energy renormalized weak solutions to the compressible Navier--Stokes equations on the three-dimensional torus, with constant viscosity coefficients and pressure law \(P(\rho)=\rho^\gamma\), for every \(\gamma>4/3\) and arbitrary finite-energy initial data. This extends the classical three-dimensional Lions--Feireisl existence theory beyond the threshold \(\gamma>3/2\), and in particular, covers the important case $ \gamma=7/5 $ for air (cf. \cite{Anderson}). Compared to the classical theory, we develop a novel approach to establish the improved pressure-integrability estimate compatible with the viscous and inertial terms, via constructing a periodic test field through barycentric optimal transport. This improved pressure-integrability estimate, combined with an effective viscous flux identity and an entropy comparison, provides the compactness needed to identify the nonlinear limits and recover the energy inequality. 
\end{abstract}

\maketitle

\section{Introduction}
\label{sec:introduction}

We consider the isentropic compressible Navier--Stokes system on the three-dimensional unit torus $\T=\R^3/\mathbb Z^3$:
\begin{equation}\label{eq:cns}
\left\{
\begin{aligned}
 &\partial_t\rho+\Div(\rho u)=0,\\
 &\partial_t(\rho u)+\Div(\rho u\otimes u)+\nabla P
 =\Div S(u), \qquad (t,x)\in(0,\infty)\times\T.
\end{aligned}
\right.
\end{equation}
Here $\rho\geq0$ denotes the density, $u\in\R^3$ the velocity, and $P=\rho^\gamma$ the pressure, whose coefficient is normalized to $1$. As in \cite{FNP}, we take into account the standard Newtonian stress tensor and viscosity conditions:
\begin{equation}\label{eq:parameters}
\begin{gathered}
S(u)=\mu\bigl(\nabla u+(\nabla u)^{\mathsf T}\bigr)
      +\lambda(\Div u)\Id,\\
\mu>0,\ \  \lambda+\frac23\mu\geq0,\ \
 \nu:=2\mu+\lambda>0.
\end{gathered}
\end{equation}
Thus $\Div S(u)=\mu\Delta u+(\lambda+\mu)\nabla\Div u$.
We write $Q_T=(0,T)\times\T$ and $m=\rho u$.

Lions \cite{Lions} established the global finite-energy weak-solution theory in three
dimensions for $\gamma\geq9/5$. His framework couples the
energy inequality with weak continuity of the effective viscous flux and
renormalization of the continuity equation. These ingredients address
different failures of weak convergence: the energy controls the size of
the approximations, whereas the flux and renormalization arguments identify
the pressure and remove density oscillations. Thus the existence theorem
also provides a compactness method for nonlinear compressible flows.
Feireisl, Novotn\'y and Petzeltov\'a 
\cite{FNP} extended the range to $\gamma>3/2$. The compactness argument works with bounded truncations of the
density and controls the oscillations before removing the truncation. Feireisl's analysis in \cite{Feireisl} for the density in $L^p$ with $ p<2 $ explains why renormalization must be recovered from this oscillation control. While the
systematic approximation and stability theory is developed in \cite{FeireislBook}.

To see how these ingredients fit together, consider approximate
solutions with the uniform energy bounds
\[
 \rho_n\in L^\infty(0,T;L^\gamma),\qquad
 u_n\in L^2(0,T;H^1),\qquad
 \sqrt{\rho_n}u_n\in L^\infty(0,T;L^2).
\]
On the torus, the velocity mean is controlled using the conserved mass
and the weighted kinetic energy; see Lemma~\ref{lem:mean}. These bounds
give weakly convergent subsequences. However, they do not give spatial derivatives
of the density or identify a weak limit of the pressure. The first
additional step is therefore to establish an integrability estimate for the density function via
testing the momentum equation. The classical test function reads
\begin{gather}\label{class-test}
\displaystyle  V=\nabla\Delta_\T^{-1}
\bigl(\rho^\theta- \langle \rho^\theta \rangle\bigr),
\qquad
\Div V=\rho^\theta- \langle \rho^\theta \rangle.
\end{gather}
Hereafter $ \langle f \rangle:=\int _{\mathbb{T}^{3}}f dx$. For $P=\rho^\gamma$ and a smooth time cutoff
$0\leq\eta\in C_c^\infty(0,T)$, its pressure pairing is
\[
 \int_{Q_T}\eta\rho^\gamma\Div V
 =\int_{Q_T}\eta\rho^{\gamma+\theta}
  -\int_0^T\eta \langle \rho^\theta \rangle
                   \int_\T\rho^\gamma\,dx\,dt.
\]
The mean correction is controlled by energy when $\theta\leq\gamma$.
The same test also produces the unsigned
convective term $\int(\rho u\otimes u):\nabla V$. With
$\rho\in L^\infty_tL^\gamma_x$ and $u\in L^2_tL^6_x$, direct
H\"older estimates require
\[
\frac1\gamma+\frac13+\frac{\theta}{\gamma}\leq 1.
\]
This, together with the viscous term, gives the classical range (cf. \cite{FeireislBook,FNP}):
\begin{equation}\label{eq:classical-theta-window}
0<\theta<\min\left\{\frac\gamma2,\frac{2\gamma}{3}-1\right\},
\qquad \gamma>\frac32,
\end{equation}
Clearly, \eqref{eq:classical-theta-window} leaves room for a positive $ \theta $ only when
$\gamma>3/2$. The time derivative of the test is handled through the
renormalized continuity equation; it is not an independent regularity
bound on $\partial_t\rho$.

The effective viscous flux combines the pressure with the compressive part
of the viscous stress. For the coefficients in \eqref{eq:parameters}, it
is $P-\nu\Div u$. Although neither factor in its product with a density
truncation need converge strongly, the momentum and continuity equations
provide a compensated weak-continuity relation
\cite{Feireisl,FNP,Lions}. On the torus one tests with a gradient
inverse Laplacian, retaining the zero-mode correction. Combining the
time derivative with the continuity equation leaves a commutator of
the velocity and the periodic Riesz transforms. Compactness of this
commutator, rather than separate strong convergence of pressure and
velocity divergence, is what the flux calculation supplies. Monotonicity of the pressure then gives a sign
to the remaining density oscillation. The renormalized equation, whose
commutator foundation is due to DiPerna and Lions \cite{DL}, converts
this information into an entropy comparison. A Young measure records
possible oscillations of the approximate densities, see \cite{Ball}. And the strict
convexity of the entropy identifies when this measure is a Dirac mass.
Concentration of an unbounded pressure product requires a separate
integrability argument and is not ruled out merely by introducing a
Young measure. More precisely, writing
$\rho_n^\gamma\rightharpoonup\overline{\rho^\gamma}$ does not yet show
$\overline{\rho^\gamma}=\rho^\gamma$. Hereafter $ \overline{g} $ denotes the weak limit of $ g $. Therefore one has to compare the initialized
renormalized equations for the sequence and the limit. It is important that the pressure
covariance has a favorable sign, while the convex entropy defect starts
from zero. For densities below $L^2$, this comparison must first be
performed with bounded density truncations and their associated
entropies, and the weak source in the limiting truncated equation is kept
until the oscillation estimate permits its removal. These are the
classical compactness mechanisms used, after the required integrability of the density function has been achieved for $ \gamma>3/2 $ in the classical theory.

Several results isolate parts of this difficulty, and thus relax the constraint upon $ \gamma $. Under spherical
symmetry, Jiang and Zhang 
\cite{JiangZhang} proved global weak existence for the Cauchy
problem in $\R^d$, $d=2,3$, with $\gamma>1$. The reduction supplied by symmetry is essential and
does not yield a general three-dimensional theorem. Hu 
\cite{Hu} constructed
finite-energy renormalized solutions outside a closed exceptional set
of zero parabolic Hausdorff measure for $6/5<\gamma\leq3/2$; the result controls concentration away from that set. At the
critical exponent $\gamma=3/2$, Plotnikov \cite{Plotnikov} obtained local logarithmic
kinetic-energy estimates for a regularized problem.
The background boundary-value formulation uses a smooth bounded domain
and regular initial data with positive bounded density. The local
estimate itself is expressed through the energy and local parameters,
independently of those initial and boundary conditions.
Such an estimate excludes the associated local concentration mechanism,
but is distinct from identification of every nonlinear limit in a global
finite-energy construction.

Recently, Bresch and Jabin \cite{BJ} developed a quantitative compactness method for the
continuity equation coupled to the momentum equation. Their
work treats thermodynamically unstable pressure laws and anisotropic
viscous stresses, for which the usual monotone-pressure effective-flux
argument does not apply in its original form. It is a separate advance
in compactness theory, rather than an unconditional low-exponent
existence theorem for the constant-viscosity system considered here.
These developments show that pressure integrability, density oscillation
control and convergence of the inertial terms are connected but distinct
tasks.

There is also a huge literature concerning the regularity of solutions or the existence of strong/classical solutions to \eqref{eq:cns} under stronger
assumptions upon the initial data. Serrin \cite{Serrin} studied uniqueness of sufficiently regular
compressible motions, Nash \cite{Nash} established a local Cauchy
theory, and Matsumura and Nishida 
\cite{MatsumuraNishida} obtained global solutions
for small smooth perturbations of a nonvacuum equilibrium. In the presence of vacuum, Huang, Li and Xin \cite{HuangLiXin} proved global classical well-posedness for smooth data of small energy,
allowing large oscillations and either vacuum or nonvacuum far-field
states. Large oscillations in this result do not mean
arbitrary initial energy. Li and Xin \cite{LiXin} obtained decay estimates in two and
three dimensions and global classical existence in two dimensions for
small-energy smooth data with far-field vacuum.
Hong, Hou, Peng and Zhu \cite{HongHouPengZhu} treated a class of large-energy classical
solutions under a near-isothermal condition and a vacuum or sufficiently
small far-field density. The near-isothermal
restriction is different from small initial energy. The blow-up results of
Xin 
\cite{Xin1998} and Xin--Yan 
\cite{XinYan2013} concern smooth solutions, with compactly supported
density or an isolated mass group in the models they study. The results in \cite{Xin1998,XinYan2013} do not exclude continuation in the weak
solution class used below.

Hoff 
\cite{Hoff1995,Hoff1995Poly,Hoff1997} established the intermediate-regularity theory for the compressible Navier--Stokes equations which allows discontinuities in the density while retaining estimates beyond the basic energy inequality. Specifically, for the isothermal and polytropic constructions, the initial density is
a small perturbation of a positive constant in $L^2\cap L^\infty$, and
the velocity has small $L^2$ norm together with additional integrability. These restrictions are substantive relative
to arbitrary finite-energy data with vacuum. One of the main ideas in \cite{Hoff1995,Hoff1995Poly,Hoff1997} is to 
use the effective flux and estimates of the material derivative to
separate the rough density from better controlled velocity quantities.
The analysis of singularity curves \cite{Hoff2002} further describes how
density and velocity-gradient jumps are transported. The uniqueness
and continuous-dependence result is established in \cite{Hoff2006} by comparing corresponding particles in Lagrangian coordinates. This result applies to a weak-solution
class with additional regularity and structural bounds, but it does not
assert unconditional uniqueness for all Lions--Feireisl finite-energy solutions.

Density-dependent viscosities furnish another mechanism. Vasseur and Yu
\cite{VasseurYu} and Li and Xin \cite{LiXinDegenerate} constructed
global weak solutions for degenerate viscosity laws. The entropy structure developed in the viscous shallow-water setting
by Bresch and Desjardins in \cite{BD2003} supplies weighted
density-gradient information. Mellet and Vasseur \cite{MelletVasseur} established an
additional logarithmic kinetic estimate and a stability theory for
degenerate barotropic flows. This estimate addresses
compactness of the kinetic quantities. Bresch, Vasseur and Yu 
\cite{BreschVasseurYu} extended the entropy-weak theory to
nonlinear density-dependent viscosities satisfying the BD relation. These estimates depend on the constitutive
relation between the two viscosity coefficients. For positive constant
$\mu$, that relation would impose $\lambda=-2\mu$, which is incompatible
with \eqref{eq:parameters}. The additional density-gradient estimate of
that theory is consequently unavailable in the present problem.

For the constant-viscosity, monotone-pressure system considered here, the key of matter is to derive a pressure gain from the basic energy bounds with $ \gamma<3/2 $. The test field must detect high density, have an
$L^2$ Hessian affordable by the viscous stress, and control its time and
convective contributions without a density-gradient estimate. Improving
only the pressure weight does not resolve these simultaneous demands.
This identifies the role of the geometric construction: it supplies a
new test field within the same weak momentum equation, before the
effective-flux and entropy comparison are used.

Our geometric starting point is Serre's compensated integrability theory for nonnegative symmetric divergence-free tensors, see \cite{Serre}. In gas dynamics the mass and momentum equations combine into a space--time tensor; its positivity allows determinant estimates to exploit the
inertial terms together. Serre's periodic extension in \cite{Serre2024}
gives estimates involving $\rho^{1/3}P$ for three-dimensional inviscid
gas flows. For the Navier--Stokes equations, placing the viscous stress
inside this tensor destroys its positivity. Leaving the stress outside
requires control of its divergence that is not supplied by the energy
estimate. More recently, Serre \cite{Serre2026} considered a mildly viscous model, 
relaxed symmetry constrtaint to positivity of the symmetric part, and obtained the
corresponding estimate in two dimensions. The three-dimensional case
of that model remains open there. Thus none of these results directly
provides the pressure estimate for \eqref{eq:cns}. 

Optimal transport suggests how to retain the geometric test while
controlling its viscous cost. Brenier's polar factorization identifies
convex gradients in quadratic transport, see \cite{Brenier}. A classical
transport map assigns one target point to each source point, but does
not in general give the uniform Hessian cap needed for pairing with an
$L^2$ stress. Transport planning permits a source point to distribute
mass according to a probability measure. Alibert, Bouchitt\'e and
Champion 
\cite{ABC} study costs depending nonlinearly on that measure, and Gozlan, Roberto, Samson and Tetali \cite{GRST} establish a general
Kantorovich duality framework. For barycentric quadratic
costs, Gozlan and Juillet \cite{GJ} relate the conditional barycenter to a convex
gradient followed by a martingale coupling. The barycenter
and the target variable therefore have different roles: the former is a deterministic field used in a differential equation, while the coupling preserves the prescribed target marginal. Strong convexity
and subgradient monotonicity in Hilbert space provide the analytic control of the optimizer, see \cite{BC}.

We use these ideas to construct a periodic field with a controlled Hessian and a density-sensitive trace. The exact target marginal is retained at the level of plans, whereas a residual transport map identifies the part of the source on which the plan becomes deterministic. The area formula is applied only to this restricted part. On its complement, saturation of the Hessian bound supplies the trace estimate directly. Periodic lifts, the exact-marginal subgradient representation, the restricted Jacobian bound and the material pairing are established below; they are additional assertions needed here, rather than conclusions quoted from the transport literature.

The resulting field treats the time and convective contributions in the momentum equation together. After the kinetic-energy correction
is separated, optimizer monotonicity makes this material contribution nonnegative. The mean Hessian and its upper bound give an $L^2$ cost of order $K^{1/2}$, while its trace detects density at scale $K$.
Balancing these two estimates yields a uniform pressure tail. This
gain supplies an integrability exponent strictly above $3/2$ whenever
$\gamma>4/3$. The effective-flux and entropy arguments can then be
completed on the torus, with no positivity assumption on the density
and no smallness condition on the initial energy.

\subsection{Main results}
\label{sec:weak-formulation}\label{sec:main-results}

Assume that the initial density $\rho_0$ and initial momentum $m_0$ satisfy
\begin{equation}\label{eq:data}
\rho_0\geq0,\qquad
\rho_0\in L^\gamma(\T),\qquad
m_0=0\quad\text{on }\{\rho_0=0\},
\end{equation}
and have finite initial total energy:
\begin{equation}\label{eq:initial-energy}
E_0:=\int_\T
\left(
\frac{|m_0|^2}{2\rho_0}
+\frac{\rho_0^\gamma}{\gamma-1}
\right)\,dx<\infty,
\end{equation}
where the initial kinetic-energy density is defined to be zero where $\rho_0=0$.

\begin{definition}[Renormalized weak solution, see \cite{Lions}]\label{def:weak-solution}
A pair $(\rho,u)$ is a renormalized weak solution of \eqref{eq:cns} with the prescribed initial data if the following conditions hold. The density satisfies $\rho\geq0$, and, for every finite $T>0$,
\[
\rho\in L^\infty(0,T;L^\gamma),\qquad
u\in L^2(0,T;H^1),\qquad
\sqrt\rho\,u\in L^\infty(0,T;L^2).
\]
For all periodic test functions $\phi\in C_c^\infty([0,\infty)\times\T)$ and $\Phi\in C_c^\infty([0,\infty)\times\T;\R^3)$, the mass and momentum equations hold in the form
\begin{align}
&\int_0^\infty\!\!\int_\T
 [\rho\phi_t+m\cdot\nabla\phi]\,dx\,dt+\int_\T\rho_0\phi(0,x)\,dx=0,
 \label{eq:initialized-continuity}\\
&\begin{aligned}[t]
 &\int_0^\infty\!\!\int_\T
  [m\cdot\Phi_t+(\rho u\otimes u):\nabla\Phi\\
 &\qquad+\rho^\gamma\Div\Phi-S(u):\nabla\Phi] \,dx\,dt
  +\int_\T m_0\cdot\Phi(0,x)\,dx=0.
 \end{aligned}
 \label{eq:initialized-momentum}
\end{align}

In addition, $\rho$ satisfies the renormalized continuity equation with its initial term: for every $b\in C^1([0,\infty))$ such that $b'(s)=0$ for all sufficiently large $s$, and every $\phi\in C_c^\infty([0,\infty)\times\T)$,
\begin{equation}\label{eq:initialized-renormalization}
\begin{aligned}
&\int_0^\infty\!\!\int_\T
\Bigl[
 b(\rho)\partial_t\phi
 +b(\rho)u\cdot\nabla\phi\\
&\qquad
 -\bigl(\rho b'(\rho)-b(\rho)\bigr)
 \Div u\,\phi
\Bigr] \,dx\,dt
+\int_\T b(\rho_0)\phi(0,x)\,dx=0.
\end{aligned}
\end{equation}
The value $b(0)$ is unrestricted. The solution also satisfies the distributional energy inequality \eqref{eq:differential-main} with initial energy $E_0$ from \eqref{eq:initial-energy}.
\end{definition}

\begin{theorem}\label{thm:main}
Assume that \(\gamma>4/3\), the viscosity coefficients satisfy \eqref{eq:parameters}, and the initial data \((\rho_0,m_0)\) satisfy \eqref{eq:data}--\eqref{eq:initial-energy}. Then system \eqref{eq:cns} admits a global finite-energy renormalized weak solution \((\rho,u)\) in the sense of Definition~\ref{def:weak-solution}.
\end{theorem}

\begin{remark}
For every $T>0$, the weak solution in Theorem~\ref{thm:main} satisfies
\[
\rho\in C_w([0,T];L^\gamma(\T)),\qquad
\rho u\in C_w\left([0,T];
L^{\frac{2\gamma}{\gamma+1}}(\T)\right).
\]
The initial conditions
\[
\rho(0)=\rho_0,\qquad (\rho u)(0)=m_0
\]
hold in these weakly continuous spaces.
Moreover, the total mass and total momentum are conserved for every $t\in[0,T]$:
\[
\int_\T\rho(t,x)\,dx=\int_\T\rho_0(x)\,dx,\qquad
\int_\T\rho u(t,x)\,dx=\int_\T m_0(x)\,dx.
\]
\end{remark}

\vspace{4mm}

\begin{remark}
The weak solution in Theorem~\ref{thm:main} satisfies the following energy inequalities. Set
\[
\begin{aligned}
 E(t)&=\int_\T
 \left(\frac{|\rho u|^2}{2\rho}+\frac{\rho^\gamma}{\gamma-1}\right)(t,x)\,dx,\\
 \mathcal D(t)&=\int_\T
 \left(\mu|\nabla u|^2+(\lambda+\mu)|\Div u|^2\right)(t,x)\,dx.
\end{aligned}
\]
The kinetic-energy density is defined to be zero on vacuum. For every nonnegative $\psi\in C_c^\infty([0,T))$,
\begin{equation}\label{eq:differential-main}
 -\int_0^T E(t)\psi'(t)\,dt
 +\int_0^T\mathcal D(t)\psi(t)\,dt\leq E_0\psi(0).
\end{equation}
In particular, for every $t\in[0,T]$,
\[
 E(t)+\int_0^t\mathcal D(\tau)\,d\tau\leq E_0.
\]
\end{remark}

\begin{remark}\label{rem-bounded}
With the constraint $ \gamma>4/3 $, the result in Theorem \ref{thm:main} holds also for the case that the fluid domain is $ \mathbb{R}^{3} $ or a bounded smooth domain in $ \mathbb{R}^{3} $, since one can construct the corresponding test fields by slight modifications. 
\end{remark}

\subsection{Proof strategy}
\label{sec:proof-outline}

The proof has two distinct compactness stages. At fixed artificial-pressure
coefficient the standard construction supplies enough integrability to use
the geometric test. The estimates used to remove that coefficient must
then depend only on mass, energy, and viscous dissipation. In particular,
the enhanced density norm that makes a test admissible is not part of its
uniform bound.

We first construct global weak solutions with the cubic artificial pressure
$P_e(\rho)=\rho^\gamma+e\rho^3$ for fixed $e>0$. A pressure estimate based
on barycentric transport then allows us to pass to the limit $e\downarrow0$.

\subsubsection{The three-level approximation}

Fix $e>0$, $G>\max\{8,\gamma+1\}$, $0<\alpha\leq1$, and $\sigma>0$.
For smooth positive initial density $\rho_*$ and smooth velocity $u_*$,
consider the periodic approximate system
\begin{equation}\label{eq:strategy-approximate-system}
\left\{
\begin{aligned}
&\partial_t \rho_N+\Div(\rho_Nu_N)=\alpha\Delta \rho_N,\\
&\partial_t(\rho_Nu_N)+\Div(\rho_Nu_N\otimes u_N)
+\nabla\widehat P(\rho_N)=\Div S(u_N)
-\alpha\sum_{j=1}^3(\partial_j\rho_N)\,\partial_ju_N,
\end{aligned}
\right.
\end{equation}
where
\[
\widehat P(\rho)=\rho^\gamma+e\rho^3+\sigma \rho^G,
\qquad (\rho_N,u_N)|_{t=0}=(\rho_*,u_N^0).
\]
Here $u_N\in X_N$, the Fourier--Galerkin space containing the constant
modes, and $u_N^0$ is the $L^2(\rho_*\,dx)$-orthogonal projection of $u_*$.
The momentum equation is projected onto $X_N$, while the density equation
is not. The artificial-pressure and density-diffusion framework follows
\cite{FeireislBook,FNP,Lions}. The added pressure improves density
integrability at each construction level. The diffusion term makes the
density equation parabolic, and its momentum correction cancels the
additional kinetic-energy term; the cancellation is written out in
Appendix~\ref{app:construction}. The choice of a cubic pressure at the
last fixed level provides exactly the memberwise bound needed to pair
$m$ with $\partial_tV$. No priority claim is made for this exponent or
for artificial pressure itself. Successive
limits $N\to\infty$, $\alpha\downarrow0$, and $\sigma\downarrow0$, with $e$
fixed, yield a global renormalized weak solution for $P_e$, including the
initial terms; see Appendix~\ref{app:construction}.

We then take $e=e_n\downarrow0$, with initial data from
\eqref{eq:periodic-preparation} and Lemma~\ref{prop:prepared-family}.
The energy bounds are uniform in $n$, but the $L^\infty_tL^3_x$ density
bounds are not. This requires a pressure estimate independent of $e_n$.

\subsubsection{From classical tests to barycentric relaxation}

The classical test in \eqref{class-test}~(with $ \theta $ obeys \eqref{eq:classical-theta-window}) locates the
loss of integrability. The geometric construction below instead arranges
that the time and convective contributions, after the kinetic correction,
have a favorable combined sign. Its Hessian must simultaneously be
affordable in the stress pairing.

Serre's compensated integrability suggests treating inertia through the
nonnegative space--time tensor associated with the Euler equations:
\begin{equation}\label{eq:positive-tensor-motivation}
\mathbb M_+
=\begin{pmatrix}
\rho&m^{\mathsf T}\\
m&\rho u\otimes u+P\Id
\end{pmatrix}
\succeq0,\ \ \
\operatorname{Div}_{t,x}\mathbb M_+=0,\ \ \
(\det\mathbb M_+)^{1/3}=\rho^{1/3}P.
\end{equation}
In the smooth setting, put $g=(\det\mathbb M_+)^{1/3}$ on a bounded
space--time domain $D\subset\R^4$. A convex potential satisfying
$\det D^2\Psi=g/\int_Dg$ gives, by the matrix arithmetic--geometric
mean inequality,
\[
\mathbb M_+:D^2\Psi\geq
4\left(\int_Dg\right)^{-1/4}g.
\]
Integration by parts converts the left-hand side into a boundary flux,
controlling $\rho^{1/3}P$, see \cite{Serre,Serre2024}.

For the Navier--Stokes equations,
\[
\operatorname{Div}_{t,x}\mathbb M_+
=\binom{0}{\Div S},
\]
and the same test incurs the viscous cost
\begin{equation}\label{eq:serre-ma-viscous-cost}
\left|\int_D S:D_x^2\Psi\,dz\right|
\leq\|S\|_{L^2(D)}
\|D_x^2\Psi\|_{L^2(D)}.
\end{equation}
The Monge--Amp\`ere determinant does not control the $L^2$ norm of the
Hessian. We therefore relax the transport constraint, retaining
a density-sensitive trace bound while controlling the Hessian.

Following weak and barycentric transport in \cite{GJ,GRST}, particularly the
mixed quadratic cost of Alibert--Bouchitt\'e--Champion (see \cite[Section~5.2]{ABC}), we allow mass splitting and use the conditional
barycenter as the test map. For $\rho\geq0$ and a fixed background
$\alpha_{\mathrm{bg}}>0$, set
\[
c=\int_\T(\rho+\alpha_{\mathrm{bg}})\,dx,
\qquad
\sigma=\frac{\rho+\alpha_{\mathrm{bg}}}{c}.
\]
Here $\sigma$ is the normalized transport density, distinct from the
auxiliary-pressure coefficient in \eqref{eq:strategy-approximate-system}.
For plans $\pi$ with source marginal $(\rho+\alpha_{\mathrm{bg}})\,dx$
and target marginal $c\,dx$, write
$\bar z_\pi(x)=\mathbb E_\pi[z\mid x]$. With periodic integer lifts, the
cost is
\[
\mathcal J_K(\pi)=\frac12\int\left[
\frac{K-1}{K}|z-x|^2
+\frac1K|\bar z_\pi(x)-x|^2\right]d\pi,
\qquad K\geq6.
\]
The full plan has the prescribed target marginal, while its barycenter
need not. The construction in Section~\ref{sec:geometry} gives a unique
periodic barycentric map, with the additional properties proved here,
\[
\left\{
\begin{alignedat}{2}
T&=\nabla\varphi=x+V, \qquad
& H&=DT=D^2\varphi,\\
0&\preceq H\preceq K\Id, \qquad
& \int_\T H\,dx&=\Id.
\end{alignedat}
\right.
\]
Since $V$ is uniformly bounded, $\nabla V=H-\Id$, and
$H^2\preceq KH$, we obtain
\begin{equation}\label{eq:stress-duality-motivation}
\left\{
\begin{aligned}
\int_\T |H-\Id|^2\,dx
&\leq 3(K-1),\\[2mm]
\left|\int_{Q_T}S:\nabla V\,dx\,dt\right|
&\leq \sqrt{3T(K-1)}\,\|S\|_{L^2(Q_T)}.
\end{aligned}
\right.
\end{equation}
To recover the pressure weight without the Jacobian identity,
use the residual map $y=(Kx-T)/(K-1)$. On the uncapped region the area
formula gives $\det H\geq\sigma$; on its complement an eigenvalue of $H$
equals $K$. Hence
\begin{equation}\label{eq:strategy-trace}
\tr H\geq\min\{K,3\sigma^{1/3}\},
\end{equation}
as proved in Lemma~\ref{lem:periodic-optimizer}.

The time-dependent field requires a further sign property. Along a
smooth particle trajectory $\dot X=u(t,X)$,
\[
\frac{d}{dt}T(t,X(t))
=\partial_tV(t,X(t))+H(t,X(t))u(t,X(t)).
\]
Applying the optimizer's monotonicity to the moving mass and stationary
background gives
\begin{equation}\label{eq:design-criteria}
\mathcal M(\eta):=\int_{Q_T}\eta
\bigl[m\cdot\partial_tV+(\rho u\otimes u):H\bigr]\geq0,
\qquad 0\leq\eta\in C_c^\infty(0,T).
\end{equation}
The convective term in the momentum equation involves $H-\Id$, so the
kinetic-energy contribution is separated before applying this sign.
Lemma~\ref{lem:periodic-material} justifies the pairing for the
approximate solutions.

\subsubsection{Pressure gain and the final limit}

For $P_n=\rho_n^\gamma+e_n\rho_n^3$, testing the momentum equation with
$\eta V$ gives
\begin{equation}\label{eq:strategy-momentum-identity}
\begin{aligned}
&\int_{Q_T}\eta P_n(\tr H-3)+\mathcal M_n(\eta)\\
&\qquad=-\int_{Q_T}\eta' m_n\cdot V
+\int_{Q_T}\eta\rho_n|u_n|^2
+\int_{Q_T}\eta S_n:(H-\Id).
\end{aligned}
\end{equation}
The material term is nonnegative, the remaining inertial terms are
controlled by energy, and \eqref{eq:stress-duality-motivation} gives
\[
\left|\int_{Q_T}\eta S_n:(H-\Id)\right|
\leq C_T\|\eta\|_\infty K^{1/2}.
\]
Consequently,
\[
\int_{Q_T}\eta P_n\tr H
\leq C_T\bigl(\|\eta'\|_{L^1}+\|\eta\|_\infty\bigr)
(1+K^{1/2}).
\]
Combining this bound with \eqref{eq:strategy-trace}, then removing the
time cutoff and choosing the cap height, yields
\begin{equation}\label{eq:strategy-cap}
\sup_n\int_{Q_T}P_n\min\{\delta^{-1},\rho_n^{1/3}\}
\leq C_T\delta^{-1/2},
\qquad 0<\delta\leq1.
\end{equation}
Taking $\delta=R^{-1/3}$ for $R\geq1$ gives
\begin{equation}\label{eq:strategy-tail}
\sup_n\int_{Q_T\cap\{\rho_n>R\}}P_n
\leq C_TR^{-1/6}.
\end{equation}
In particular,
\[
|Q_T\cap\{\rho_n>R\}|
\leq C_TR^{-\gamma-1/6},
\]
so $\rho_n$ is uniformly bounded in $L^r(Q_T)$ for
$1\leq r<\gamma+1/6$. The detailed pressure estimate and the
artificial-pressure limit are proved in Lemmas~\ref{lem:pressure-test}
and \ref{lem:pressure-consequences}.

When $\gamma>4/3$, we may choose
\begin{equation}\label{eq:strict}
\left\{
\begin{aligned}
\frac32<q
&<\min\left\{2,\gamma+\frac16\right\},\\
\kappa(q)
&=\frac{2q-3}{q+1}>0,\\
p(q)
&=\frac{6(q+1)}{q+11}>\frac65.
\end{aligned}
\right.
\end{equation}
These exponents control the high-frequency effective-flux remainder and
give strong momentum compactness. To identify the pressure limit, the
spatially integrated effective viscous flux is combined with the
renormalized continuity equation. For the truncation $T_k$ and its convex
entropy $L_k$, the resulting defect relation is
\begin{equation}\label{eq:prelim-entropy-closure}
\mathscr D_k'(t)+\frac1\nu\int_\T C_k\,dx
=\int_\T Q_k\Div u\,dx,
\qquad \mathscr D_k(0)=0.
\end{equation}
Here $\mathscr D_k\geq0$ is the truncated entropy defect,
$C_k\geq0$ is the pressure--truncation covariance, and
$Q_k=T_k(\rho)-\langle T_k(z)\rangle$. Since $Q_k\to0$ in
$L^2(Q_T)$, the common initial entropy and strict convexity of
$z\log z$ eliminate density oscillations; see
Lemmas~\ref{prop:flux} and \ref{prop:renorm}.

Finally, the momentum bound with $p(q)>6/5$ and the time regularity from
the momentum equation yield $m_n\to\rho u$ in $L^2_tH^{-1}_x$.
Together with weak velocity convergence in $L^2_tH^1_x$, this identifies
convection and the kinetic-energy limit. Proposition~\ref{prop:stability}
then gives the limiting momentum equation and energy inequality.

\vspace{2mm}

The rest of this paper is organized as follows: In Section \ref{sec-appro-main-esties}, we present the approximate solutions, and state the main estimates used to prove the main theorem. In Section \ref{sec:geometry}, we establish the spatial geometry, the material-pairing sign, and the complete pressure and tail estimates. In Section \ref{sec:scalar}, we establish the desired estimates of the effective viscous flux, and then show the nonlinear convergence.

\section{Approximate Solutions and Main Estimates}\label{sec-appro-main-esties}

We state the three results used to prove Theorem~\ref{thm:main}.
Lemma~\ref{prop:prepared-family}, proved in
Appendix~\ref{app:construction}, supplies solutions with cubic
artificial pressure and the uniform energy bound \eqref{eq:budget}.
Proposition~\ref{thm:cap} gives the pressure and tail estimates; its proof
in Section~\ref{sec:geometry} consists of the barycentric construction
(see Lemma~\ref{lem:periodic-optimizer}), the material-pairing estimate
(see Lemma~\ref{lem:periodic-material}), and the momentum test
(see Lemma~\ref{lem:pressure-test}). Proposition~\ref{prop:stability} is
proved in Section~\ref{sec:scalar}: the pressure estimate first supplies
the integrability for the effective-flux and entropy arguments, which
identify the density and pressure limits; momentum compactness then
gives the convective and energy limits. We first state these results and
recall the periodic tools.

The approximation has a classical part and a requirement specific to
the geometric test. The classical part consists of a Fourier--Galerkin
momentum equation, artificial density diffusion, and artificial pressure, see \cite{FeireislBook,FNP,Lions,NovotnyStraskraba}. Here the last retained
artificial term is $e\rho^3$. An auxiliary term $\sigma\rho^G$, with
$G>\max\{8,\gamma+1\}$, is used while the Galerkin dimension and
diffusion are removed. Neither exponent is a new constitutive law for
the limiting fluid. The resulting fixed-$e$ solution must solve the full
initialized equations and have enough integrability for the
time-dependent barycentric field to be an admissible momentum test.

The order of construction in Appendix~\ref{app:construction} is as follows.
With $e,\sigma,\alpha$ fixed, the limit $N\to\infty$ uses parabolic
density compactness. The uniform $L^{5G/3}(Q_T)$ bound supplies strong
density convergence below this exponent and identifies $\widehat P$.
Strong convergence of $\nabla\rho_N$ at fixed $\alpha$ is separately
needed for the diffusion correction in momentum; it follows from
\eqref{eq:galerkin-weighted-density-gradient} and the terminal-layer
estimate, not from weak gradient convergence.

Next, with $e,\sigma$ fixed, $\alpha\downarrow0$ removes diffusion and
its momentum correction. The classical test $A\rho_\alpha$ gives a
uniform $L^{G+1}$ density bound. The effective flux and initialized
entropy comparison then give strong density convergence below $G+1$
and identify the full pressure. Finally, with $e$ fixed,
$\sigma\downarrow0$ uses the retained cubic pressure and
\eqref{eq:sigma-gain}. The auxiliary pressure vanishes strongly, and a
truncated flux and entropy argument identifies $P_e$. The constants in
this last step may depend on $e$; they are not the constants used in
Proposition~\ref{thm:cap}.

At every level the initial terms are included before taking limits.
The weighted Galerkin projection recovers the initial kinetic energy;
the smooth initial data stay fixed in the two subsequent limits.
Strong convergence of the energy densities and lower semicontinuity of
the dissipation yield the initialized energy inequality at that level.
Only after these three limits are complete do we choose $e=e_n\to0$
and use the finite-energy preparation
\eqref{eq:periodic-preparation}. This prevents the final density
compactness statement from being used to construct its own input sequence.

The hypotheses of Proposition~\ref{thm:cap} deliberately separate two
uses of regularity. For each member of the sequence,
\[
 \|\rho_nu_n\|_{L^2_tL^2_x}
 \leq\|\rho_n\|_{L^\infty_tL^3_x}
       \|u_n\|_{L^2_tL^6_x}<\infty
\]
makes the time-dependent momentum test legitimate. In the resulting
inequality this norm disappears: only $E_*$ remains. Thus the order is
to remove the smoothing at fixed $n$, establish the common pressure
budget, and only then send $n\to\infty$. Proposition~\ref{prop:stability}
uses that budget without returning to the memberwise $L^3$ bound.

\begin{lemma}\label{prop:prepared-family}
Assume \eqref{eq:data}--\eqref{eq:initial-energy}
with total mass $\mathfrak m>0$.
For any prescribed sequence $e_n>0$ with $e_n\to0$,
there exist smooth positive initial densities,
smooth initial velocities, and global periodic
approximate solutions $(\rho_n,u_n)$ with
\[
m_n=\rho_nu_n,\qquad
P_n=\rho_n^\gamma+e_n\rho_n^3.
\]
These solutions satisfy the continuity, momentum,
and renormalized continuity equations, including
their initial terms, and conserve the mass $\mathfrak m$.
Their initial data satisfy
\begin{gather}
\rho_{0n}\to\rho_0 \quad\text{in }L^\gamma,\qquad
m_{0n}\to m_0 \quad\text{in }L^{r_m},\qquad
\sqrt{\rho_{0n}}u_{0n}\to W_0 \quad\text{in }L^2,
\label{eq:prepared-data}\\
e_n\int_\T\rho_{0n}^3\,dx\to0,\qquad
E_{n0}\to E_0,
\nonumber
\end{gather}
where $r_m=2\gamma/(\gamma+1)$ and
$W_0=m_0/\sqrt{\rho_0}$, extended by zero on vacuum.

For every finite $T>0$ and nonnegative
$\psi\in C_c^\infty([0,T))$,
\[
-\int_0^T E_n\psi'\,dt
+\int_0^T\mathcal D_n\psi\,dt
\le E_{n0}\psi(0),
\]
where
\[
\begin{aligned}
E_n&=\int_\T\left(
\frac{|m_n|^2}{2\rho_n}
+\frac{\rho_n^\gamma}{\gamma-1}
+\frac{e_n}{2}\rho_n^3\right)dx,\\
\mathcal D_n&=\mu\|\nabla u_n\|_2^2
+(\lambda+\mu)\|\Div u_n\|_2^2.
\end{aligned}
\]
In particular, $\rho_n\in L^\infty(0,T;L^3(\T))$
for each $n$, and
\begin{equation}\label{eq:budget}
\begin{aligned}
\esssup_{0<t<T}\int_\T
(\rho_n+\rho_n^\gamma+\rho_n|u_n|^2+P_n)\,dx+\int_0^T
(\|\nabla u_n\|_2^2+\|u_n\|_6^2)\,dt
\le C_T,
\end{aligned}
\end{equation}
where $C_T$ is independent of $n$.
\end{lemma}

\begin{proposition}\label{thm:cap}
Let $\rho_n,P_n\geq0$ and $m_n=\rho_nu_n$ satisfy
\[
\left\{
\begin{aligned}
&\partial_t\rho_n+\Div m_n=0,\\
&\partial_tm_n+\Div(\rho_nu_n\otimes u_n)
+\nabla P_n=\Div S_n.
\end{aligned}
\right.
\]
Assume that
\[
\sup_n\left[
\esssup_{0<t<T}\int_\T
(\rho_n+\rho_n|u_n|^2+P_n)\,dx
+\int_{Q_T}|S_n|^2\,dx\,dt
\right]\le E_*,
\]
and, for each fixed $n$,
\[
\rho_n\in L^\infty(0,T;L^3(\T)),
\qquad
u_n\in L^2(0,T;H^1(\T)).
\]
Then, for some constant $C_T$ independent of $n$,
\begin{align}
\sup_n\int_{Q_T}P_n
\min\{\delta^{-1},\rho_n^{1/3}\}\,dx\,dt
&\le C_T\delta^{-1/2},
&&0<\delta\le1,
\label{eq:cap}\\
\sup_n\int_{Q_T\cap\{\rho_n>R\}}P_n\,dx\,dt
&\le C_TR^{-1/6},
&&R\ge1.
\label{eq:tail}
\end{align}
The constant depends only on $E_*$, $T$, and the torus.
\end{proposition}

\begin{proposition}\label{prop:stability}
For the approximate solutions of Lemma~\ref{prop:prepared-family}, there is a single subsequence valid on every finite time interval such that $u_n\rightharpoonup u$ in $L^2(0,T;H^1(\T))$ and the following conclusions hold.

\textup{(i)} The density Young measure is concentrated at $\rho(t,x)$ almost everywhere, that is, $Y_{t,x}=\delta_{\rho(t,x)}$, and
\begin{align*}
 \rho_n&\to\rho\quad\text{in }L^r(Q_T),
 &&1\leq r<\gamma+1/6,\\
 \rho_n^\gamma&\to\rho^\gamma\quad\text{in }L^s(Q_T),
 &&1\leq s<1+1/(6\gamma),\\
 P_n&\to\rho^\gamma\quad\text{in }L^s(Q_T),
 &&1\leq s<1+1/(6\max\{3,\gamma\}).
\end{align*}
The artificial pressure satisfies the $L^s$ convergence estimate \eqref{eq:artificial-rate}. Moreover, $m_n\to m=\rho u$ in distributions, including tests whose support meets $t=0$. The density satisfies $\rho\in C_w([0,T];L^\gamma(\T))$ and $\rho(0)=\rho_0$. The continuity and renormalized continuity equations retain the initial terms $\rho_0$ and $b(\rho_0)$, respectively. The kinetic energy satisfies
\[
 \esssup_{0<t<T}\int_\T\rho|u|^2\,dx
 \leq\sup_n\esssup_{0<t<T}\int_\T\rho_n|u_n|^2\,dx.
\]

\textup{(ii)} Along the same subsequence as in \textup{(i)},
\begin{gather*}
 m_n\to\rho u\quad\text{in }
 C([0,T];H^{-4}(\T))\cap L^2(0,T;H^{-1}(\T)),\\
 \sqrt{\rho_n}u_n\to\sqrt\rho u\quad\text{in }L^2(Q_T),\\
 \rho_nu_n\otimes u_n\to\rho u\otimes u
 \quad\text{in }L^1(Q_T).
\end{gather*}
The total-energy density converges strongly in $L^1(Q_T)$, and its spatial integral satisfies \eqref{eq:total-energy-conv}, allowing energy tests with sign-changing time derivatives. The limiting momentum equation and energy inequality retain the initial data $(\rho_0,m_0)$ and energy $E_0$, with weakly continuous representatives and the stated conservation laws.
\end{proposition}

Unless otherwise stated, spatial integrals are over $\T$. We abbreviate $L^pL^q=L^p(0,T;L^q(\T))$ and $C_tH^{-1}=C([0,T];H^{-1}(\T))$; notation with subscripts $t,x$ has the same meaning. We use $\nabla V=(\partial_jV_i)_{ij}$ and $A:B=\sum_{ij}A_{ij}B_{ij}$, with summation over repeated indices. The space $H^{-j}$ is the dual of the inhomogeneous Sobolev space $H^j$, including the constant mode. We write $\bar f=\int_\T f\,dx$ for the spatial average and $\av{F}=\int F(z)Y_{t,x}(\dd z)$ for expectation with respect to the Young measure.

The effective-viscous-flux method of \cite{Feireisl,FeireislBook,Lions}
uses a gradient right inverse of divergence. On the torus the divergence
of a periodic vector field has zero mean, so the density test must retain
its mean correction. Define the periodic inverse Laplacian to be zero
at frequency zero, and define
\[
 A_i=\partial_i\Delta_\T^{-1},\qquad
 R_{ij}=\partial_jA_i,\qquad
 \sum_iR_{ii}f=f-\bar f.
\]
The periodic Calder\'on--Zygmund estimate (see \cite{Stein}) gives
\begin{equation}\label{eq:periodic-potential-bounds}
 \|R_{ij}f\|_r+\|Af\|_{W^{1,r}}
 \leq C_r\|f\|_r,\qquad1<r<\infty.
\end{equation}
To obtain this estimate, decompose the periodic Green function locally into the Newtonian kernel and a smooth remainder. Lemma~\ref{lem:kernel} in Appendix~\ref{app:auxiliary-estimates} gives the high-frequency estimates.

Choose a smooth, nondecreasing, concave function $\mathsf T$ that equals the identity on $[0,1]$, equals $2$ on $[3,\infty)$, and satisfies $0\leq\mathsf T'\leq1$. Define
\begin{equation}\label{eq:truncations}
\left\{
\begin{alignedat}{2}
T_k(z)&=k\mathsf T(z/k), \qquad
& b_k&=T_k/k,\\
a_k(z)&=zT_k'(z)-T_k(z), \qquad
& c_k&=a_k/k.
\end{alignedat}
\right.
\end{equation}
Then $0\leq b_k\leq2$ and $-2\leq c_k\leq0$. We use $k$ for density truncation and $K$ for the Hessian bound of the convex potential; the corresponding optimizer is $T=T_{K,\alpha_{\mathrm{bg}}}$.

\begin{lemma}\label{lem:compactness}
If $r>1$, then $L^r(\T)\Subset H^{-4}(\T)$; if $p>6/5$, then $L^p(\T)\Subset H^{-1}(\T)$. A sequence with a time-uniform bound in the corresponding spatial norm and a common modulus of time continuity in the negative Sobolev space has a convergent subsequence in the corresponding space of continuous time-dependent functions. Moreover, for $6/5<p<2$ and every $\varepsilon>0$,
\[
 \|f\|_{H^{-1}}
 \leq\varepsilon\|f\|_p+C_\varepsilon\|f\|_{H^{-4}}.
\]
\end{lemma}
\begin{proof}
By Rellich's theorem \cite{Brezis},
\[
 H^4(\T)\Subset L^{r'}(\T)\quad(r>1),\qquad
 H^1(\T)\Subset L^{p'}(\T)\quad(p>6/5).
\]
Duality gives the spatial embeddings. The common time modulus and Arzel\`a--Ascoli yield compactness in time; compact interpolation gives the final inequality \cite{Simon}.
\end{proof}

\begin{lemma}\label{lem:mean}
Suppose that $r\geq0$, $\int_\T r\,dx=\mathfrak m>0$, $r\in L^{6/5}(\T)$, $\sqrt r\,v\in L^2(\T)$, and $v\in H^1(\T)$. Then
\[
 \|v\|_6\leq\mathfrak m^{-1/2}\|\sqrt r\,v\|_2
 +C\left(1+\frac{\|r\|_{6/5}}{\mathfrak m}\right)
 \|\nabla v\|_2.
\]
\end{lemma}
\begin{proof}
Set $\langle v \rangle=\int_\T v\,dx$. Applying Poincar\'e--Sobolev, weighted Cauchy--Schwarz, and H\"older's inequality yields
\begin{align}
 \|v-\langle v \rangle\|_6&\leq C\|\nabla v\|_2,
  \nonumber \\ 
  \mathfrak m|\langle v \rangle|
 &=\left|\int_\T rv\,dx-\int_\T r(v-\langle v \rangle)\,dx\right| 
  \nonumber \\ 
  &\leq\mathfrak m^{1/2}\|\sqrt r\,v\|_2
       +\|r\|_{6/5}\|v-\langle v \rangle\|_6.\label{mean-esti}
\end{align}
Dividing \eqref{mean-esti} by $\mathfrak m$ and using $\|v\|_6\leq|\langle v \rangle|+\|v-\langle v \rangle\|_6$, we then finish the proof of the present lemma.
\end{proof}

We will apply the DiPerna--Lions commutator (see \cite{DL}) to bounded
truncations with a weakly identified source. The following form keeps
that source and the initial trace explicit; it does not identify the
source with a nonlinear expression in the limiting density.

\begin{lemma}\label{lem:forced-renormalization}
Let $v\geq0$ be bounded on $Q_T$, let $u\in L^2(0,T;H^1(\T))$ and $h\in L^2(Q_T)$, and suppose that
\[
 \partial_tv+\Div(vu)+h=0,
\]
in distributions, with bounded initial data $v_0$ and with the initial term included. Then, for every $b\in C^1([0,\infty))$,
\begin{equation}\label{eq:bounded-source-chain}
 \partial_tb(v)+\Div(b(v)u)
 +(vb'(v)-b(v))\Div u+b'(v)h=0.
\end{equation}
The initial value is $b(v_0)$. The case $b(0)\ne0$ is allowed.
\end{lemma}
\begin{proof}
Extend $v=v_0$ and $u=h=0$ to negative times. The initial term makes the extended equation valid across $t=0$. Spatial convolution gives
\[
 \partial_tv_\varepsilon+\Div(v_\varepsilon u)
       +h_\varepsilon=r_\varepsilon,
\]
where
\[
 r_\varepsilon=v_\varepsilon\Div u+
 \int\nabla\omega_\varepsilon(y)\cdot
 [u(x)-u(x-y)]v(x-y)\dd y.
\]
On each finite time interval of the extension, the line-segment formula gives
\[
 \|r_\varepsilon\|_{L^2_{t,x}}
 \leq C\|v\|_{L^\infty_{t,x}}\|\nabla u\|_{L^2_{t,x}}.
\]
For smooth $u$, translation continuity of $v$ in $L^2$ and
$\int y_j\partial_i\omega_\varepsilon=-\delta_{ij}$ cancel the local linear part against $v_\varepsilon\Div u$. Then, in view of the DiPerna--Lions commutator argument in \cite{DL}, approximation of $u$ in $L^2_tH^1_x$ and the preceding bound give
$r_\varepsilon\to0$ in $L^2_{t,x}$.

For fixed $\varepsilon$, the ordinary chain rule gives
\[
 \begin{aligned}
 \partial_tb(v_\varepsilon)+\Div(b(v_\varepsilon)u)
 &+\bigl(v_\varepsilon b'(v_\varepsilon)-b(v_\varepsilon)\bigr)\Div u\\
 &+b'(v_\varepsilon)h_\varepsilon
 =b'(v_\varepsilon)r_\varepsilon.
 \end{aligned}
\]
On finite time intervals, the bounded compositions converge strongly in every finite $L^r$, and $h_\varepsilon\to h$ and $r_\varepsilon\to0$ in $L^2$. These limits yield \eqref{eq:bounded-source-chain}. The extension equals $b(v_0)$ on negative times, giving the initial term upon restriction to $(0,T)$.
\end{proof}

\section{Barycentric Transport and Pressure Estimates}
\label{sec:geometry}
This section proves Proposition~\ref{thm:cap}. Lemmas~\ref{lem:periodic-optimizer},
\ref{lem:periodic-material}, and \ref{lem:pressure-test} respectively
establish the spatial geometry, the material-pairing sign, and the complete
pressure and tail estimates.

The classical pressure test uses
\[
 V=\nabla\Delta_\T^{-1}
       (\rho^\theta- \langle \rho^\theta \rangle),
 \qquad \Div V=\rho^\theta- \langle \rho^\theta \rangle.
\]
It produces the positive term $\int\rho^{\gamma+\theta}$, but its
convective contribution requires control of
$(\rho u\otimes u):\nabla V$. The energy estimates give
$\rho u\otimes u\in L^1_tL^{3\gamma/(\gamma+3)}_x$ when the
displayed spatial exponent exceeds one. This requires $\gamma>3/2$;
the available time and space exponents then restrict $\theta$ as in
\eqref{eq:classical-theta-window}. This calculation concerns the classical
test under the energy budget; it is not an obstruction to every
possible density-dependent test.

We instead seek $T=x+V=\nabla\varphi$ with $H=DT$ satisfying
$0\preceq H\preceq K\Id$ and $\int_\T H=\Id$.
The upper bound makes the test spatially admissible, while the mean
identity improves its viscous cost to order $K^{1/2}$. These two
conditions alone do not detect high density: the additional trace
lower bound in Lemma~\ref{lem:periodic-optimizer} supplies that
information. The time-dependent construction must also respect the
material sign in Lemma~\ref{lem:periodic-material}.

\subsection{Periodic Barycentric Map}
\label{sec:barycentric-spatial}\label{dp:app-residual}
Let $\rho\in L^1(\T)$ be nonnegative and fix a background density $\alpha_{\mathrm{bg}}>0$. Set
\begin{equation}\label{eq:periodic-source}
\mathfrak m=\int_\T \rho,\qquad c=\mathfrak m+\alpha_{\mathrm{bg}},
\qquad\sigma=\frac{\rho+\alpha_{\mathrm{bg}}}{c},\qquad\int_\T\sigma=1.
\end{equation}
Since $|\T|=1$, the source $\mu=(\rho+\alpha_{\mathrm{bg}})\dd x$ and uniform target $\beta=c\dd x$ both have mass $c$. The background is part of the transport construction and allows the physical density $\rho$ to vanish. As above, $K\geq6$ bounds the derivative, while $k$ indexes the density truncation $T_k$.

For the quadratic Monge problem, a deterministic transport map sends
the source directly to the target (see \cite{Brenier}). A Jacobian constraint
then controls the product of its eigenvalues, but gives no upper
bound on individual eigenvalues. In the viscous momentum test, one needs
an $L^2$ bound for the whole Hessian. Barycentric transport separates
the prescribed target law from the field used in that test, see 
\cite{ABC,GJ,GRST}.

Write $q:\R^3\to\T$ for the quotient map. Disintegrate a plan as
$\pi(dx,dz)=\mu(dx)\pi_x(dz)$, using integer lifts of $z$.
Its two relevant quantities are
\[
 \int\pi_x(q^{-1}(B))\,\mu(dx)=c|B|,
 \qquad T(x)=\int z\,\pi_x(dz).
\]
The first identity prescribes the full target marginal. The second is a
single vector at each source point and can serve as a momentum test.
For instance, splitting between $z_1,z_2$ with weights $\theta,1-\theta$
gives $T=\theta z_1+(1-\theta)z_2$ and
\[
 \int|z-x|^2\,\pi_x(dz)
 =|T-x|^2+\int|z-T|^2\,\pi_x(dz).
\]
Thus the mixed cost in \eqref{dp:barycentric-cost} penalizes both
displacement and conditional variance while permitting splitting.
It does not impose $T_\#\mu=c\,dx$. The exact marginal will instead be
used after restricting the residual plan to its deterministic part.

The uniform background has a second role beyond normalization. The
strict positivity of $\rho+\alpha_{\mathrm{bg}}$ transfers identities
from source-almost-everywhere to Lebesgue-almost-everywhere statements.
In Lemma~\ref{lem:periodic-material}, its stationary labels will also
control $\partial_tV$ in unweighted $L^2$. None of these uses imposes a
positive lower bound on the physical density.

\begin{lemma}\label{lem:periodic-optimizer}
Let $\rho\in L^1(\T)$ be nonnegative, $\alpha_{\mathrm{bg}}>0$, and $K\geq6$. Optimization over transport plans with the periodic uniform measure as target marginal determines a unique barycentric map
\[
 T=T_{K,\alpha_{\mathrm{bg}}}[\rho]=x+V=\nabla\varphi,
 \qquad \varphi-|x|^2/2\text{ is periodic}.
\]
Here $\varphi$ is convex and $C^{1,1}$ on $\R^3$, and almost everywhere
\begin{align}
 &\begin{aligned}[t]
  &0\leq H:=DT\leq K\Id,\qquad \int_\T V=0,\\
  &\|V\|_\infty\leq D,\qquad
    D=\operatorname{diam}(\T)=\sqrt3/2,
 \end{aligned}
 \label{eq:periodic-geo-a}\\
 &\tr H\geq\min\{K,3\sigma^{1/3}\}.
 \label{eq:periodic-trace}
\end{align}
Moreover, the periodic Hessian satisfies
\begin{equation}\label{eq:periodic-hessian}
 \int_\T|H-\Id|^2\,dx\leq3(K-1).
\end{equation}
These estimates allow $\rho$ to vanish on any set.
\end{lemma}

\begin{proof}
We use Lemma~\ref{lem:periodic-envelope} and the marginal
representation in Lemma~\ref{lem:exact-marginal-barycenter}, proved in
Appendix~\ref{app:barycentric-foundations}, to construct the potential
and obtain a determinant bound where the Hessian cap is inactive.
The convex-analytic ingredients are strong convexity, subgradient
monotonicity, and conjugacy (cf. \cite{BC}). Exact periodic lifts and the
restricted determinant estimate below are additional parts of the
present construction.
	
	\medskip
	\noindent\emph{Step 1: Construction of the convex potential.}
	Let $\mu=(\rho+\alpha_{\mathrm{bg}})\,dx$ and represent its positions by
	$X\in\mathscr H=L^2(\mathcal A,\ell;\R^3)$ on an atomless label space of
	mass $c$. The convex envelope $G$ in Appendix~\ref{app:hilbert-envelope}
turns the plan problem into minimization of
\[
 G(Y)+\frac{K-1}{2}\|Y\|_{\mathscr H}^2
 -K\langle X,Y\rangle_{\mathscr H}
 +\frac K2\|X\|_{\mathscr H}^2.
\]
The quadratic part is strongly convex with modulus $K-1$; hence its
unique minimizer is characterized by a subgradient inclusion. The
exact-marginal lemma represents that subgradient by an optimal plan,
rather than merely by a limit of approximate target laws.
For clarity, the monotonicity used below is a statement on the whole
label space. If two inputs give differences $\Delta X$, $\Delta Y$,
and $\Delta\mathscr T=K\Delta X-(K-1)\Delta Y$, then
\[
 K\langle\Delta\mathscr T,\Delta X\rangle
 =\|\Delta\mathscr T\|^2
 +(K-1)\langle\Delta\mathscr T,\Delta Y\rangle
 \geq\|\Delta\mathscr T\|^2.
\]
Here the last term is nonnegative by monotonicity of $\partial G$.
This is both the $K$-Lipschitz estimate and the sign that will survive
when physical labels move and background labels remain stationary.

The exact-marginal representation requires a separate argument from
this convex-analytic inequality. Appendix~\ref{app:exact-marginal-representation}
uses the kernel topology tested by $L^1(\mathcal A;C(Q))$ to retain
the target marginal for measurable residual positions. The nonnegative
nearest-lift defect makes the lift constraint closed in that topology.
Directional derivatives then identify the supporting functional in
every bounded direction; separation gives
$\partial G(y)=\mathcal B_y$. Thus the plan used in the Jacobian
argument is an actual plan with the exact uniform target marginal.

The Hilbert-space construction gives the residual position
	$y=Y_K(X)$ and a unique periodic barycentric map $T=x+V$.
	By \eqref{dp:subgradient} and the exact-marginal representation,
	there is an optimal plan $\pi$ for which
	\begin{equation}\label{eq:periodic-residual}
	\begin{gathered}
	y(x)=\frac{Kx-T(x)}{K-1},\qquad
	T(x)=\mathbb E_\pi[z\mid x],\\
	(qy,qz)_\#\pi\text{ is optimal with target marginal }c\,dx.
	\end{gathered}
	\end{equation}
	Here $z$ is a nearest integer lift to $y(x)$, so $|z-y(x)|\le D$.
	The target marginal is exact even for bounded measurable $y$.
	
	Choose conjugate Kantorovich potentials $a,b$ for $d_\T^2/2$.
	For nearest lifts the dual inequality reads
	\[
	a(qy)+b(qz)\le\tfrac12|y-z|^2,
	\]
	with equality on the support of $\pi$. The $c$-transform therefore gives
	the convex function
	\[
	h(z)=\tfrac12|z|^2-b(qz)
	=\sup_{y\in\R^3}\bigl\{y\cdot z-\tfrac12|y|^2+a(qy)\bigr\}.
	\]
	For fixed $y(x)$, all conditional target points belong to the convex
	contact set of $h(z)-y(x)\cdot z$. Their barycenter $T(x)$ belongs to
	the same set, and hence $y(x)\in\partial h(T(x))$.
	
	Define
	\[
	h_K(z)=\frac{|z|^2}{2K}+\frac{K-1}{K}h(z),
	\qquad \varphi=h_K^*.
	\]
	Using \eqref{eq:periodic-residual}, we have
	\[
	x=\frac{T(x)}K+\frac{K-1}{K}y(x)
	\in\partial h_K(T(x)).
	\]
	Since $h_K$ is $1/K$-strongly convex, $\varphi$ is convex and
	$C^{1,1}$, with
	\begin{equation}\label{dp:spatial-gradient}
	T=\nabla\varphi,\qquad H=D^2\varphi,\qquad
	0\le H\le K\Id.
	\end{equation}
	
	The periodicity of $b$ implies
	\[
	h_K(z+n)=h_K(z)+n\cdot z+\tfrac12|n|^2,
	\qquad n\in\mathbb Z^3.
	\]
	Taking Legendre transforms shows that
	$\psi=\varphi-|x|^2/2$ is periodic. In particular, for each $i,j$,
\[
 \int_\T\partial_i\psi\,dx=0,
 \qquad \int_\T\partial_j\partial_i\psi\,dx=0.
\]
These identities hold distributionally and hence for the Lipschitz
gradient of $\psi$. They give the mean Hessian without any assumption
on the spatial mean of the fluid velocity. Consequently,
	\begin{equation}\label{eq:periodic-div-zero}
	\begin{aligned}
	V&=\nabla\psi,\qquad \nabla V=H-\Id,\qquad
	\int_\T V\,dx=0,\\
	\int_\T H\,dx&=\Id,\qquad \int_\T\tr H\,dx=3.
	\end{aligned}
	\end{equation}
	
	\medskip
	\noindent\emph{Step 2: The trace lower bound.}
The argument uses absolute continuity twice, for different measures.
First the residual map transfers the source restricted to the
unsaturated set into an absolutely continuous measure; this permits
differentiation of the transport potential. Only after the restricted
plan has become deterministic do we use its target domination to
estimate the Jacobian of $T$. Reversing these steps would assume the
determinism that the proof must establish.
	Let
	\[
	U=\{x:H(x)\text{ exists and }\lambda_{\max}(H(x))<K\}.
	\]
	On $U$ the residual map is Lipschitz and satisfies
	\begin{equation}\label{dp:residual-jacobian}
	Dy=\frac{K\Id-H}{K-1}>0\qquad\text{a.e. on }U.
	\end{equation}
The set $U$ is a measurable subset of differentiability points; it
need not be open. The area formula for Lipschitz maps
\cite[Chapter~3]{EG} applies to such restrictions. If $N\subset\T$
is null and $E=U\cap(qy)^{-1}(N)$, it gives
\[
 \int_E\det Dy(x)\,dx
 =\int_N\#\{x\in E:qy(x)=z\}\,dz=0.
\]
Here periodic charts, or equivalently integer lifts, give the same
counting formula. Since $\det Dy>0$ on $U$, $|E|=0$. Thus
$\mu(E)=0$ and $(qy)_\#(\ind_U\mu)\ll dx$.
This is the absolute continuity needed at the residual position,
before any assertion about the Jacobian of $T$.
	
	Restrict the optimal plan to source points in $U$. Its dual contact
	relation is unchanged. At almost every residual position, the
	Kantorovich potential $a$ is differentiable; at a contact pair,
	\[
	a(q(y+w))\le\tfrac12|y+w-z|^2-b(qz)
	\]
	with equality at $w=0$. Differentiation gives
	\[
	z=y-\nabla a(qy).
	\]
	The differentiability exceptional set of $a$ has zero Lebesgue
measure. The preceding residual absolute continuity makes its inverse
image negligible for the restricted source. At every remaining contact
point the displayed identity determines the Euclidean lift $z$, even
if several nearest lifts were initially allowed. Thus the target lift
is uniquely determined by the residual position.
	On $U$ the conditional plan is concentrated at its barycenter $T(x)$.
	The exact target marginal of $\pi$ therefore yields
	\begin{equation}\label{dp:restricted-target-domination}
	(qT)_\#\bigl[(\rho+\alpha_{\mathrm{bg}})\ind_U\,dx\bigr]
	\le c\,dx.
	\end{equation}
	
	Apply the area formula to the Lipschitz map $T$
	\cite[Chapter~3]{EG}. Since its pushforward in
	\eqref{dp:restricted-target-domination} is absolutely continuous,
	$\det H>0$ at almost every source point of $U$. Indeed, for
$E=U\cap\{\det H=0\}$ the area formula gives $|qT(E)|=0$.
Domination in \eqref{dp:restricted-target-domination} then gives
$\mu(E)\leq\mu(U\cap(qT)^{-1}(qT(E)))=0$, and the positive
background yields $|E|=0$. Null images may be replaced by measurable
null supersets in this argument. The density formula then implies
	\[
	\sum_{\substack{x\in U\\qT(x)=z}}
	\frac{\rho(x)+\alpha_{\mathrm{bg}}}{\det H(x)}\le c
	\qquad\text{for a.e. }z\in\T.
	\]
The domination also makes the preimage of the exceptional target-null
set a set of zero $\mu$-measure. Since $\rho+\alpha_{\mathrm{bg}}>0$,
each nonnegative summand is at most $c$ for almost every $x\in U$. Hence
	\begin{equation}\label{eq:periodic-restricted-det}
	\det H(x)\ge\frac{\rho(x)+\alpha_{\mathrm{bg}}}{c}
	=\sigma(x)\qquad\text{for a.e. }x\in U.
	\end{equation}
	The periodic area formula counts all preimages, using finitely many
	integer lifts because $T-x$ is bounded. Global injectivity is
unnecessary: every preimage contributes a nonnegative summand, so the
bound for the sum bounds each source contribution. Nor does the
argument give equality of Jacobians: the restricted target is only
dominated by $c\,dx$, and the barycentric map outside $U$ need not
push the source onto that target.
	
	On $U$, the arithmetic--geometric mean inequality gives
	$\tr H\ge3\sigma^{1/3}$. Outside $U$, an eigenvalue of $H$ equals $K$
	almost everywhere, so $\tr H\ge K$. Therefore
	\begin{equation}\label{dp:trace-detector}
	\tr H\ge\min\{K,3\sigma^{1/3}\}
	\qquad\text{a.e. on }\T.
	\end{equation}
	
	\medskip
	\noindent\emph{Step 3: The Hessian estimate.}
	From \eqref{dp:spatial-gradient}, $H^2\preceq KH$. Together with
	\eqref{eq:periodic-div-zero}, this gives
	\[
	\begin{aligned}
 \int_\T|H-\Id|^2\,dx
 &=\int_\T\bigl[\tr(H^2)-2\tr H+3\bigr]\,dx\\
 &=\int_\T\tr(H^2)\,dx-3\\
 &\leq K\int_\T\tr H\,dx-3=3(K-1).
 \end{aligned}
	\]
	This is \eqref{eq:periodic-hessian}. The matrix inequality follows
eigenvalue by eigenvalue from $0\leq\lambda_i\leq K$ and
$\lambda_i^2\leq K\lambda_i$. A pointwise bound alone would give
an $L^2$ cost of order $K$; the fixed mean trace reduces that cost
to order $K^{1/2}$. Finally,
	\eqref{dp:bounded-displacement} gives $\|V\|_\infty\le D$.
\end{proof}

The resulting estimate is independent of any density derivative.
For every stress $S\in L^2(Q_T)$ it gives the explicit payment
\[
 \left|\int_{Q_T}S:\nabla V\,dx\,dt\right|
 \leq \sqrt{3T(K-1)}\,\|S\|_{L^2(Q_T)}.
\]
Indeed, the spatial estimate holds at almost every time and
$\nabla V=H-\Id$. The factor $T^{1/2}$ comes only from integrating
that estimate in time. Time admissibility of $V$ is the next issue.

\subsection{Material Pairing and Admissibility}
\label{sec:material-pairing}\label{dp:app-material}
To make the field in Lemma~\ref{lem:periodic-optimizer} admissible
in time, we represent physical mass by moving labels and the background
by stationary labels. The time term $m\cdot\partial_tV$ by itself
has no sign: $V$ changes when the entire density distribution moves.
Monotonicity concerns the derivative along those moving positions,
which is $\partial_tV+Hu$. It therefore controls the sum with the
convective term contracted against $H$. The momentum equation instead
contains $H-\Id$; the resulting difference is the actual kinetic energy.

\begin{lemma}\label{lem:periodic-material}
	Fix one approximate solution $(\rho,u,P,S)$ from Proposition~\ref{thm:cap}, set $m=\rho u$, and choose $K\geq6$ and a time-independent $\alpha_{\mathrm{bg}}>0$. Let $V=V_{K,\alpha_{\mathrm{bg}}}[\rho]$ and $H=\Id+\nabla V$. Then $V\in H^1_{\mathrm{loc}}((0,T);L^2(\T))$, and for every $0\leq\eta\in C_c^\infty(0,T)$,
	\begin{equation}\label{eq:periodic-material-sign}
	\mathcal M_{K,\alpha_{\mathrm{bg}}}(\eta)
	:=\int_{Q_T}\eta[m\cdot\partial_tV+(\rho u\otimes u):H]\geq0.
	\end{equation}
\end{lemma}

\begin{proof}
	If the physical mass vanishes, the uniform background gives $V=0$ and
	$H=\Id$, so the assertion is immediate. We henceforth assume positive mass.
	
\medskip
\noindent\emph{Step 1: Space--time regularization.}  Choose intervals $I\Subset I'\Subset(0,T)$ with
$\supp\eta\subset I$. Let
$\chi\in C_c^\infty((-1,1))$ be nonnegative,
with $\int_{\mathbb R}\chi=1$, and set
\[
\chi_h(s)=h^{-1}\chi(s/h).
\]
For the spatial regularization, use the periodic heat kernel
\[
\Theta_h(x)=\frac{1}{(4\pi h^2)^{3/2}}
\sum_{k\in\mathbb Z^3}
\exp\left(-\frac{|x+k|^2}{4h^2}\right).
\]
Then $\Theta_h>0$ and $\int_\T\Theta_h\,dx=1$. For $F\in\{\rho,m,\rho u\otimes u,P,S\}$, define
\[
F_h(t,x)=
\int_{\mathbb R}\chi_h(s)
\int_\T\Theta_h(x-y)F(t-s,y)\,dy\,ds,
\]
where $h>0$ is sufficiently small that
$t-s\in I'$ for $t\in I$ and $s\in\supp\chi_h$.

Since convolution commutes with differentiation,
the smoothed fields satisfy
\[
\begin{cases}
\partial_t\rho_h+\Div m_h=0,\\
\partial_t m_h+\Div(\rho u\otimes u)_h
+\nabla P_h=\Div S_h.
\end{cases}
\]
The positivity of the heat kernel and conservation of
positive total mass imply $\rho_h>0$. Thus
$v_h=m_h/\rho_h$ is well defined and smooth. The nonnegative convolution kernel also preserves
the positivity of the inertia matrix:
\[
\begin{pmatrix}
\rho_h&m_h^{\mathsf T}\\
m_h&(\rho u\otimes u)_h
\end{pmatrix}\succeq0.
\]
Taking its Schur complement gives
\[
R_h:=(\rho u\otimes u)_h
-\frac{m_h\otimes m_h}{\rho_h}\succeq0.
\]

For any $a\in\R$ and $b\in\R^3$, the original matrix has
quadratic form $\rho(a+b\cdot u)^2$. Using one and the same positive
kernel for its entries preserves this inequality. More explicitly, at
each smoothing point $(t,x)$,
\[
 \begin{aligned}
 R_h(t,x)={}&\int\chi_h(s)\Theta_h(x-y)\rho(t-s,y)\\
 &\quad\times[u(t-s,y)-v_h(t,x)]
 \otimes[u(t-s,y)-v_h(t,x)]\,dy\,ds.
 \end{aligned}
\]
The cross terms cancel because the weighted mean velocity is $v_h$.
Different smoothing kernels for mass, momentum, and convection would
not give this covariance identity. The tensor $R_h$ is only a
smoothing remainder; it adds no assumption to the original equations.

	\medskip\noindent\emph{Step 2: The signed material pairing.}
	Represent $\rho_h\dd x$ by labels moving along the smooth velocity $v_h$,
	and $\alpha_{\mathrm{bg}}\dd x$ by stationary labels. Both groups enter
	the same Hilbert-space optimization \eqref{eq:periodic-hilbert-envelope}.
	The firm monotonicity of Lemma~\ref{lem:periodic-envelope}, applied along
	these paths as in \eqref{eq:periodic-path}, gives
	\[
	\langle\dot{\mathscr T}_K,\dot X\rangle_{\mathscr H}
	\geq\frac1K\|\dot{\mathscr T}_K\|_{\mathscr H}^2.
	\]
	In fact the moving labels have finite squared speed
\[
 \|\dot X(t)\|_{\mathscr H}^2
 =\int_\T\rho_h|v_h|^2\,dx
 \leq\int_\T(\rho|u|^2)_h\,dx,
\]
whereas the background labels have zero speed. Local trajectory lifts
may be chosen continuously on $\R^3$; integer-lift invariance makes
the resulting displacement independent of this choice.
Restricting the Hilbert-space $H^1$ map $\mathscr T_K(X(t))$
to the stationary background gives
$V_h\in H^1_{\mathrm{loc}}((0,T);L^2(\T))$: those labels have fixed
position, and their measure is $\alpha_{\mathrm{bg}}\,dx$.
Together with the spatial Lipschitz bound, this justifies the chain
rule on the moving labels. To see this without assuming a time-smooth
optimizer, let $\Phi_h$ be the smooth flow of $v_h$ on a compact time
interval. Its Jacobian and inverse Jacobian are bounded for this
fixed $h$. Mollify $V_h$ in space and time, apply the classical chain
rule to $V_h(t,\Phi_h(t,a))$, and pass to the limit in $L^2$ after
the change of variables $x=\Phi_h(t,a)$. The convergence follows from
the $L^2$ time derivative, the bounded spatial derivative, and the
smooth fixed-$h$ velocity. Hence
\[
 \frac{d}{dt}T_h(t,\Phi_h(t,a))
 =\bigl(\partial_tV_h+H_hv_h\bigr)(t,\Phi_h(t,a)).
\]
The stationary derivative is $\partial_tV_h$. Bounds for the flow
in this chain-rule justification may depend on $h$; they are not
used in the estimates below, which come directly from the Hilbert
inequality. Consequently,
	\begin{equation}\label{cn:dp:two-group-material}
	\begin{aligned}
	&\int_\T\left(m_h\cdot\partial_tV_h
	+\frac{m_h\otimes m_h}{\rho_h}:H_h\right)\dd x\\
	&\qquad\geq\frac1K\int_\T\rho_h
	|\partial_tV_h+H_hv_h|^2\dd x
	+\frac{\alpha_{\mathrm{bg}}}{K}
	\int_\T|\partial_tV_h|^2\dd x.
	\end{aligned}
	\end{equation}
	The background contributes nothing to the left pairing because its
velocity is zero, but contributes
$\alpha_{\mathrm{bg}}\|\partial_tV_h\|_2^2/K$ on the right. This
provides time regularity also on vacuum regions. Since $H_h\succeq0$
and $R_h\succeq0$, their contraction is nonnegative:
\[
 R_h:H_h=\tr(R_h^{1/2}H_hR_h^{1/2})\geq0.
\]
	Adding it to \eqref{cn:dp:two-group-material} and integrating against
	$\eta\geq0$ yields
	\begin{equation}\label{cn:dp:smooth-complete-sign}
	\mathcal M_h(\eta)
	\geq\frac{\alpha_{\mathrm{bg}}}{K}
	\int_{Q_T}\eta|\partial_tV_h|^2\dd x\dd t
	+\int_{Q_T}\eta R_h:H_h\dd x\dd t\geq0.
	\end{equation}
This inequality includes the covariance $R_h$ and hence the full
smoothed convective tensor $(\rho u\otimes u)_h$.
	
	The $K$-Lipschitz bound on $\mathscr T_K$ also controls the time
	derivative independently of the smoothing scale:
	\[
	\begin{aligned}
	\alpha_{\mathrm{bg}}\int_\T|\partial_tV_h|^2\dd x
	&\leq\|\dot{\mathscr T}_K\|_{\mathscr H}^2
	\leq K^2\|\dot X\|_{\mathscr H}^2\\
	&=K^2\int_\T\frac{|m_h|^2}{\rho_h}\dd x
	\leq K^2\int_\T(\rho|u|^2)_h\dd x.
	\end{aligned}
	\]
	The energy bound therefore yields a local $L^2$ bound for
	$\partial_tV_h$.
	
	\medskip\noindent\emph{Step 3: Removal of smoothing and admissibility.}
	Keep the approximate solution, $K$, and $\alpha_{\mathrm{bg}}$ fixed.
	For two nonnegative densities of equal mass, couple their common part
	identically and transport the remaining mass at distance at most $D$.
	Realizing this coupling on common atomless labels, with identical
	stationary-background labels, and applying
	\eqref{eq:periodic-firm}, we obtain
	\begin{equation}\label{eq:periodic-law-stability}
	\begin{aligned}
	\alpha_{\mathrm{bg}}
	\|V[\rho]-V[\widetilde\rho]\|_{L^2(\T)}^2
	&\leq K^2W_{2,\T}^2(\rho\,dx,\widetilde\rho\,dx)\\
	&\leq\frac{K^2D^2}{2}
	\|\rho-\widetilde\rho\|_{L^1(\T)}.
	\end{aligned}
	\end{equation}
The transport cost is unnormalized. Common smoothing preserves mass
at each time, and $\rho_h\to\rho$ locally in $L^1(Q_T)$, so
\eqref{eq:periodic-law-stability} yields $V_h\to V$ strongly in local
space--time $L^2$. The spatial Lipschitz and time-derivative bounds then give
	\[
	\partial_tV_h\rightharpoonup\partial_tV
	\quad\text{locally in }L^2,
	\qquad
	H_h\mathrel{\mathop\rightharpoonup\limits^*}H
	\quad\text{in }L^\infty.
	\]
	
	For this fixed solution, the $L^3$ density bound and Sobolev embedding
	imply
	\begin{equation}\label{cn:dp:member-momentum-L2}
	\|m\|_{L^2_tL^2_x}
	\leq\|\rho\|_{L^\infty_tL^3_x}
	\|u\|_{L^2_tL^6_x}<\infty.
	\end{equation}
	On compact time intervals, the mollified fields converge as follows:
	\[
	\begin{gathered}
	m_h\to m\quad\text{in }L^2,\
	(\rho u\otimes u)_h\to\rho u\otimes u,
	\quad P_h\to P\quad\text{in }L^1,
	\qquad S_h\to S\quad\text{in }L^2.
	\end{gathered}
	\]
Pair the strong limit of $m_h$ with the weak limit of $\partial_tV_h$,
and the strong $L^1$ limit of $(\rho u\otimes u)_h$ with the weak-star
limit of $H_h$. Passing to the limit in \eqref{cn:dp:smooth-complete-sign}
proves \eqref{eq:periodic-material-sign} and
$V\in H^1_{\mathrm{loc}}((0,T);L^2(\T))$.
For example, the temporal pairing splits as
\[
 \begin{aligned}
 \int\eta(m_h\cdot\partial_tV_h-m\cdot\partial_tV)
 &=\int\eta(m_h-m)\cdot\partial_tV_h\\
 &\quad+\int\eta m\cdot(\partial_tV_h-\partial_tV).
 \end{aligned}
\]
The first term tends to zero by the uniform smoothing-scale $L^2$
bound for $\partial_tV_h$; the second does so by weak convergence.
For the convective pairing the identical splitting uses strong $L^1$
and weak-star $L^\infty$ convergence. Thus no strong convergence of
the optimizer Hessians is required in this step.
	
To justify $\eta V$ as a momentum test, mollify the limiting field $V$
in space and time. Its time derivative converges in $L^2$ and pairs with
$m\in L^2$. Its spatial gradient converges in every finite $L^p$ and
remains bounded in $L^\infty$. An almost-everywhere convergent subsequence
then permits dominated convergence against the $L^1$ pressure and
convective tensors, while strong $L^2$ convergence handles $S\in L^2$.
This proves admissibility. The preceding sequence $H_h$ converges
weak-star; the almost-everywhere convergence here concerns mollifications
of the limiting field. All smoothing limits are taken at a fixed
approximate solution, before $n\to\infty$.
\end{proof}

\subsection{Pressure Estimate}
\label{sec:momentum-pressure}

The Hessian bound from Lemma~\ref{lem:periodic-optimizer} and the signed
pairing from Lemma~\ref{lem:periodic-material} give the pressure and
tail estimates through the weak momentum equation.

\begin{lemma}\label{lem:pressure-test}
Under the assumptions of Proposition~\ref{thm:cap}, for every finite
$T>0$ there exists $C_T$, depending only on $E_*$, $T$ and the torus,
such that for all $n$, $0<\delta\leq1$, and $R\geq1$,
\begin{align*}
 \int_{Q_T}P_n\min\{\delta^{-1},\rho_n^{1/3}\}\,dx\,dt
 &\leq C_T\delta^{-1/2},\\
 \int_{Q_T\cap\{\rho_n>R\}}P_n\,dx\,dt
 &\leq C_TR^{-1/6}.
\end{align*}
For the family in Lemma~\ref{prop:prepared-family}, these constants are
independent of $e_n$ and the individual $L^\infty_tL^3_x$ norms.
The estimates give \eqref{eq:cap}--\eqref{eq:tail} of Proposition~\ref{thm:cap}.
\end{lemma}

\begin{proof}
	Fix $n$ and $K\geq6$. Write
	$\mathfrak m_n=\int_\T\rho_n\,dx$ and choose
	\begin{equation}\label{eq:fixed-background-choice}
	\alpha_{\mathrm{bg}}=\max\{1,\sup_n\mathfrak m_n\},
	\qquad
	c_n=\mathfrak m_n+\alpha_{\mathrm{bg}},
	\qquad
	\Lambda_n=Kc_n^{1/3}.
	\end{equation}
	Then
	$1\leq\alpha_{\mathrm{bg}}\leq c_n
	\leq2\alpha_{\mathrm{bg}}$, uniformly in $n$.
	
	\medskip
\noindent\emph{Step 1: Momentum testing.} Let $0\leq\eta\in C_c^\infty(0,T)$.
	By Lemma~\ref{lem:periodic-material}, the field
	$\eta V$ is admissible in the weak momentum equation.
	Using $\nabla V=H-\Id$ and
	$\Div V=\tr H-3$, we obtain
	\begin{equation}\label{cn:dp:raw-momentum-test}
	\begin{aligned}
	0={}&\int_{Q_T}\eta' m_n\cdot V
	+\int_{Q_T}\eta m_n\cdot\partial_tV\\
	&+\int_{Q_T}\eta
	(\rho_nu_n\otimes u_n):(H-\Id)\\
	&+\int_{Q_T}\eta P_n(\tr H-3)
	-\int_{Q_T}\eta S_n:(H-\Id).
	\end{aligned}
	\end{equation}
	Since
	\[
	(\rho_nu_n\otimes u_n):(H-\Id)
	=(\rho_nu_n\otimes u_n):H-\rho_n|u_n|^2,
	\]
separating the kinetic term leaves the material pairing from
Lemma~\ref{lem:periodic-material}. Thus
	\begin{equation}\label{eq:physical-test}
	\begin{aligned}
	&\int_{Q_T}\eta P_n(\tr H-3)
	+\mathcal M_{K,\alpha_{\mathrm{bg}}}(\eta)\\
	&\quad=-\int_{Q_T}\eta'm_n\cdot V
	+\int_{Q_T}\eta\rho_n|u_n|^2
	+\int_{Q_T}\eta S_n:(H-\Id),
	\end{aligned}
	\end{equation}
	where
	$\mathcal M_{K,\alpha_{\mathrm{bg}}}(\eta)\geq0$.
	
	\medskip
\noindent\emph{Step 2: Energy and Hessian bounds.}
To estimate the right-hand side of \eqref{eq:physical-test}, use
mass conservation and Cauchy--Schwarz:
	\[
	\|m_n(t)\|_1
	\leq\mathfrak m_n^{1/2}
	\left(\int_\T\rho_n|u_n|^2\,dx\right)^{1/2}.
	\]
	Consequently, the uniform energy bounds and
	$\|V\|_\infty\leq D$ give
	\[
	\begin{aligned}
	\left|\int_{Q_T}\eta'm_n\cdot V\right|
	&\leq D\|\eta'\|_1
	\esssup_t\|m_n(t)\|_1
	\leq C_T\|\eta'\|_1,\\
	\int_{Q_T}\eta(\rho_n|u_n|^2+3P_n)
	&\leq T\|\eta\|_\infty
	\esssup_t\int_\T(\rho_n|u_n|^2+3P_n)\,dx\\
	&\leq C_T\|\eta\|_\infty.
	\end{aligned}
	\]
	
	For the viscous term, the Hessian estimate
	\eqref{eq:periodic-hessian} yields
	\begin{equation}\label{eq:periodic-stress-cost}
	\begin{aligned}
	\left|\int_{Q_T}\eta S_n:(H-\Id)\right|
	&\leq\|\eta^{1/2}S_n\|_2
	\|\eta^{1/2}(H-\Id)\|_2\\
	&\leq\sqrt{3T(K-1)}\,\|\eta\|_\infty
	\|S_n\|_{L^2(Q_T)}\\
	&\leq C_T\|\eta\|_\infty K^{1/2}.
	\end{aligned}
	\end{equation}
	Inserting these bounds into
	\eqref{eq:physical-test} and moving the term
	$-3P_n$ to the right gives
	\begin{equation}\label{eq:periodic-budget}
	\begin{aligned}
	&\int_{Q_T}\eta P_n\tr H
	+\mathcal M_{K,\alpha_{\mathrm{bg}}}(\eta)\\
	&\qquad\leq
	C_T(\|\eta'\|_1+\|\eta\|_\infty)(1+K^{1/2}).
	\end{aligned}
	\end{equation}
One can keep the budget constant explicit at this point:
\[
 \begin{aligned}
 &\int_{Q_T}\eta P_n\tr H
       +\mathcal M_{K,\alpha_{\mathrm{bg}}}(\eta)\\
 &\quad\leq D E_*\|\eta'\|_1
 +\left(3TE_*+\sqrt{3T(K-1)E_*}\right)\|\eta\|_\infty.
 \end{aligned}
\]
Indeed, both the mass and kinetic energy are at most $E_*$,
$\int(\rho_n|u_n|^2+3P_n)\leq3E_*$ at almost every time, and
$\|S_n\|_2^2\leq E_*$. Moreover,
$c_n\leq2\max\{1,E_*\}$. Thus the normalizations and the small-cap
case below only change constants through $E_*$ and $T$; neither the
individual cubic-density bound nor the time-regularity constant of
the optimizer enters the final estimate.

By the trace lower bound \eqref{eq:periodic-trace} and
	$\rho_n+\alpha_{\mathrm{bg}}\geq\rho_n$,
	\begin{equation}\label{eq:periodic-detector}
	\begin{aligned}
	\min\{\Lambda_n,\rho_n^{1/3}\}
	&\leq\min\{Kc_n^{1/3},
	3(\rho_n+\alpha_{\mathrm{bg}})^{1/3}\}\\
	&\leq c_n^{1/3}\tr H.
	\end{aligned}
	\end{equation}
	Multiplying \eqref{eq:periodic-budget} by
	$c_n^{1/3}$ and using
	$K^{1/2}=\Lambda_n^{1/2}c_n^{-1/6}$,
	we obtain
	\begin{equation}\label{eq:cap-budget}
	\begin{aligned}
	&\int_{Q_T}\eta P_n
	\min\{\Lambda_n,\rho_n^{1/3}\}
	+c_n^{1/3}\mathcal M_{K,\alpha_{\mathrm{bg}}}(\eta)\\
	&\qquad\leq
	C_T(\|\eta'\|_1+\|\eta\|_\infty)
	(1+\Lambda_n^{1/2}).
	\end{aligned}
	\end{equation}
	The constant is independent of $n$, $K$,
	and the individual $L^\infty_tL^3_x$
	bounds of $\rho_n$.
	
	\medskip
The normalized trace cap $\Lambda_n=Kc_n^{1/3}$ is the density
weight that can be detected by a given Hessian cap. Once a desired
weight $L$ is prescribed, the least admissible large cap is
$K=Lc_n^{-1/3}$. A larger cap increases the viscous cost without
improving $\min\{L,\rho_n^{1/3}\}$. Since $1\leq c_n\leq
2\alpha_{\mathrm{bg}}$, this is the choice $K\sim L$ with uniform
comparison constants. The power $L^{1/2}$ below comes from the
Hessian cost, whereas the density power $1/3$ comes from the
three-dimensional determinant inequality.

\noindent\emph{Step 3: Cutoff removal and tail decay.} Choose $0\leq\eta_\varepsilon\leq1$ in
	$C_c^\infty(0,T)$ such that
	$\eta_\varepsilon\to1$ almost everywhere
	and $\|\eta_\varepsilon'\|_1\leq2$. One may use a smooth rise from
zero to one near $t=0$, a constant middle part, and a smooth fall near
$t=T$; their total variations add to two. Their derivatives need not
converge in $L^1$. The bound in \eqref{eq:cap-budget} uses only this
total variation and therefore controls both time endpoints uniformly.
	Dropping the nonnegative material pairing
	in \eqref{eq:cap-budget} and applying
	Fatou's lemma gives
	\[
	\int_{Q_T}P_n\min\{\Lambda_n,\rho_n^{1/3}\}
	\leq C_T(1+\Lambda_n^{1/2}).
	\]
	
	Let $L=\delta^{-1}\geq1$.
	If $L\geq6c_n^{1/3}$, choose
	$K=Lc_n^{-1/3}\geq6$. Then
	$\Lambda_n=L$, and
	\[
	\int_{Q_T}P_n\min\{L,\rho_n^{1/3}\}
	\leq C_T(1+L^{1/2})
	\leq2C_TL^{1/2}.
	\]
	If $1\leq L<6c_n^{1/3}$, the energy bound
	instead gives
	\[
	\begin{aligned}
	\int_{Q_T}P_n\min\{L,\rho_n^{1/3}\}
	&\leq L\int_{Q_T}P_n\\
	&\leq C_TL
	\leq C_TL^{1/2},
	\end{aligned}
	\]
	where the last constant is uniform in $n$
	because $c_n$ is uniformly bounded.
	Combining the two cases yields
	\begin{equation}\label{cn:geo:global-cap}
	\int_{Q_T}P_n
	\min\{\delta^{-1},\rho_n^{1/3}\}
	\leq C_T\delta^{-1/2},
	\qquad 0<\delta\leq1.
	\end{equation}
	This proves \eqref{eq:cap}.
	
	Finally, for $R\geq1$, take
	$\delta=R^{-1/3}$. On $\{\rho_n>R\}$,
	the minimum in \eqref{cn:geo:global-cap}
	equals $R^{1/3}$. Hence
	\[
	\begin{aligned}
	\int_{Q_T\cap\{\rho_n>R\}}P_n
	&\leq R^{-1/3}
	\int_{Q_T}P_n\min\{R^{1/3},\rho_n^{1/3}\}\\
	&\leq C_TR^{-1/6},
	\end{aligned}
	\]
	which proves \eqref{eq:tail}. In this last step the powers are
$R^{-1/3}(R^{1/3})^{1/2}=R^{-1/6}$.
The resulting uniform pressure bound depends only on the stated
energy budget, although test admissibility used the fixed-member
regularity. The constitutive law
$P_n=\rho_n^\gamma+e_n\rho_n^3$ enters the subsequent conversion
of this pressure tail into density moments.
\end{proof}

Proposition~\ref{thm:cap} follows from Lemma~\ref{lem:pressure-test}.

\section{Effective Viscous Flux and Nonlinear Convergence}
\label{sec:scalar}

We prove Proposition~\ref{prop:stability} by combining the pressure tail
estimate with effective viscous flux and entropy comparison. The tail
estimate gives integrability; a high-frequency commutator bound then
yields the spatially integrated flux identity. Renormalization and the
common initial entropy eliminate density oscillations, after which
momentum compactness gives the convective and energy limits.

Fix a finite $T>0$ and write
\[
 Q_T=(0,T)\times\T,\qquad
 m_n=\rho_nu_n,\qquad
 P_n=\rho_n^\gamma+e_n\rho_n^3.
\]
The effective viscous flux method originates in the compactness theory of
Lions and Feireisl \cite{Feireisl,FeireislBook,Lions}. Its use here has two
separate prerequisites. The pressure tail must exclude concentration in
the pressure--truncation product, and the momentum must have enough
integrability to pair with a small high-frequency commutator. Neither
prerequisite identifies the pressure by itself. The entropy comparison
below performs that final identification.

The tail estimate \eqref{eq:tail} yields the following bounds.

\begin{lemma}
	\label{lem:pressure-consequences}
	Let $(\rho_n,u_n)$ be the family in Lemma~\ref{prop:prepared-family}, with
	$m_n=\rho_nu_n$ and $P_n=\rho_n^\gamma+e_n\rho_n^3$. Set
	$q_c=\gamma+1/6$. For every finite $T>0$, the following estimates hold
	with constants independent of $n$.
	
	\textup{(i)} For $1\leq q<q_c$ and $0<\beta<1/6$,
	\[
	\sup_n\int_{Q_T}\rho_n^q\,dx\,dt\leq C_{T,q},\qquad
	\sup_n\int_{Q_T}P_n\rho_n^\beta\,dx\,dt\leq C_{T,\beta}.
	\]
	
	\textup{(ii)} For $1<s<19/18$, the artificial pressure satisfies
	\begin{equation}\label{eq:artificial-rate}
	\|e_n\rho_n^3\|_{L^s(Q_T)}
	\leq C_{T,s}e_n^{1-1/s}.
	\end{equation}
	Moreover,
	\[
	\sup_n\|P_n\|_{L^s(Q_T)}\leq C_{T,s},
	\qquad 1<s<1+\frac1{6\max\{3,\gamma\}}.
	\]
	In particular, $e_n\rho_n^3\to0$ in the indicated $L^s$ spaces,
	and $\{P_n\}$ is uniformly integrable on $Q_T$.
	
	\textup{(iii)} Since $\gamma>4/3$, choose
	$3/2<q<\min\{2,q_c\}$, and set
	\[
	p(q)=\frac{6(q+1)}{q+11},\qquad
	\kappa(q)=\frac{2q-3}{q+1}.
	\]
	Then
	\[
	\sup_n\|m_n\|_{L^2(0,T;L^{p(q)}(\T))}\leq C_{T,q},
	\qquad p(q)>\frac65,\quad \kappa(q)>0.
	\]
	Consequently $L^{p(q)}(\T)\Subset H^{-1}(\T)$, and $\kappa(q)$
	is the positive exponent used in the high-frequency estimate.
\end{lemma}

\begin{proof}
	\medskip
\noindent\emph{Step 1: Weighted moments.} Let $P_*:=\sup_n\esssup_{0<t<T}\int_\T P_n\,dx$,
	which is finite by \eqref{eq:budget}. From $P_n\geq\rho_n^\gamma$
	and \eqref{eq:tail}, for $R\geq1$,
	\[
	R^\gamma|Q_T\cap\{\rho_n>R\}|
\leq\int_{Q_T\cap\{\rho_n>R\}}\rho_n^\gamma
\leq\int_{Q_T\cap\{\rho_n>R\}}P_n
\leq C_TR^{-1/6}.
	\]
	The distribution-function formula therefore gives
	\begin{equation}\label{eq:moments}
	\begin{aligned}
	&|Q_T\cap\{\rho_n>R\}|\leq C_TR^{-q_c},\ \ \ \ \ \ \ \ \
	R\geq1,\\
	&\int_{Q_T}\rho_n^q\,dx\,dt
	\leq T+C_Tq\int_1^\infty R^{q-1-q_c}\,dR
	=T+\frac{C_Tq}{q_c-q},\ \
	1\leq q<q_c.
	\end{aligned}
	\end{equation}
Applying the layer-cake formula with weight $P_n$ yields
	\begin{equation}\label{eq:pressure-weighted-moment-proof}
	\begin{aligned}
	\int_{Q_T}P_n\rho_n^\beta\,dx\,dt
	&\leq TP_*+\beta\int_1^\infty R^{\beta-1}
	\int_{Q_T\cap\{\rho_n>R\}}P_n\,dx\,dt\,dR\\
	&\leq TP_*+\frac{C_T\beta}{1/6-\beta},
	\qquad 0<\beta<\frac16.
	\end{aligned}
	\end{equation}
	The strict inequality $q<q_c$ is required by the integral at infinity:
at $q=q_c$ its integrand is $R^{-1}$. Thus the tail estimate supplies
all subcritical moments, without asserting a strong $L^{q_c}$ endpoint.
This proves \textup{(i)}.
	
	\medskip
\noindent\emph{Step 2: Pressure integrability.} Take $1<s<19/18$ and put $\beta=3(s-1)<1/6$.
	Since $e_n\rho_n^3\leq P_n$, estimate
	\eqref{eq:pressure-weighted-moment-proof} gives
	\[
	\begin{aligned}
	\int_{Q_T}(e_n\rho_n^3)^s\,dx\,dt
	&=e_n^{s-1}\int_{Q_T}e_n\rho_n^3\rho_n^{3(s-1)}\,dx\,dt\\
	&\leq e_n^{s-1}\int_{Q_T}P_n\rho_n^{3(s-1)}\,dx\,dt
	\leq C_{T,s}e_n^{s-1}.
	\end{aligned}
	\]
	This proves \eqref{eq:artificial-rate}. Similarly,
	\[
	\int_{Q_T}(\rho_n^\gamma)^s\,dx\,dt
	\leq\int_{Q_T}P_n\rho_n^{\gamma(s-1)}\,dx\,dt
	\leq C_{T,s},\qquad 1<s<1+\frac1{6\gamma}.
	\]
	Combining these estimates with
	$P_n^s\leq2^{s-1}[(\rho_n^\gamma)^s+(e_n\rho_n^3)^s]$
	proves the asserted uniform $L^s$ bound for $P_n$.
	Since $s>1$, this also gives uniform integrability.
	
	\medskip
\noindent\emph{Step 3: Momentum interpolation.} Fix $3/2<q<\min\{2,q_c\}$ and set
	$\vartheta=2/(q+1)$. On the vacuum set, both sides of the
	following identity are taken to be zero:
	\[
	|m_n|=|\sqrt{\rho_n}u_n|^\vartheta
	|u_n|^{1-\vartheta}\rho_n^{1-\vartheta/2}.
	\]
	The spatial H\"older exponents satisfy
	\[
	\frac1p=\frac\vartheta2+\frac{1-\vartheta}6
	+\frac{1-\vartheta/2}{q}
	=\frac{q+11}{6(q+1)},
	\]
	while $(1-\vartheta)+(2-\vartheta)/q=1$ gives the time
	exponents. Applying H\"older in space and then in time, we obtain
	\begin{equation}\label{eq:momentum-interpolation}
	\begin{aligned}
	\|m_n\|_{L^2(0,T;L^p)}\leq{}&
	\left(\esssup_t\int_\T\rho_n|u_n|^2\,dx\right)^{1/(q+1)}\\
	&\times\left(\int_0^T\|u_n\|_6^2\,dt\right)^{(q-1)/(2(q+1))}
	\left(\int_{Q_T}\rho_n^q\,dx\,dt\right)^{1/(q+1)}.
	\end{aligned}
	\end{equation}
The right-hand side is uniformly bounded in $n$ by \eqref{eq:budget}
and \eqref{eq:moments}; Lemma~\ref{lem:mean} controls velocities with
arbitrary spatial mean. Finally,
	\[
	p=\frac{6(q+1)}{q+11}>\frac65,\qquad
	\frac52-\frac3p=\kappa(q)=\frac{2q-3}{q+1}>0.
	\]
	Indeed,
\[
 p-\frac65=\frac{6(4q-6)}{5(q+11)}>0,
 \qquad q_c-\frac32=\gamma-\frac43>0.
\]
The second inequality provides a nonempty choice of $q$; the first
converts that choice into the spatial compactness needed for momentum.
By Rellich compactness and duality \cite{Brezis},
$H^1\Subset L^{p'}$ and $L^p\Subset H^{-1}$.
This completes \textup{(iii)}. The quantity $5/2-3/p$ measures the
remaining integrability available for the high-frequency kernel; its
positivity will give a power of the smoothing scale.
\end{proof}

The following high-frequency effective-flux estimate is uniform in both
the approximation index $n$ and the density truncation level $k$.
Let $\chi\in C_c^\infty(B_1)$ be even, nonnegative, and have unit integral.
For $0<a<1/8$, write
\[
S_ag=\chi_a*_{\T}g,\qquad H_a=I-S_a,\qquad
A_{a,i}=A_iH_a,\qquad L_{ij,a}=R_{ij}H_a.
\]
The inverse Laplacian acts on mean-zero functions, so that
\[
 \Div Af=f-\bar f,\qquad
 \Div A_af=H_af,\qquad S_a\bar f=\bar f.
\]
Thus $A$ tests pressure through a prescribed divergence, while $A_a$
isolates the part that must vanish with $a$. The complementary part
$S_af$ retains the spatial mean. Removing that mean from both parts
would change the flux calculation.
The finite-exponent bounds for the Riesz transforms follow from
singular-integral theory \cite{Stein}. The quantitative periodic kernel
bounds needed here are proved in Lemma~\ref{lem:kernel} of
Appendix~\ref{app:auxiliary-estimates}.

\begin{lemma}[Uniform high-frequency effective-flux remainder]
	\label{lem:uniform-high-frequency}
	Let $(\rho_n,u_n)$ be the family in Lemma~\ref{prop:prepared-family}.
	Choose $q$ as in Lemma~\ref{lem:pressure-consequences} \textup{(iii)} and put
	\[
	p=p(q)=\frac{6(q+1)}{q+11},\qquad
	\kappa(q)=\frac{2q-3}{q+1}>0.
	\]
	For every $\eta\in C_c^\infty(0,T)$, $0<a<1/8$, and truncation
	level $k\geq1$, the field $\eta A_ab_k(\rho_n)$ is an admissible test
	in the weak momentum equation of each fixed approximate solution.
	Moreover, setting
	\[
	\mathcal R_{n,k,a}
	=\int_{Q_T}\eta(P_n-\nu\Div u_n)H_ab_k(\rho_n)\,dx\,dt,
	\]
	we have
	\begin{equation}\label{eq:normalized-flux-remainder}
	\sup_{n,k}|\mathcal R_{n,k,a}|
	\leq C_{T,q}\left[
	\|\eta\|_\infty\bigl(a^{\kappa(q)}+a^2\bigr)
	+a\|\eta'\|_1\right]\longrightarrow0
	\qquad\text{as }a\downarrow0.
	\end{equation}
	The constant is independent of $n$, $k$, $a$, and the individual
	$L^\infty_tL^3_x$ bounds of $\rho_n$.
\end{lemma}

\begin{proof}
\medskip\noindent\emph{Step 1: The complete commutator bound.} Fix $n,k,a$ and write
	\[
	\begin{gathered}
	d_n=\Div u_n,\qquad
	g=b_k(\rho_n)=\frac{T_k(\rho_n)}{k},\\
	h=c_k(\rho_n)d_n=\frac{a_k(\rho_n)}{k}d_n.
	\end{gathered}
	\]
	The renormalized continuity equation takes the form
	\begin{equation}\label{eq:two-source}
	\partial_tg+\Div(gu_n)+h=0.
	\end{equation}
	Its spatial mean satisfies $\bar h=-\partial_t\bar g$.
	Since $A_a$ annihilates constants, $A_ah=A_a(h-\bar h)$;
	the source term must therefore be retained in the high-frequency pairing.
This source is produced by the DiPerna--Lions renormalization
\cite{DL}, through $zg'(z)-g(z)=c_k(z)$. It is not an error in the
continuity equation and cannot be omitted when differentiating $A_ag$.
In particular, the differentiated test involves the complete expression
\[
 \partial_t(A_ag)_i=-L_{ij,a}(gu_{n,j})-A_{a,i}h.
\]
	
	By Lemma~\ref{lem:kernel}, for $g\in L^\infty$ and
	$|h|\leq C_h|\nabla u|$,
	\begin{equation}\label{cn:stb:complete-kernel}
	\left|L_{ij,a}(gu_j)-u_jL_{ij,a}g+A_{a,i}h\right|
	\leq C(\|g\|_\infty+C_h)
	\mathcal K_a*_{\T}|\nabla u|,
	\end{equation}
	where repeated indices are summed and
	\[
	\|\mathcal K_a\|_s\leq C_s(a^{3/s-2}+a^2),
	\qquad 1\leq s<\frac32.
	\]
	Near the diagonal, the velocity difference in the commutator removes
one power of the order-zero singular kernel. Integration along the
line segment leaves the majorant
$r^{-2}(1+|\log(r/a)|)$ for $r<a$. Its $s$th power is locally
integrable exactly when $s<3/2$. Away from the diagonal, evenness of
$\chi$ cancels the first-order Taylor term, leaving $a^2r^{-4}$
after the same line-segment integration. The $a^2$ term in the periodic
bound also includes the nonsingular lattice translates. These are the
two contributions established in Lemma~\ref{lem:kernel}; no
uniform bound on a spatial derivative of $g$ is used.
Set
	\[
	\begin{aligned}
	B_{n,k,a,i}
	&=L_{ij,a}(gu_{n,j})-u_{n,j}L_{ij,a}g+A_{a,i}h,\\
	J_{n,k,a}(\eta)&=\int_{Q_T}\eta\,m_n\cdot B_{n,k,a}\,dx\,dt.
	\end{aligned}
	\]
	Let $p'$ be the conjugate exponent of $p$. The Young and H\"older
	exponents can be chosen as
	\begin{equation}\label{eq:highpass-exponents}
	\begin{gathered}
	\frac1p=\frac{q+11}{6(q+1)},\quad
	p'=\frac{6(q+1)}{5(q-1)},\quad
	s=\frac{3(q+1)}{4q-1},\\
	\frac1s=\frac32-\frac1p,\quad
	1+\frac1{p'}=\frac12+\frac1s,\quad
	\frac3s-2=\kappa(q)>0.
	\end{gathered}
	\end{equation}
	The convolution is paired first in space and then in time:
\[
 L^s_x*L^2_x\longrightarrow L^{p'}_x,
 \qquad L^2_tL^p_x\cdot L^2_tL^{p'}_x\longrightarrow L^1_{t,x}.
\]
Since $3/2<q<2$, we have $1<s<3/2$ and $p'<6$.
At $q=3/2$, one would have $s=3/2$ and $\kappa(q)=0$;
the near-diagonal integral would cease to be integrable. The strict
inequality on $q$ therefore enters both the kernel estimate and its
decay as $a\downarrow0$.
	Using $\|b_k\|_\infty,\|c_k\|_\infty\leq2$ and
	$|d_n|\leq\sqrt3|\nabla u_n|$, the kernel estimate and Young's
	inequality give
	\[
	\|B_{n,k,a}(t)\|_{p'}
	\leq C(a^{\kappa(q)}+a^2)\|\nabla u_n(t)\|_2.
	\]
	The momentum bound in Lemma~\ref{lem:pressure-consequences}
	and Cauchy--Schwarz in time now yield
	\begin{equation}\label{eq:work-small}
	\begin{aligned}
	|J_{n,k,a}(\eta)|
	&\leq C\|\eta\|_\infty
	\|m_n\|_{L^2_tL^p_x}\|\nabla u_n\|_{L^2(Q_T)}
	(a^{\kappa(q)}+a^2)\\
	&\leq C_{T,q}\|\eta\|_\infty(a^{\kappa(q)}+a^2).
	\end{aligned}
	\end{equation}
	
\medskip\noindent\emph{Step 2: Test admissibility.} For the fixed solution, the cubic-pressure bound gives
\[
 m_n\in L^2(0,T;L^2),\qquad
 \rho_nu_n\otimes u_n\in L^1(0,T;L^{3/2}).
\]
Indeed, pair $\rho_n\in L^\infty_tL^3_x$ with $u_n\in L^2_tL^6_x$
for the first bound and with $|u_n|^2\in L^1_tL^3_x$ for the second.
The resulting constants may depend on $n$.
Take $V=A_ag$. Equation~\eqref{eq:two-source} implies
\[
\left\{
\begin{aligned}
\partial_tV_i
&=-L_{ij,a}(gu_{n,j})-A_{a,i}h
\in L^2_tL^2_x,\\[1mm]
\nabla V
&=(L_{ij,a}g)_{ij}
\in L^\infty_tL^r_x,\qquad 1<r<\infty,\\[1mm]
V&\in L^\infty(Q_T),\\[1mm]
\Div V&=H_ag\in L^\infty(Q_T),\\[1mm]
\|H_ag\|_\infty&\leq4.
\end{aligned}
\right.
\]
	To use $\eta V$ in the weak momentum equation, smooth $V$ in
	space and time, and then remove these regularizations at fixed $n,k,a$.
	The time derivatives converge in $L^2$ and pair with $m_n\in L^2$.
	The gradients converge in $L^3$ almost everywhere in time and remain
	bounded in $L^\infty_tL^3_x$, which permits passage in the convective
	term. Their $L^2$ convergence handles the stress term, while the
	uniformly bounded divergences converge almost everywhere and pair
	with $P_n\in L^1$. Thus $\eta V$ is admissible. These smoothing limits are taken before
any limit in $n$. The displayed cubic-density norms may depend on $n$;
they justify testing a fixed equation. Only
\eqref{eq:work-small} and the time-cutoff estimate below are used
uniformly along the sequence.
	
\medskip\noindent\emph{Step 3: The flux remainder.} Write $V=\nabla\phi$, where $\phi=\Delta_\T^{-1}H_ag$.
	Periodic integration by parts gives
	\[
	\begin{aligned}
	\int_\T\partial_ju_{i,n}\partial_{ij}\phi\,dx
	&=-\int_\T u_{i,n}\partial_i\Delta\phi\,dx
	=\int_\T d_n\Delta\phi\,dx,\\
	\int_\T S(u_n):\nabla V\,dx
	&=\nu\int_\T d_nH_ag\,dx,
	\qquad\nu=2\mu+\lambda.
	\end{aligned}
	\]
	Moreover, the time and convective terms combine exactly:
	\begin{equation}\label{eq:complete-inertial-regroup}
	\begin{aligned}
	&m_{i,n}\partial_tV_i
	+\rho_nu_{i,n}u_{j,n}\partial_jV_i\\
	&\quad=-m_{i,n}L_{ij,a}(gu_{n,j})-m_{i,n}A_{a,i}h
	+m_{i,n}u_{j,n}L_{ij,a}g\\
	&\quad=-m_{i,n}B_{n,k,a,i}.
	\end{aligned}
	\end{equation}
	Testing the weak momentum equation with $\eta V$ therefore yields
	\begin{equation}\label{eq:periodic-work-identity}
	J_{n,k,a}(\eta)
	=\int_{Q_T}\eta(P_n-\nu d_n)H_ag\,dx\,dt
	+\int_{Q_T}\eta'm_n\cdot A_ag\,dx\,dt.
	\end{equation}
	Finally, Lemma~\ref{lem:kernel} gives
	$\|A_ag\|_\infty\leq Ca\|g\|_\infty\leq Ca$.
	Mass conservation and the kinetic-energy bound imply
	\[
	\left|\int_{Q_T}\eta'm_n\cdot A_ag\,dx\,dt\right|
	\leq Ca\|\eta'\|_1
	\operatorname*{ess\,sup}_{0<t<T}\|m_n(t)\|_1
	\leq C_Ta\|\eta'\|_1.
	\]
	Since $g=b_k(\rho_n)$, the first term on the right of
	\eqref{eq:periodic-work-identity} is $\mathcal R_{n,k,a}$.
	Combining this identity with \eqref{eq:work-small} proves
	\eqref{eq:normalized-flux-remainder}.
\end{proof}

The bound \eqref{eq:normalized-flux-remainder} permits the limit
$n\to\infty$ at fixed $k$ and smoothing scale $a$, followed by
$a\downarrow0$. The final limit $k\to\infty$ belongs to the entropy
argument. In particular, multiplication by $k$ after the flux limit
causes no loss: $k$ is fixed throughout the first two limits.

\begin{lemma}\label{prop:flux}
	After passing to a common subsequence, let $Y_{t,x}$ be the
	density Young measure, with $\rho=\av z$, and suppose that
	$u_n\rightharpoonup u$ in $L^2H^1$. Set $d=\Div u$ and
	\[
	\begin{gathered}
	\Pi=\av{z^\gamma},\qquad
	\beta_k=\av{T_k},\qquad
	\Gamma_k=\av{z^\gamma T_k},\qquad
	C_k=\Gamma_k-\Pi\beta_k,\\
	H_k=\operatorname{weak\,lim}_n T_k(\rho_n)d_n,\qquad
	h_k=\operatorname{weak\,lim}_n a_k(\rho_n)d_n,
	\end{gathered}
	\]
	where the last two limits are weak in $L^2(Q_T)$ for fixed $k$.
	Then $C_k\geq0$ almost everywhere, and
	\begin{equation}\label{eq:averaged-flux}
	\int_\T C_k(t,x)\dd x
	=\nu\int_\T(H_k-\beta_kd)(t,x)\dd x
	\quad\text{for almost every }t.
	\end{equation}
	Moreover,
	\[
	T_k(\rho_n)\to\beta_k
	\quad\text{in }C([0,T];H^{-1}(\T)),
	\]
	and
	\begin{equation}\label{eq:beta-equation}
	\partial_t\beta_k+\Div(\beta_ku)+h_k=0,
	\qquad \beta_k(0)=T_k(\rho_0).
	\end{equation}
\end{lemma}

\begin{proof} Extract a common subsequence on integer time intervals and
	countably many truncation levels. The uniform density bounds
	give a probability-valued Young measure $Y_{t,x}$, see \cite{Ball}.
Here and below $\av{F}=\int_{[0,\infty)}F(z)\,dY_{t,x}(z)$.
For bounded continuous $F$, the fundamental theorem identifies
$F(\rho_n)$ with this average in the weak limit. For unbounded $F$,
one must additionally control its tails; a Young measure alone does
not account for a possible concentration of $F(\rho_n)$.
Set
	\[
	B=\sup_n\esssup_{0<t<T}\int_\T\rho_n^\gamma\,dx.
	\]
	For a measurable $I\subset(0,T)$, apply Young-measure convergence
	to the bounded truncations
	$F_L(z)=\min\{z^\gamma,L\}$ and
	$G_R(z)=\min\{z,R\}$. Passing to the limits
	$L,R\to\infty$ gives
	\begin{equation}\label{eq:young-slice-moments}
	\esssup_{0<t<T}\int_\T\av{z^\gamma}\,dx\leq B,
	\qquad
	\int_\T\av z(t,x)\,dx=\mathfrak m
	\quad\text{for almost every }t.
	\end{equation}
	In particular, Jensen's inequality yields
	$\esssup_t\int_\T\rho^\gamma\leq B$.
	
	The pressure tail estimate allows us to pass from bounded
	truncations to $z^\gamma$ and $z^\gamma T_k(z)$, since for fixed $k$,
\[
 \sup_n\int_{\{\rho_n>R\}}
 \bigl(\rho_n^\gamma+\rho_n^\gamma T_k(\rho_n)\bigr)
 \leq(1+2k)C_TR^{-1/6}\longrightarrow0.
\]
Together with the $L^s$, $s>1$, pressure bound, this excludes a
concentration measure in both limits. There is no corresponding claim
for the untruncated product $\rho_n^{\gamma+1}$.
	Moreover, Lemma~\ref{lem:pressure-consequences} gives
	\[
	\|e_n\rho_n^3T_k(\rho_n)\|_{L^1(Q_T)}
	\leq2k\|e_n\rho_n^3\|_{L^1(Q_T)}
	\longrightarrow0.
	\]
	Consequently,
	\[
	P_n\rightharpoonup\Pi,\qquad
	P_nT_k(\rho_n)\rightharpoonup\Gamma_k
	\quad\text{weakly in }L^1(Q_T).
	\]
	Since $T_k$ and $a_k$ are bounded by $2k$,
	the products $T_k(\rho_n)d_n$ and $a_k(\rho_n)d_n$
	admit the weak $L^2$ limits $H_k$ and $h_k$. By symmetry,
	\[
	C_k=\frac12\iint
	(z^\gamma-w^\gamma)(T_k(z)-T_k(w))
	\dd Y(z)\dd Y(w)\geq0,
	\]
	where the inequality follows from the monotonicity of
	$z^\gamma$ and $T_k$.
	
	For fixed $k$, the renormalized continuity equation gives
	\[
	\partial_tT_k(\rho_n)
	+\Div(T_k(\rho_n)u_n)+a_k(\rho_n)d_n=0.
	\]
	Hence
	\[
	\begin{aligned}
	\|\partial_tT_k(\rho_n)\|_{L^2H^{-1}}
	&\leq2k\bigl(\|u_n\|_{L^2L^2}+\|d_n\|_2\bigr)
	\leq C_Tk,\\
	\|T_k(\rho_n)(t)-T_k(\rho_n)(s)\|_{H^{-1}}
	&\leq C_Tk|t-s|^{1/2}.
	\end{aligned}
	\]
	The compact embedding $L^2\Subset H^{-1}$ and
	the time equicontinuity yield compactness in
	$C([0,T];H^{-1})$, by the compactness argument recalled in
Lemma~\ref{lem:compactness}; compare \cite{Simon}. The initial data satisfy
	\[
	\|T_k(\rho_{0n})-T_k(\rho_0)\|_2^2
	\leq2k\|\rho_{0n}-\rho_0\|_1\longrightarrow0.
	\]
	Thus
	\[
	T_k(\rho_n)\longrightarrow\beta_k
	\quad\text{in }C([0,T];H^{-1}),
	\qquad
	\beta_k(0)=T_k(\rho_0).
	\]
	For every smooth vector field $\varphi$,
	\[
	\begin{aligned}
	&\left|\int_0^T
	\langle T_k(\rho_n)-\beta_k,
	\varphi\cdot u_n\rangle\dd t\right|\\
	&\qquad\leq C_\varphi
	\|T_k(\rho_n)-\beta_k\|_{C_tH^{-1}}
	\|u_n\|_{L^2H^1}\longrightarrow0.
	\end{aligned}
	\]
	Together with $u_n\rightharpoonup u$,
	this identifies the flux limit as $\beta_ku$.
	Passing to the limit in the initialized truncated equation
	proves \eqref{eq:beta-equation}.
	
	Fix $k$ and $a>0$. The preceding strong convergence gives
	\[
	\begin{aligned}
	&\|S_a[T_k(\rho_n)-\beta_k]\|_{C_{t,x}}\\
	&\qquad\leq
	\|\chi_a^{\rm per}\|_{H^1}
	\|T_k(\rho_n)-\beta_k\|_{C_tH^{-1}}
	\longrightarrow0.
	\end{aligned}
	\]
	Taking $n\to\infty$ in the remainder of
	Lemma~\ref{lem:uniform-high-frequency} yields
	\[
	\lim_n\mathcal R_{n,k,a}
	=\frac1k\int_{Q_T}\eta
	\bigl[\Gamma_k-\nu H_k
	-(\Pi-\nu d)S_a\beta_k\bigr].
	\]
	
	We next let $a\downarrow0$, keeping $k$ fixed.
	Since $S_a\beta_k\to\beta_k$ in $L^1(Q_T)$
	and $|S_a\beta_k-\beta_k|\leq4k$, for every
	$F\in L^1(Q_T)$ and $N>0$,
	\[
	\begin{aligned}
	\int_{Q_T}|F|\,|S_a\beta_k-\beta_k|
	&\leq4k\int_{\{|F|>N\}}|F|+N\|S_a\beta_k-\beta_k\|_1.
	\end{aligned}
	\]
	Sending first $a\downarrow0$ and then $N\to\infty$
	gives
	\[
	\int_{Q_T}F(S_a\beta_k-\beta_k)\longrightarrow0.
	\]
	
	By \eqref{eq:normalized-flux-remainder},
	the high-frequency remainder tends to zero uniformly
	in $n$ and $k$. Therefore, after multiplying
	the limiting identity by $k$,
	\[
	\int_0^T\eta(t)\int_\T
	\bigl[C_k-\nu(H_k-\beta_kd)\bigr]
	\,dx\,dt=0.
	\]
	As $\eta\in C_c^\infty(0,T)$ is arbitrary,
	this proves \eqref{eq:averaged-flux}. This is an identity of spatial
integrals for almost every time. No spatially pointwise equality
$C_k=\nu(H_k-\beta_kd)$ is inferred. The distinction is harmless for
the subsequent entropy argument because its transport divergence
also vanishes after integration over the torus. The covariance
$C_k\geq0$, in contrast, is a pointwise consequence of the
one-variable monotonicity calculation above.
\end{proof}

\begin{lemma}\label{prop:renorm}
	For every $T>0$, there exists $C_T$ independent of $k$ such that
	\begin{equation}\label{eq:oscillation}
	\lim_n\int_{Q_T}|T_k(\rho_n)-T_k(\rho)|^{\gamma+1}
	+\lim_n\int_{Q_T}|T_k(\rho_n)-\beta_k|^{\gamma+1}
	\leq C_T.
	\end{equation}
	Moreover, for $L\geq0$ and $k>L$,
	\begin{equation}\label{eq:source-decay}
	\int_{Q_T\cap\{\beta_k\leq L\}}|h_k|
	\leq C_Tk(k-L)^{-(\gamma+1)/2}.
	\end{equation}
	The mass flux converges to $\rho u$, including the initial terms.
	The limiting density satisfies the continuity and renormalized
	continuity equations with the prescribed initial data, for every
	$b\in C^1([0,\infty))$ with compactly supported derivative.
\end{lemma}

\begin{proof}
We use the nonnegative covariance $C_k$ and the averaged flux identity
\eqref{eq:averaged-flux}. In \eqref{eq:beta-equation}, $h_k$ is only the
weak limit of the approximate sources. In particular, neither
$h_k=a_k(\rho)d$ nor $\beta_k=T_k(\rho)$ is available at this stage.
The oscillation estimate first gives decay of $h_k$ on bounded
sublevel sets of $\beta_k$; this is exactly where a renormalization
with boundedly supported derivative will sample the source.
	
	\medskip
\noindent\emph{Step 1: Uniform oscillation control.} Let $Y=Y_{t,x}$ be the Young measure from Lemma~\ref{prop:flux}.
	For $a,b\geq0$, pressure monotonicity and the $1$-Lipschitz
	property of $T_k$ give
	\[
	(a^\gamma-b^\gamma)(T_k(a)-T_k(b))
	\geq |T_k(a)-T_k(b)|^{\gamma+1}.
	\]
	Since $T_k$ is concave, Jensen's inequality gives
	$\beta_k\leq T_k(\rho)$ and $\Pi\geq\rho^\gamma$.
	Taking $b=\rho$ and averaging with respect to $Y$, we obtain
	\[
	\begin{aligned}
	\av{|T_k(z)-T_k(\rho)|^{\gamma+1}}
	&\leq C_k+(\Pi-\rho^\gamma)
	(\beta_k-T_k(\rho))\\
	&\leq C_k.
	\end{aligned}
	\]
	Set
	\[
	X_k=\int_{Q_T}\av{|T_k(z)-T_k(\rho)|^{\gamma+1}},
	\qquad
	V_k=\int_{Q_T}\av{|T_k(z)-\beta_k|^2}.
	\]
	As $\beta_k=\av{T_k(z)}$, H\"older's inequality implies
	\[
	X_k\leq\int_{Q_T}C_k,
	\qquad
	V_k\leq |Q_T|^{1-2/(\gamma+1)}X_k^{2/(\gamma+1)}.
	\]
	Writing $D_T=\sup_n\|\Div u_n\|_{L^2(Q_T)}$, the averaged flux
	identity and Cauchy--Schwarz yield
	\[
	\begin{aligned}
	X_k
	&\leq \nu\lim_n\int_{Q_T}
	(T_k(\rho_n)-\beta_k)\Div u_n\\
	&\leq \nu D_T V_k^{1/2}\\
	&\leq \nu D_T|Q_T|^{1/2-1/(\gamma+1)}
	X_k^{1/(\gamma+1)}.
	\end{aligned}
	\]
	Thus, including the case $X_k=0$,
	\[
	X_k\leq
	\left(\nu D_T|Q_T|^{1/2-1/(\gamma+1)}\right)^{(\gamma+1)/\gamma}.
	\]
	Jensen's and the triangle inequalities also give
	\[
	\begin{aligned}
	|\beta_k-T_k(\rho)|^{\gamma+1}
	&\leq\av{|T_k(z)-T_k(\rho)|^{\gamma+1}},\\
	\int_{Q_T}\av{|T_k(z)-\beta_k|^{\gamma+1}}
	&\leq2^{\gamma+1}X_k.
	\end{aligned}
	\]
	For fixed $k$, the corresponding Young-measure limits hold with the
	measurable centers $T_k(\rho)$ and $\beta_k$: approximate these bounded
	centers by simple functions and use Lipschitz continuity of the powers
	on $[0,2k]$. An initial time layer of width $\varepsilon$ contributes
	at most $C(k)\varepsilon$. Sending $\varepsilon\downarrow0$ proves
	\eqref{eq:oscillation}, uniformly in $k$.
	
	\medskip
\noindent\emph{Step 2: Decay of the truncated source.} Fix $L\geq0$, $k>L$, and let
	$F=Q_T\cap\{\beta_k\leq L\}$. On
	$F\cap\{\rho_n>k\}$, the truncation properties imply
	\[
	|a_k(\rho_n)|\leq T_k(\rho_n)
	\leq\frac{k}{k-L}(T_k(\rho_n)-\beta_k),
	\qquad T_k(\rho_n)-\beta_k\geq k-L.
	\]
	Since $a_k=0$ on $[0,k]$ and $\gamma+1>2$, we have
	\[
	\ind_F|a_k(\rho_n)|^2
	\leq\frac{k^2}{(k-L)^{\gamma+1}}
	|T_k(\rho_n)-\beta_k|^{\gamma+1}.
	\]
	Recall that $a_k(\rho_n)\Div u_n\rightharpoonup h_k$ in $L^2$
	for fixed $k$. Testing this convergence with
	$\ind_F\operatorname{sign}h_k$, using Cauchy--Schwarz, and applying
	\eqref{eq:oscillation}, we obtain
	\[
	\begin{aligned}
	\int_F|h_k|
	&=\lim_n\int_F a_k(\rho_n)\Div u_n
	\operatorname{sign}h_k\\
	&\leq D_T\limsup_n\|\ind_Fa_k(\rho_n)\|_2\\
	&\leq C_Tk(k-L)^{-(\gamma+1)/2}.
	\end{aligned}
	\]
	This is \eqref{eq:source-decay}. In particular, for
	$k\geq2\max\{1,L\}$,
	\[
	\|\ind_{\{\beta_k\leq L\}}h_k\|_{L^1(Q_T)}
	\leq C_{T,L}k^{-(\gamma-1)/2}\longrightarrow0.
	\]
	Hence $b'(\beta_k)h_k\to0$ in $L^1(Q_T)$ whenever
	$b'=0$ on $[L,\infty)$.
	
	\medskip
\noindent\emph{Step 3: Mass flux and initialized renormalization.} First recover the kinetic bound for the limit. Put
	\[
	B=\sup_n\esssup_t\int_\T\rho_n^\gamma,
	\qquad
	K_{\mathrm{kin}}=\sup_n\esssup_t\int_\T\rho_n|u_n|^2.
	\]
	For each fixed $j$, Lemma~\ref{prop:flux} gives
	$T_j(\rho_n)u_n\rightharpoonup\beta_ju$, while
	\[
	\int_\T T_j(\rho_n)|u_n|^2\leq K_{\mathrm{kin}},
	\qquad
	\int_\T|T_j(\rho_n)u_n|^2\leq2jK_{\mathrm{kin}}.
	\]
	For completeness, weighted lower semicontinuity follows from the
convex perspective representation
\[
 \frac{|w|^2}{r}
 =\sup_{\xi\in\mathbb R^3}\{2w\cdot\xi-r|\xi|^2\},
 \qquad r>0,
\]
extended by $0$ at $(r,w)=(0,0)$ and by $+\infty$ at
$r=0$, $w\ne0$. Apply this formula with
$r=T_j(\rho_n)$ and $w=T_j(\rho_n)u_n$.
The joint weak limits are $\beta_j$ and $\beta_ju$; testing the
supremum with smooth vector fields and then taking the supremum
proves the lower semicontinuity. On any measurable time set
	$I\subset(0,T)$, it yields
	\[
	\int_I\!\int_\T\beta_j|u|^2
	\leq\liminf_n\int_I\!\int_\T T_j(\rho_n)|u_n|^2
	\leq K_{\mathrm{kin}}|I|.
	\]
	By concavity of the truncation,
	$T_j(z)=j\mathsf T(z/j)\uparrow z$ as $j\to\infty$;
	indeed,
	$\partial_j[j\mathsf T(z/j)]
	=\mathsf T(z/j)-(z/j)\mathsf T'(z/j)\geq0$.
	Thus $\beta_j\uparrow\rho$, and monotone convergence gives
	\begin{equation}\label{eq:preliminary-periodic-kinetic}
	\esssup_t\int_\T\rho|u|^2\leq K_{\mathrm{kin}}.
	\end{equation}
	
	The time-slice density moments also imply, for almost every $t$,
	\[
	\int_\T(\rho_n-T_j(\rho_n))\leq Bj^{1-\gamma},
	\qquad
	\int_\T(\rho-\beta_j)\leq Bj^{1-\gamma}.
	\]
	Weighted Cauchy--Schwarz now gives
	\[
	\begin{aligned}
	&\|(\rho_n-T_j(\rho_n))u_n\|_1
	+\|(\rho-\beta_j)u\|_1\\
	&\qquad\leq2\sqrt{K_{\mathrm{kin}}B}\,
	j^{-(\gamma-1)/2}.
	\end{aligned}
	\]
	Let $n\to\infty$ at fixed $j$ in the truncated mass flux, and then
let $j\to\infty$ in the preceding tail bound. For any smooth bounded
vector field $\varphi$ on $[0,T]\times\T$, this gives
\[
 \lim_n\int_{Q_T}m_n\cdot\varphi
 =\lim_{j\to\infty}\int_{Q_T}\beta_ju\cdot\varphi
 =\int_{Q_T}\rho u\cdot\varphi.
\]
Hence the weak momentum limit is $m=\rho u$. We now pass to the
limit in the original initialized, untruncated continuity equation,
using $\rho_{0n}\to\rho_0$ in $L^1$. This avoids any assertion of
global decay for $h_j$ and recovers
	\[
	\partial_t\rho+\Div(\rho u)=0,
	\qquad \rho(0)=\rho_0.
	\]
	Moreover, $\beta_j\uparrow\rho\in L^\gamma(Q_T)$ implies
	$\beta_j\to\rho$ in $L^r(Q_T)$ for $1\leq r\leq\gamma$.

Set $r_m=2\gamma/(\gamma+1)$. The density bound and
\eqref{eq:preliminary-periodic-kinetic} give
$m\in L^\infty(0,T;L^{r_m})$.
The initialized continuity equation yields
$\rho\in C([0,T];H^{-4})$, which, together with
the uniform $L^\gamma$ bound, implies
\[
\rho\in C_w([0,T];L^\gamma(\T)),
\qquad \rho(0)=\rho_0.
\]

For $b\in C^1([0,\infty))$ with compactly supported
derivative, Lemma~\ref{lem:forced-renormalization}
applied to \eqref{eq:beta-equation} gives the source version of
the DiPerna--Lions chain rule \cite{DL},

\[
\partial_tb(\beta_k)+\Div(b(\beta_k)u)
+(\beta_kb'(\beta_k)-b(\beta_k))\Div u
+b'(\beta_k)h_k=0.
\]
As $k\to\infty$, bounded convergence handles the
transport, divergence, and initial terms, while
\eqref{eq:source-decay} gives
$b'(\beta_k)h_k\to0$ in $L^1(Q_T)$.
Thus $\rho$ satisfies the renormalized continuity
equation with initial data $\rho_0$.
The case $b(0)\ne0$ follows from
$\Div(b(0)u)-b(0)\Div u=0$.
\end{proof}

\medskip\noindent\emph{Momentum compactness before the entropy comparison.}
The following estimates use the pressure-derived momentum bound and the
initialized momentum equation. They are available before strong density
convergence; the identification $m=\rho u$ was obtained by the truncation
argument in Lemma~\ref{prop:renorm}.

The approximate energy bound gives $m_n$ bounded in $L^\infty L^{r_m}$, where $r_m=2\gamma/(\gamma+1)>1$. Test the momentum equation, with its initial term, against $H^4$ vector fields. Since $H^4(\T)\hookrightarrow W^{1,\infty}$,
\begin{align*}
 \|m_n(t)-m_n(s)\|_{H^{-4}}
 &\leq C\int_s^t(\|\rho_n(r)u_n(r)\otimes u_n(r)\|_1+\|P_n(r)\|_1)\dd r\\
 &\quad+C\int_s^t\|S(u_n(r))\|_2\dd r\\
 &\leq C_T\bigl(|t-s|+|t-s|^{1/2}\bigr).
\end{align*}
This uses the $L^\infty_tL^1_x$ bounds for $(\rho_nu_n\otimes u_n),P_n$ and the $L^2(Q_T)$ stress bound. The initial term extends the estimate to $s=0$, with $m_n(0)=m_{0n}$. Lemma~\ref{lem:compactness} yields
\begin{equation}\label{cn:stb:momentum-C}
 m_n\to m=\rho u\quad\text{in }C([0,T];H^{-4}).
\end{equation}
The time compactness uses the initialized momentum equation and the
energy bounds, not an energy equality for the limit. The
time-uniform $L^{r_m}$ bound and density of smooth test functions also imply
\[
 m\in C_w([0,T];L^{r_m}),\qquad m(0)=m_0.
\]

For the chosen $p>6/5$, we have $p'<6$ and
\[
 H^1(\T)\Subset L^{p'}(\T),\qquad
 L^p(\T)\Subset H^{-1}(\T).
\]
The limit $m$ has the same $L^2L^p$ bound. The compact interpolation estimate in Lemma~\ref{lem:compactness} gives
\[
 \|m_n-m\|_{L^2H^{-1}}
 \leq\varepsilon\|m_n-m\|_{L^2L^p}
 +C_\varepsilon\sqrt T\,
   \|m_n-m\|_{C_tH^{-4}}.
\]
First let $n\to\infty$, then $\varepsilon\downarrow0$, to obtain
\begin{equation}\label{cn:stb:momentum-H}
 m_n\to\rho u\quad\text{in }L^2(0,T;H^{-1}).
\end{equation}
The inhomogeneous norms retain the full constant Fourier mode and allow arbitrary mean momentum.

\medskip\noindent\emph{Closing the density oscillations.}
We now compare the initialized entropy identities. The averaged flux
identity controls the spatial integral of the pressure covariance;
Lemma~\ref{prop:renorm} supplies the limiting renormalized equation
without identifying $h_k$ prematurely. These are distinct inputs:
the former supplies a nonnegative term, while the latter permits an
entropy comparison with the same initial value. The remaining error is
the truncation defect $Q_k$, which will tend to zero in $L^2$ before
strict convexity is invoked.

\begin{proof}[Proof of Proposition~\ref{prop:stability} (i)]
Following \cite[Section 8]{Feireisl}, choose truncated entropies with divergence coefficient $T_k$ and compare their spatial integrals through \eqref{eq:averaged-flux}. For $k\geq1$, set
\[
 \begin{aligned}
 L_k(z)&=z\int_1^zT_k(s)s^{-2}\dd s\quad(z>0),
 \qquad L_k(0)=0,\\
 \ell_k&=\av{L_k},\qquad Q_k=T_k(\rho)-\beta_k\geq0.
 \end{aligned}
\]
We have
\[
 zL_k'-L_k=T_k,\qquad L_k''=\frac{T_k'}z\geq0.
\]
Thus $L_k$ is convex, agrees with $z\log z$ on $[0,k]$, and is affine
on $[3k,\infty)$. The affine continuation preserves convexity, which
will give a nonnegative Jensen defect. A cutoff of the entropy itself
would not in general preserve that property.

The derivative of $z\log z$ is singular at zero. Lemma~\ref{lem:vacuum-entropy}
regularizes this derivative, subtracts the affine high-density part
to obtain an admissible renormalization, and restores that part by the
continuity equation. At fixed $k$, the vacuum regularization has errors
of order $\varepsilon$ in the entropy, the divergence coefficient, and
the initial integral. They vanish for both the approximate and the
limiting equations before the comparison below. No positive lower
bound on either density is required.

Integrating the resulting renormalized equations over the torus removes the entropy-flux divergence. Affine growth of $L_k$ and the $\gamma$-moment bound, $\gamma>1$, give
\[
 L_k(\rho_n)\rightharpoonup\ell_k\text{ in }L^1(Q_T),\qquad
 L_k(\rho_{0n})\to L_k(\rho_0)\text{ in }L^1(\T).
\]
For every $\psi\in C_c^\infty([0,T))$, we obtain
\begin{align}
 -\int_0^T\psi'\int_\T\ell_k
 +\int_0^T\psi\int_\T H_k
 &=\psi(0)\int_\T L_k(\rho_0),
 \label{eq:entropy-member-limit}\\
 -\int_0^T\psi'\int_\T L_k(\rho)
 +\int_0^T\psi\int_\T T_k(\rho)d
 &=\psi(0)\int_\T L_k(\rho_0).
 \label{eq:entropy-limit-companion}
\end{align}
Here $H_k$ is the weak limit of $T_k(\rho_n)d_n$, whereas $h_k$ is the
source in the truncated continuity equation. Their roles differ:
$zL_k'(z)-L_k(z)=T_k(z)$, but
$zT_k'(z)-T_k(z)=a_k(z)$. The identical right-hand sides in
\eqref{eq:entropy-member-limit}--\eqref{eq:entropy-limit-companion}
come from strong convergence of the prepared initial data. Equality
of the initial masses alone would not give equality of these entropy
integrals.

Set
\[
 \mathscr D_k(t)=\int_\T\bigl[\ell_k-L_k(\rho)\bigr]\geq0.
\]
Subtracting \eqref{eq:entropy-limit-companion} from \eqref{eq:entropy-member-limit} gives
\[
 \mathscr D_k'+\int_\T(H_k-T_k(\rho)d)=0,\qquad
 \mathscr D_k(0)=0.
\]
The derivative is integrable, so $\mathscr D_k$ has an absolutely continuous representative. Substituting \eqref{eq:averaged-flux} into the spatial integral yields
\begin{equation}\label{eq:integrated-defect}
 \mathscr D_k'+\frac1\nu\int_\T C_k
 =\int_\T Q_kd,
 \qquad\mathscr D_k(0)=0.
\end{equation}
The covariance term controls the sign after spatial integration, while
the remaining term contains the truncation defect $Q_k$, rather than
an unidentified product involving $h_k$. To estimate this right-hand
side, observe that
\[
 0\leq Q_k\leq\rho-\beta_k,
 \qquad\|Q_k\|_{L^1(Q_T)}\leq TBk^{1-\gamma}.
\]
Jensen's inequality and \eqref{eq:oscillation} also give
\[
 \int_{Q_T}|Q_k|^{\gamma+1}
 \leq\int_{Q_T}\av{|T_k(\rho)-T_k(z)|^{\gamma+1}}
 \leq C_T.
\]
Since
\[
 \frac12=\theta+\frac{1-\theta}{\gamma+1},
 \qquad \theta=\frac{\gamma-1}{2\gamma},
\]
interpolation yields
\begin{equation}\label{cn:stb:Q-small}
 \|Q_k\|_2
 \leq\|Q_k\|_1^{(\gamma-1)/(2\gamma)}
       \|Q_k\|_{\gamma+1}^{(\gamma+1)/(2\gamma)}
 \leq C_Tk^{-(\gamma-1)^2/(2\gamma)}\longrightarrow0.
\end{equation}

To recover the full entropy $z\log z$, set $h_{\rm ent}(z)=z\log z$, with value zero at $z=0$. The error $E_k=h_{\rm ent}-L_k$ vanishes on $[0,k]$ and is nonnegative and convex, since
\[
 E_k''=\frac{1-T_k'}z\geq0.
\]
For $z>k$,
\[
 0\leq E_k(z)\leq z\log(z/k)
 \leq\frac{k^{1-\gamma}}{\exp(1)(\gamma-1)}z^\gamma.
\]
Define the full entropy defect
\[
 \mathscr D(t)=\int_\T
 \bigl[\av{h_{\rm ent}}-h_{\rm ent}(\rho)\bigr]\geq0.
\]
For almost every time,
\[
 0\leq\mathscr D-\mathscr D_k
 =\int_\T\bigl[\av{E_k}-E_k(\rho)\bigr]
 \leq C_\gamma Bk^{1-\gamma}.
\]
Given $0\leq\zeta\in C_c^\infty(0,T)$, use the backward primitive $\psi(t)=\int_t^T\zeta(s)\dd s$ in \eqref{eq:integrated-defect}. The zero initial difference, $C_k\geq0$, and \eqref{cn:stb:Q-small} give
\[
 0\leq\int_0^T\zeta\mathscr D_k
 \leq\|\psi\|_\infty\|Q_k\|_2\|d\|_2\longrightarrow0.
\]
The test $\psi$ is admissible up to $t=0$ and vanishes near $T$.
Its initial contribution is zero because the two initialized entropy
identities have the same initial term. Moreover,
\[
 \left|\int_0^T\zeta(\mathscr D-\mathscr D_k)\right|
 \leq C_\gamma B\|\zeta\|_1k^{1-\gamma}\longrightarrow0.
\]
Letting $k\to\infty$ gives $\mathscr D=0$. The nonnegative Jensen defect therefore vanishes almost everywhere. Strict convexity gives $Y_{t,x}=\delta_\rho$ where $\rho>0$; nonnegative support and zero mean give $\delta_0$ where $\rho=0$. Hence
\begin{equation}\label{cn:stb:young-dirac}
 Y_{t,x}=\delta_{\rho(t,x)}\quad\text{almost everywhere in }Q_T.
\end{equation}

Set $J(z)=z/(1+z)$. Testing with $J,J^2$ and the bounded multiplier $J(\rho)$ gives
\[
 \|J(\rho_n)-J(\rho)\|_{L^2(Q_T)}^2
 =\int_{Q_T}\bigl[J(\rho_n)^2-2J(\rho_n)J(\rho)+J(\rho)^2\bigr]
 \longrightarrow0.
\]
The inverse of $J$ is Lipschitz on bounded density ranges; the tail bound controls their complements. Thus $\rho_n\to\rho$ in measure. Set $q_c=\gamma+1/6$. For every $1\leq r<q_c$, the distribution-function formula and \eqref{eq:moments} give
\[
 \int_{\{\rho_n>R\}}\rho_n^r
 \leq C_T\left(1+\frac r{q_c-r}\right)R^{r-q_c}.
\]
Thus $\rho_n^r$ is uniformly integrable, and Vitali's theorem yields
\[
 \rho_n\to\rho\quad\text{in }L^r(Q_T),
 \qquad1\leq r<\gamma+\frac16.
\]
The same uniform-integrability argument and the artificial-pressure estimate give
\[
 \begin{aligned}
 \rho_n^\gamma&\to\rho^\gamma
 &&\text{in }L^s(Q_T),\quad 1\leq s<1+\frac1{6\gamma},\\
 P_n&\to\rho^\gamma
 &&\text{in }L^s(Q_T),\quad
 1\leq s<1+\frac1{6\max\{3,\gamma\}}.
 \end{aligned}
\]
The continuity equation supplies the initial trace.
\end{proof}

The negative-norm momentum compactness was proved before the entropy
comparison. We now use it to identify convection, and use strong density
convergence to identify the weak limit of $\sqrt{\rho_n}u_n$. Convergence
of the kinetic norms then gives strong convergence of these weighted
velocities. Only after these steps do we pass to the energy inequality
using $E_{n0}\to E_0$.
\begin{proof}[Proof of Proposition~\ref{prop:stability} (ii)]
We use the momentum convergence already established in
\eqref{cn:stb:momentum-C}--\eqref{cn:stb:momentum-H}. The entropy argument
has now identified the density and pressure, so these compactness estimates
can be combined with the strong density convergence.

Let $\Psi$ be a smooth periodic matrix field compactly supported in time in $(0,T)$. Since multiplication by $\Psi$ preserves $H^1$, weak convergence of $\Psi u_n$ in $L^2H^1$ and \eqref{cn:stb:momentum-H} give
\[
 \int_{Q_T}(\rho_nu_n\otimes u_n):\Psi
 =\int_0^T\langle m_n,\Psi u_n\rangle\dd t
 \longrightarrow\int_{Q_T}(\rho u\otimes u):\Psi.
\]
Indeed, the pairing is split into
\[
 \int_0^T\langle m_n-m,\Psi u_n\rangle\,dt
 +\int_0^T\langle m,\Psi(u_n-u)\rangle\,dt,
\]
whose first term tends to zero by strong negative-norm convergence
and whose second tends to zero by weak $L^2H^1$ convergence.
Thus $\rho_nu_n\otimes u_n\to\rho u\otimes u$ in distributions.
This argument precedes the passage to the limiting energy inequality.

Set $W_n=\sqrt{\rho_n}u_n$ and $W=\sqrt\rho u$, with value zero on vacuum. Strong density convergence gives
\[
 \|\sqrt{\rho_n}-\sqrt\rho\|_{2\gamma}^{2\gamma}
 \leq\|\rho_n-\rho\|_\gamma^\gamma\longrightarrow0,
\]
and hence the square roots also converge strongly in $L^2_tL^{6/5}_x$. For a smooth vector test function $\varphi$,
\[
 \left|\int_{Q_T}(\sqrt{\rho_n}-\sqrt\rho)u_n\cdot\varphi\right|
 \leq\|\varphi\|_\infty
 \|\sqrt{\rho_n}-\sqrt\rho\|_{L^2L^{6/5}}
 \|u_n\|_{L^2L^6}\longrightarrow0.
\]
Weak convergence of $u_n$ handles the fixed coefficient $\sqrt\rho$. The kinetic-energy bound then gives $W_n\rightharpoonup W$ in $L^2(Q_T)$.

Take $\Psi=\eta(t)\Id$, $\eta\geq0$, in the identified convective limit to obtain
\[
 \int_{Q_T}\eta|W_n|^2\longrightarrow
 \int_{Q_T}\eta|W|^2,
 \qquad\eta\in C_c^\infty(0,T).
\]
Weak convergence and convergence of the weighted norms yield strong convergence on interior time intervals. For $0<\varepsilon<T/4$, choose $0\leq\eta_\varepsilon\leq1$ supported in $(0,T)$ and equal to one on $[\varepsilon,T-\varepsilon]$. Set
$K_{\mathrm{kin}}=\sup_n\esssup_t\int_\T\rho_n|u_n|^2$.
The approximate bound and \eqref{eq:preliminary-periodic-kinetic} give
\[
 \int_{Q_T}(1-\eta_\varepsilon)(|W_n|^2+|W|^2)
 \leq2K_{\mathrm{kin}}\,| (0,\varepsilon)\cup(T-\varepsilon,T) |
 =4K_{\mathrm{kin}}\varepsilon.
\]
Let $n\to\infty$ and then $\varepsilon\downarrow0$. The norms converge on $Q_T$, and hence
\begin{equation}\label{cn:stb:weighted-strong}
 \|W_n-W\|_2^2
 =\|W_n\|_2^2+\|W\|_2^2-2\int_{Q_T}W_n\cdot W
 \longrightarrow0.
\end{equation}
Consequently,
\begin{equation}\label{cn:stb:convection-strong}
 \|W_n\otimes W_n-W\otimes W\|_1
 \leq(\|W_n\|_2+\|W\|_2)\|W_n-W\|_2\longrightarrow0.
\end{equation}
Strong convergence of the pressure, convection, and initial data, together with weak convergence of the stress, yields the limiting momentum equation. Conservation of mass and momentum follows by spatial integration and extends to all times by weak continuity.

The same convergence gives strong $L^1$ convergence of the full energy density:
\begin{equation}\label{cn:stb:energy-density}
 \begin{gathered}
 \frac12|W_n|^2+\frac{\rho_n^\gamma}{\gamma-1}
       +\frac12e_n\rho_n^3
 \longrightarrow
 \frac12|W|^2+\frac{\rho^\gamma}{\gamma-1}\\
 \text{in }L^1(Q_T).
 \end{gathered}
\end{equation}
Integrating the energy density in space gives
\begin{equation}\label{eq:total-energy-conv}
 E_n\to E\qquad\text{in }L^1(0,T),
\end{equation}
where $E_n$ is defined in Lemma~\ref{prop:prepared-family}. Thus
$\int E_n\psi'\to\int E\psi'$ also for sign-changing $\psi'$. For every $0\leq\psi\in C_c^\infty([0,T))$, lower semicontinuity of the dissipation and $E_{n0}\to E_0$ yield
\[
 \begin{aligned}
 -\int_0^TE\psi'+\int_0^T\mathcal D\psi
 &\leq\liminf_n
 \left[-\int_0^TE_n\psi'+\int_0^T\mathcal D_n\psi\right]\\
 &\leq E_0\psi(0).
 \end{aligned}
\]
Here $\mathcal D(t)=\mu\|\nabla u(t)\|_2^2+(\lambda+\mu)\|\Div u(t)\|_2^2$, so the distributional energy inequality follows. The distributional energy inequality gives
\[
E(t)+\int_0^t\mathcal D\leq E_0
\quad\text{for almost every }t.
\]
For each fixed $t$, choose times $t_j\to t$ from the full-measure
set where the inequality holds. The weakly continuous representatives
of $\rho$ and $m$ give weak convergence at these times. The pressure
potential is convex, and the kinetic energy has the convex
perspective form $|m|^2/(2\rho)$ used above. Consequently,
\[
 E(t)\leq\liminf_{j\to\infty}E(t_j),\qquad
 \int_0^{t_j}\mathcal D(s)\,ds\longrightarrow
 \int_0^t\mathcal D(s)\,ds.
\]
This extends the inequality to every $t\geq0$, with the energy
defined by these representatives. On vacuum the kinetic density is
zero; finiteness of the perspective functional enforces $m=0$ there.
\end{proof}

\begin{proof}[Proof of Theorem~\ref{thm:main}]
For zero initial mass, the zero solution suffices.
Otherwise, Lemma~\ref{prop:prepared-family} provides
approximate solutions satisfying \eqref{eq:budget}
with $E_{n0}\to E_0$.
Proposition~\ref{thm:cap} and
Lemma~\ref{lem:pressure-consequences}
give the uniform estimates required by
Proposition~\ref{prop:stability}, which yields
the desired weak solution on each finite time interval.

A diagonal subsequence over integer time intervals,
truncation levels, and test functions gives a global
solution on $[0,\infty)$.
\end{proof}

\appendix
\section{Cubic Artificial-Pressure Approximation}
\label{app:construction}

We prove Lemma~\ref{prop:prepared-family} by preparing
the initial data and applying the fixed-$e$ existence
result below. The Galerkin scheme, artificial density diffusion, and
auxiliary pressure follow the artificial-pressure approximation framework
of \cite{FeireislBook,FNP,Lions}. We give the periodic implementation
because the later material test requires the fixed-member bounds
$\rho\in L^\infty L^3$ and $m\in L^2L^2$.
The exponent three is used for these pairings; no assertion of priority
for artificial cubic pressure is needed.

The successive limits are
\[
 N\longrightarrow\infty\quad(\alpha,e,\sigma\text{ fixed}),
 \qquad \alpha\downarrow0\quad(e,\sigma\text{ fixed}),
 \qquad \sigma\downarrow0\quad(e\text{ fixed}).
\]
Only after these limits does Section~\ref{sec:scalar} send $e=e_n$ to
zero. In particular, the $e$-dependent bounds proved here justify tests
on an individual approximate solution; they do not enter the uniform
pressure budget of Proposition~\ref{thm:cap}.

This construction does not use the barycentric pressure estimate.
At fixed $\sigma>0$, the auxiliary exponent $G$ provides the
integrability for the untruncated effective-flux test. At fixed $e>0$,
the cubic term provides the corresponding integrability while
$\sigma$ is removed. Thus the pressure gains
\eqref{eq:classical-pressure-gain} and \eqref{eq:sigma-gain} are
classical divergence-test estimates with constants allowed to depend
on the parameters that remain fixed. The later uniform estimate
\eqref{eq:tail} is used only after this supply of approximate solutions
has been constructed.

\subsection{Preparation of Initial Data}

Let $\mathfrak m=\int_\T\rho_0>0$,
$A_0=\sqrt{\rho_0}$, and
$W_0=m_0/\sqrt{\rho_0}$, extended by zero on vacuum.
The factorization $m_0=A_0W_0$ separates the two quantities controlled by
finite energy. Smoothing $A_0$ and $W_0$ separately preserves the kinetic
energy through an unweighted $L^2$ norm, including on the vacuum set.
Truncating $A_0$ also makes the initial artificial energy vanish.
For any prescribed $e_n>0$ with $e_n\to0$, set
\[
L_n=
\begin{cases}
\max\{1,e_n^{-1/[4(3-\gamma)]}\},&\gamma<3,\\
\infty,&\gamma\geq3.
\end{cases}
\]
With a nonnegative periodic mollifier $\omega_{1/n}$,
define
\begin{equation}\label{eq:periodic-preparation}
\begin{gathered}
\widetilde a_n
=\omega_{1/n}*\min(A_0,L_n)+n^{-1},
\qquad
a_n=\frac{\sqrt{\mathfrak m}}{\|\widetilde a_n\|_2}
\widetilde a_n,\\
b_n=\omega_{1/n}*W_0,\qquad
\rho_{0n}=a_n^2,\qquad
m_{0n}=a_nb_n,\qquad
u_{0n}=b_n/a_n.
\end{gathered}
\end{equation}

Then $\rho_{0n}>0$ is smooth,
$\int_\T\rho_{0n}=\mathfrak m$, and
\[
a_n\to A_0\quad\text{in }L^2\cap L^{2\gamma},
\qquad
b_n\to W_0\quad\text{in }L^2.
\]
Consequently,
$\rho_{0n}\to\rho_0$ in $L^\gamma$ and
$m_{0n}\to m_0$ in $L^{r_m}$, where
$r_m=2\gamma/(\gamma+1)$. The initial kinetic energy converges since
\[
\frac12\int_\T\frac{|m_{0n}|^2}{\rho_{0n}}
=\frac12\|b_n\|_2^2\longrightarrow
\frac12\|W_0\|_2^2.
\]
For $\gamma<3$,
\[
e_n\int_\T a_n^6
\leq Ce_nL_n^{6-2\gamma}
\int_\T a_n^{2\gamma}
\leq Ce_n^{1/2}\longrightarrow0.
\]
For $\gamma\geq3$, the same term is $O(e_n)$.
The normalization factor in \eqref{eq:periodic-preparation} tends to one
because $\|A_0\|_2^2=\mathfrak m$. The estimates above use only an upper
bound for this factor, while the lower bound $a_n>0$ is allowed to depend
on $n$. Moreover, H\"older's inequality gives the explicit momentum bound
\[
\|a_nb_n-A_0W_0\|_{r_m}
\leq\|a_n-A_0\|_{2\gamma}\|b_n\|_2
 +\|A_0\|_{2\gamma}\|b_n-W_0\|_2\longrightarrow0.
\]
Thus $E_{n0}\to E_0$, without a positive lower bound for the limiting
density.

\begin{lemma}\label{prop:supply}
	Let $e>0$, $\rho_*>0$ and $u_*$ be smooth
	periodic initial data. Set $m_*=\rho_*u_*$,
	$b_*=\max\{3,\gamma\}$, and
	\[
	P_e(\rho)=\rho^\gamma+e\rho^3,\qquad
	H_e(\rho)=\frac{\rho^\gamma}{\gamma-1}
	+\frac e2\rho^3.
	\]
	There exists a global periodic finite-energy weak
	solution with initial data $(\rho_*,m_*)$, satisfying
	the renormalized continuity equation with its initial
	term for every $C^1$ function with compactly
	supported derivative. For every $T>0$,
	
	\[
	\left\{
	\begin{aligned}
	&\rho\in L^\infty L^{b_*}
	\cap L^{b_*+1}(Q_T)\cap C_wL^{b_*},\\
	&u\in L^2H^1,\qquad
	\sqrt{\rho}\,u\in L^\infty L^2,\\
	&\rho u\in C_wL^{2b_*/(b_*+1)}
	\cap L^2L^2,\\
	&\rho u\otimes u\in L^1L^{3/2}.
	\end{aligned}
	\right.
	\]
	
	Total mass and momentum are conserved. Moreover,
	\begin{equation}\label{eq:supply-energy}
	-\int_0^T\mathcal E_e\psi'\,dt
	+\int_0^T\mathcal D_u\psi\,dt
	\leq\mathcal E_{e,0}\psi(0),
	\qquad 0\leq\psi\in C_c^\infty([0,T)),
	\end{equation}
	where
	\[
	\begin{aligned}
	\mathcal E_e(t)
	&=\int_\T\left(
	\frac12\rho|u|^2+H_e(\rho)\right)dx,\\
	\mathcal D_u(t)
	&=\mu\|\nabla u(t)\|_2^2
	+(\lambda+\mu)\|\Div u(t)\|_2^2,\\
	\mathcal E_{e,0}
	&=\int_\T\left(
	\frac{|m_*|^2}{2\rho_*}+H_e(\rho_*)\right)dx.
	\end{aligned}
	\]
\end{lemma}

\begin{proof}
	Fix $e>0$ and smooth positive initial data.
	
	\medskip\noindent\emph{Step 1: Galerkin approximation and the limit $N\to\infty$.}
	Fix $G>\max\{8,\gamma+1\}$, $\sigma>0$, and $0<\alpha\leq1$, and set
	\[
	\widehat P(\rho)=\rho^\gamma+e\rho^3+\sigma \rho^G,\qquad
	H(\rho)=H_e(\rho)+\frac{\sigma}{G-1}\rho^G.
	\]
	Let $X_N$ be the space of real vector Fourier modes with $|k|\leq N$,
	including the constant modes. For an $L^2$-orthonormal basis $(\varphi_a)$,
	define
	\begin{equation}\label{eq:supply-mass-matrix}
	\mathcal M[\rho]_{ab}=\int_\T \rho\varphi_a\cdot\varphi_b,\qquad
	(\inf \rho)|c|^2\leq c^{\mathsf T}\mathcal M[\rho]c
	\leq(\sup \rho)|c|^2.
	\end{equation}
	For positive $\rho$, the matrix is invertible. Choose
	$u_N^0=\operatorname{Proj}_{X_N}^{L^2(\rho_*\dd x)}u_*$, so that
	$\int_\T\rho_*u_N^0\cdot\varphi_a=\int_\T m_*\cdot\varphi_a$.
	In particular, the constant modes preserve the initial total momentum.
	The weighted orthogonality also gives
	\[
	\int_\T\rho_*|u_*-u_N^0|^2
	=\int_\T\rho_*|u_*|^2-\int_\T\rho_*|u_N^0|^2
	\longrightarrow0,
	\]
	since the Fourier spaces are dense and $\rho_*$ is bounded above and
	below at this fixed stage. This both bounds the initial Galerkin
	energy and identifies its limit.
	The density diffusion and its momentum correction have the form used
	in \cite[p.~361]{FNP}. The Galerkin system is
	\begin{equation}\label{eq:inner-system}
	\left\{
	\begin{aligned}
	&\partial_t\rho_N+\Div(\rho_Nu_N)
	=\alpha\Delta\rho_N,\\
	&\partial_t(\rho_Nu_N)
	+\Div(\rho_Nu_N\otimes u_N)
	+\nabla\widehat P(\rho_N)=\Div S(u_N)-D_N,
	\end{aligned}
	\right.
	\end{equation}
	where
	\[
	D_N=\alpha\sum_{j=1}^3
	(\partial_j\rho_N)\partial_j u_N,
	\qquad
	(\rho_N,u_N)|_{t=0}=(\rho_*,u_N^0).
	\]
	The momentum equation is imposed in $X_N$, while the density equation is
	unprojected. For a prescribed Galerkin velocity, the maximum principle gives
	\[
	(\min \rho_*)e^{-\int_0^t\|\Div u_N\|_\infty}
	\leq \rho_N(t,x)\leq
	(\max \rho_*)e^{\int_0^t\|\Div u_N\|_\infty}.
	\]
	For velocities in a bounded subset of $X_N$, the difference of the density
	equations satisfies
	$\|\rho[u]-\rho[w]\|_{C_tL^2}
	\leq C_{N,\alpha,\rho_*}\sqrt{\tau}\|u-w\|_{C_tX_N}$.
	The inverse mass matrix is locally Lipschitz, and the Galerkin coefficients
	satisfy
	\[
	u(t)=\mathcal M[\rho[u](t)]^{-1}
	\left(m_*|_{X_N}+\int_0^t\mathcal F[\rho[u],u](s)\dd s\right).
	\]
	For $w\in X_N$, integration by parts gives
	\[
	-\alpha\int_\T\partial_j\rho_N\,\partial_ju_{N,i}w_i
	=\alpha\int_\T\rho_N(\Delta u_{N,i}w_i
	+\partial_ju_{N,i}\partial_jw_i).
	\]
	The resulting fixed-point map has Lipschitz constant
	$C(\sqrt{\tau}+\tau)$ on a short interval and yields a local solution
	as in \cite{FNP}.
	
	Set $K_N=\frac12\int_\T\rho_N|u_N|^2$. Testing the Galerkin momentum
	equation with $u_N$ gives
	\[
	\int_\T u_N\cdot[\partial_t(\rho_Nu_N)+\Div(\rho_Nu_N\otimes u_N)]
	=K_N'-\int_\T D_N\cdot u_N.
	\]
	Indeed, the extra term from the density equation is
	\[
	\frac12\int_\T|u_N|^2\alpha\Delta\rho_N
	=-\alpha\int_\T\partial_j\rho_N\,
	                  u_N\cdot\partial_ju_N
	=-\int_\T D_N\cdot u_N.
	\]
	The correction cancels the corresponding term in the momentum equation.
	For the internal energy, the periodic chain rule gives
	\[
	\frac{\dd}{\dd t}\int_\T H(\rho_N)
	+\int_\T\widehat P(\rho_N)\Div u_N
	=-\alpha\int_\T H''(\rho_N)|\nabla\rho_N|^2.
	\]
	Thus both artificial terms have a specified role: diffusion regularizes
	the density, and its momentum correction restores the kinetic balance.
	Using $\rho H'(\rho)-H(\rho)=\widehat P(\rho)$ in the density equation,
	we obtain
	\begin{equation}\label{eq:galerkin-energy}
	\begin{aligned}
	&\frac{\dd}{\dd t}\int_\T
	\left(\frac12\rho_N|u_N|^2+H(\rho_N)\right)
	+\mu\|\nabla u_N\|_2^2
	+(\lambda+\mu)\|\Div u_N\|_2^2\\
	&\qquad+\alpha\int_\T H''(\rho_N)|\nabla \rho_N|^2=0.
	\end{aligned}
	\end{equation}
	Here $H''(s)=\gamma s^{\gamma-2}+3es+\sigma Gs^{G-2}>0$, and
	\[
	\int_\T S(u_N):\nabla u_N
	=\mu\|\nabla u_N\|_2^2
	+(\lambda+\mu)\|\Div u_N\|_2^2.
	\]
	Since $\lambda+\mu\geq\mu/3>0$, the energy identity and
	Lemma~\ref{lem:mean} control $u_N$ in $L^2(0,T;H^1)$.
	For fixed $N$, finite-dimensional norm equivalence bounds
	$\int_0^T\|\Div u_N\|_\infty$; positivity of $\rho_N$ then permits
	global continuation.
	
	Now let $N\to\infty$ with $\alpha,e,\sigma,G$ fixed.
	The energy identity and the equation for $\rho_N^2$ give
	\begin{gather}
	\rho_N\in L^\infty L^G,\quad \rho_N^{G/2}\in L^2H^1,\quad
	\rho_N\in L^{5G/3}(Q_T),\quad
	\widehat P(\rho_N)\in L^{5/3}(Q_T),\label{eq:inner-gain}\\
	\frac12\frac{\dd}{\dd t}\|\rho_N\|_2^2
	+\alpha\|\nabla \rho_N\|_2^2
	=-\frac12\int_\T \rho_N^2\Div u_N.\label{eq:inner-square}
	\end{gather}
	Interpolation yields the space--time density bound, while the second
	identity gives a bound for $\sqrt\alpha\nabla\rho_N$ uniform in $N$
	and $0<\alpha\leq1$. At fixed $\alpha>0$,
	$m_N=\rho_Nu_N\in L^2L^2$, and
	\[
	\|\partial_t \rho_N\|_{L^2H^{-1}}
	\leq C\bigl(\|m_N\|_{L^2L^2}
	+\alpha\|\nabla \rho_N\|_{L^2L^2}\bigr)\leq C_\alpha.
	\]
	For the Fourier projection $\mathbb P_L$,
	\[
	\|(I-\mathbb P_L)\rho_N\|_{L^2(Q_T)}^2
	\leq\lambda_{L+1}^{-1}\|\nabla \rho_N\|_{L^2(Q_T)}^2.
	\]
	The time-compactness argument is a Fourier realization of the compactness
	criterion in \cite{Simon}; the compact spatial embeddings used here are
	those of \cite{Brezis}. The time derivative controls each finite set of
	modes, and the displayed estimate removes the remaining modes at fixed
	$\alpha$. Finite-mode compactness and this tail bound yield, along a subsequence,
	\[
	\rho_N\to \rho\quad\text{in }L^q(Q_T),\quad 1\leq q<5G/3,
	\qquad u_N\rightharpoonup u\quad\text{in }L^2H^1.
	\]
	Consequently, $\rho_Nu_N\to\rho u$ in distributions and
	$\widehat P(\rho_N)\to\widehat P(\rho)$,
	$H(\rho_N)\to H(\rho)$ in $L^1(Q_T)$.
	The limiting density belongs to $L^2H^1$ with
	$\partial_t\rho\in L^2H^{-1}$, hence
	$\rho\in C([0,T];L^2)$ and $\rho(0)=\rho_*$.
	
	To control the momentum correction, the assumption $G>8$ yields
	\begin{equation}\label{eq:supply-time-interpolation}
	\|m_N\|_{L^{26/3}_tL^2_x}
	\leq\|m_N\|_{L^\infty_tL^{16/9}_x}^{10/13}
	\|m_N\|_{L^2_tL^{24/7}_x}^{3/13},\qquad m_N=\rho_Nu_N.
	\end{equation}
	On the zero-mean space, set $\mathcal L=-\Delta$ and
	$Z_N=\mathcal L^{-1/2}(\rho_N-\mathfrak m)$.
	The continuity equation becomes
	$Z_N'+\alpha\mathcal LZ_N=\mathcal Nm_N$, with
	$\mathcal N:L^2(\T;\R^3)\to L^2_0(\T)$ bounded.
	Maximal regularity \cite{deSimon} gives
	\[
	\|\nabla \rho_N\|_{L^{26/3}_tL^2_x}
	\leq T^{3/26}\|\nabla \rho_*\|_2+C_\alpha\|m_N\|_{L^{26/3}_tL^2_x},
	\qquad
	\|D_N\|_{L^{13/8}_tL^1_x}\leq C_\alpha.
	\]
	Thus the correction has integral $O(|I|^{5/13})$ over a time interval $I$.
	
	Multiplying \eqref{eq:inner-square} by $T-t$ gives
	\begin{equation}\label{eq:galerkin-weighted-density-gradient}
	\begin{aligned}
	\alpha\int_0^T(T-t)\|\nabla \rho_N\|_2^2\dd t
	={}&\frac T2\|\rho_*\|_2^2-\frac12\int_0^T\|\rho_N\|_2^2\dd t\\
	&-\frac12\int_{Q_T}(T-t)\rho_N^2\Div u_N.
	\end{aligned}
	\end{equation}
	The limiting continuity equation has the same identity.
	Since $\rho_N\to\rho$ in $L^4(Q_T)$ and
	$\Div u_N\rightharpoonup\Div u$ in $L^2$, comparison of the two
	identities yields
	\[
	\sqrt{T-t}\,\nabla \rho_N\longrightarrow\sqrt{T-t}\,\nabla \rho
	\quad\text{in }L^2(Q_T).
	\]
	Moreover, maximal regularity gives
	\[
	\sup_N\int_{T-h}^T\|\nabla \rho_N\|_2^2\dd t
	+\int_{T-h}^T\|\nabla \rho\|_2^2\dd t
	\leq C_\alpha h^{10/13}.
	\]
	For $0<h<T$, split the norm at $T-h$:
	\[
	\begin{aligned}
	\|\nabla\rho_N-\nabla\rho\|_{L^2(Q_T)}^2
	&\leq h^{-1}\int_0^T(T-t)
	       \|\nabla\rho_N-\nabla\rho\|_2^2\dd t
	       +C_\alpha h^{10/13}.
	\end{aligned}
	\]
	First let $N\to\infty$, then $h\downarrow0$. This proves
	$\nabla\rho_N\to\nabla\rho$ in $L^2(Q_T)$ at fixed $\alpha$.
	The terminal-layer estimate is essential: the factor $T-t$ alone
	would only give strong convergence away from $T$.
	
	Testing the momentum equation with a fixed Fourier mode
	$\varphi_k\in X_N$ yields, for $0\leq s<t\leq T$,
	\begin{equation}\label{eq:galerkin-momentum-mode-modulus}
	\left|\int_\T[m_N(t)-m_N(s)]\cdot\varphi_k\right|
	\leq C_{k,\alpha}\bigl(|t-s|+|t-s|^{1/2}+|t-s|^{5/13}\bigr).
	\end{equation}
	The Fourier tail estimate
	\begin{equation}\label{eq:periodic-spectral-tail}
	\|(I-\mathbb P_L)f\|_{H^{-1}}
	\leq C\lambda_{L+1}^{-13/32}\|f\|_{16/9}
	\end{equation}
	follows by interpolation. Together with Arzel\`a--Ascoli, it gives
	$m_N\to\rho u$ in $C([0,T];H^{-1})$, with initial value $m_*$.
	For a smooth matrix field $\Psi$,
	\begin{equation}\label{eq:construction-convection-pairing}
	\begin{split}
	\int_{Q_T}(\rho_Nu_N\otimes u_N-\rho u\otimes u):\Psi
	={}&\int_0^T\langle m_N-\rho u,\Psi u_N\rangle\dd t\\
	&+\int_0^T\langle \rho u,\Psi(u_N-u)\rangle\dd t\longrightarrow0.
	\end{split}
	\end{equation}
	The $H^{-1}$--$H^1$ pairing is valid for smooth $\Psi$.
	Identifying the weak limit of $\sqrt{\rho_N}u_N$ and taking
	$\Psi=\Id$ gives convergence of its $L^2$ norm. Hence
	\begin{equation}\label{eq:construction-strong-products}
	\sqrt{\rho_N}u_N\to\sqrt \rho\,u\quad\hbox{in }L^2(Q_T),\qquad
	\rho_Nu_N\otimes u_N\to \rho u\otimes u\quad\hbox{in }L^1(Q_T).
	\end{equation}
	Approximation of smooth tests by Fourier polynomials recovers the full
	momentum equation and its initial term. The strong convergence of
	$\nabla\rho_N$ and weak convergence of $\nabla u_N$ identify the
	correction as $D_i=\alpha\partial_j\rho\partial_ju_i$.
	
	For $b\in C^2([0,\infty))$ with compactly supported $b'$, the diffusive
	chain rule passes to the limit:
	\begin{equation}\label{eq:periodic-diffusive-renormalization}
	\partial_tb(\rho)+\Div(b(\rho)u)+(\rho b'(\rho)-b(\rho))\Div u
	=\alpha\Div(b'(\rho)\nabla \rho)-\alpha b''(\rho)|\nabla \rho|^2,
	\end{equation}
	with initial value $b(\rho_*)$. The case $b(0)\ne0$ follows by adding
	the constant transport identity. Finally, convergence of the initial
	kinetic energy and of the energy densities, followed by lower
	semicontinuity of dissipation in \eqref{eq:galerkin-energy}, gives the
	energy inequality with the prescribed initial energy.

\medskip
\noindent\emph{Step 2: Density diffusion
	$\alpha\downarrow0$.} Fix $e,\sigma,G$ and let $\alpha\downarrow0$.
The energy bounds give
\[
\rho_\alpha\in L^\infty L^8,\qquad
m_\alpha\in L^\infty L^{16/9}\cap L^2L^2,
\qquad
\|D_\alpha\|_1+\|\alpha\nabla\rho_\alpha\|_2
\leq C\sqrt\alpha.
\]
For example,
\[
\|D_\alpha\|_1
\leq\sqrt\alpha\,
 \|\sqrt\alpha\nabla\rho_\alpha\|_2
 \|\nabla u_\alpha\|_2,
\qquad
\|\alpha\nabla\rho_\alpha\|_2
=\sqrt\alpha\|\sqrt\alpha\nabla\rho_\alpha\|_2.
\]
The square-density identity \eqref{eq:inner-square} bounds
$\sqrt\alpha\nabla\rho_\alpha$ using
$\rho_\alpha\in L^\infty L^8$ and $\Div u_\alpha\in L^2$.
Thus the artificial diffusion and momentum correction vanish.
These estimates do not require a uniform bound for
$\nabla\rho_\alpha$ itself.
Set $q_\alpha=A\rho_\alpha$ and
$d_\alpha=\Div u_\alpha$.
The effective-flux test follows \cite{FeireislBook,FNP,Lions}.
Here $A=\nabla\Delta^{-1}$ is the periodic inverse with zero Fourier
coefficient at the origin, so
\[
\Div A\rho_\alpha=\rho_\alpha-\mathfrak m.
\]
This identity produces the spatial mean term on the right below.
Testing the momentum equation with the periodic potential
$\eta q_\alpha$, where
$\eta\in C_c^\infty(0,T)$, yields
\begin{equation}\label{eq:inner-flux}
\begin{split}
\int\eta(\widehat P(\rho_\alpha)-\nu d_\alpha)
\rho_\alpha
={}&\mathfrak m\int\eta\widehat P(\rho_\alpha)
-\int\eta'm_\alpha\cdot q_\alpha\\
&+\int\eta u_\alpha\cdot B_\alpha
+\mathcal E_\alpha,\\
B_{\alpha,i}
={}&\rho_\alpha R_{ij}m_{\alpha,j}
-m_{\alpha,j}R_{ij}\rho_\alpha,\\
\mathcal E_\alpha
={}&\int\eta D_\alpha\cdot q_\alpha
-\alpha\int\eta m_\alpha\cdot\nabla\rho_\alpha.
\end{split}
\end{equation}
All integrals without specified domains are over $Q_T$.
At fixed $\alpha$, the integrability in
\eqref{eq:inner-gain} justifies the test by smooth
approximation. The energy bounds and
$A:L^8\to W^{1,8}\hookrightarrow L^\infty$
control the time and convective terms, while
the stress term is bounded in $L^1$ and
$\mathcal E_\alpha=O(\sqrt\alpha)$.
Removing the time cutoff gives
\begin{equation}\label{eq:classical-pressure-gain}
\int_{Q_T}\widehat P(\rho_\alpha)\rho_\alpha
\leq C,\qquad
\rho_\alpha\in L^{G+1}(Q_T),\qquad
\widehat P(\rho_\alpha)
\in L^{(G+1)/G}(Q_T),
\end{equation}
uniformly in $\alpha$.

The continuity equation bounds
$\partial_t\rho_\alpha$ in $L^2H^{-1}$.
Together with the uniform density bounds and
Fourier compactness, this gives
$\rho_\alpha\to\rho$ in $C([0,T];H^{-1})$.
For each fixed Fourier mode $\varphi_k$,
\begin{equation}\label{eq:diffusion-momentum-mode-modulus}
\begin{aligned}
\sup_{0\leq t\leq T}
\left|\int_0^t\int_\T
D_\alpha\cdot\varphi_k\right|
&\leq C(k)\sqrt\alpha,\\
\left|\int_\T
[m_\alpha(t)-m_\alpha(s)]\cdot\varphi_k\right|
&\leq C(k)\left(
|t-s|+|t-s|^{1/2}+\sqrt\alpha\right).
\end{aligned}
\end{equation}
Combining these bounds with
\eqref{eq:periodic-spectral-tail}
and the uniform $L^\infty L^{16/9}$ estimate
yields
\[
m_\alpha\to m
\quad\text{in }C([0,T];H^{-1}),
\qquad m(0)=m_*.
\]
Since $u_\alpha\rightharpoonup u$ in $L^2H^1$,
the density convergence gives, for every
smooth vector field $\phi$,
\begin{equation}\label{eq:construction-prepressure-identification}
\begin{aligned}
\int_{Q_T}(\rho_\alpha u_\alpha-\rho u)\cdot\phi
={}&\int_0^T
\langle\rho_\alpha-\rho,u_\alpha\cdot\phi\rangle\,dt\\
&+\int_0^T
\langle\rho,(u_\alpha-u)\cdot\phi\rangle\,dt
\longrightarrow0.
\end{aligned}
\end{equation}
Hence $m=\rho u$. The same
$H^{-1}$--$H^1$ pairing identifies the
convective term as in
\eqref{eq:construction-convection-pairing}.
Passing to the limit in the initialized equations
also recovers conservation of total momentum.

We next pass to the effective-flux identity.
Lemma~\ref{lem:construction-commutator-compactness},
applied with $f_j=\rho_\alpha$,
$g_j=m_\alpha$, $a=8$, and $b=16/9$,
gives
\[
B_\alpha\to B
\quad\text{strongly in }
L^2(0,T;H^{-1}(\T)),
\]
since $1/8+9/16<5/6$.
This convergence permits passage to the
term involving $u_\alpha\cdot B_\alpha$.

Write
\[
z=\operatorname{wlim}_{\alpha\downarrow0}
\rho_\alpha d_\alpha,\qquad
\Pi=\operatorname{w^*\!lim}_{\alpha\downarrow0}
\widehat P(\rho_\alpha)\rho_\alpha,
\]
where $\Pi$ is a nonnegative measure.
Let $Y_{t,x}$ be the density Young measure, obtained from the
fundamental theorem for Young measures (see \cite{Ball}). The product
$\widehat P(\rho_\alpha)\rho_\alpha$ is only bounded in $L^1$;
its measure limit may contain a concentration part. Consequently the
first relation below is an inequality, even though the separate pressure
limit is represented by $Y_{t,x}$.
By uniform integrability,
\[
\rho=\av{s},\qquad
\overline{\widehat P}=\av{\widehat P(s)}.
\]
Pressure monotonicity gives, as measures,
\[
\begin{aligned}
\Pi&\geq\av{\widehat P(s)s}\,dt\,dx,\\
\av{\widehat P(s)s}
-\av{\widehat P(s)}\av{s}
&=\frac12\iint
[\widehat P(s)-\widehat P(t)](s-t)
\,dY(s)\,dY(t)\geq0.
\end{aligned}
\]
Passing to the limit in \eqref{eq:inner-flux}
and subtracting the corresponding identity
for the limiting momentum equation yields
\[
\int\eta\,d\Pi
-\int\eta\overline{\widehat P}\,\rho\,dt\,dx
=\nu\int\eta(z-\rho d)\,dt\,dx.
\]
The left-hand side is nonnegative for
$\eta\geq0$. Since $\nu>0$,
\begin{equation}\label{eq:inner-integrated-sign}
\int_\T(z-\rho d)\,dx\geq0
\quad\text{for almost every }t.
\end{equation}

We now compare the entropy equations.
Regularize $s\log s$ by
\[
b_\varepsilon(s)
=(s+\varepsilon)\log(s+\varepsilon)
-\varepsilon\log\varepsilon
\]
and truncate its second derivative at high densities. More explicitly,
choose $0\leq\theta_R\leq1$ equal to one on $[0,R]$ and zero on
$[2R,\infty)$, and replace
$b_\varepsilon''(s)=(s+\varepsilon)^{-1}$ by
$\theta_R(s)/(s+\varepsilon)$. Twice integration, with the same value
and derivative at zero, gives a convex function affine for $s\geq2R$.
Subtracting this affine part makes its derivative compactly supported.
Convexity preserves the sign of the diffusion term
$-\alpha b''|\nabla\rho_\alpha|^2$.
At this stage $\rho_\alpha\in L^\infty L^8$ and
$\Div u_\alpha\in L^2$, so the high-density errors in the entropy and
source vanish as $R\to\infty$; then $\varepsilon\downarrow0$.
The initial data are fixed throughout these two limits.
Applying \eqref{eq:periodic-diffusive-renormalization},
restoring the affine part through the continuity
equation, and removing the truncations gives
\begin{equation}\label{eq:construction-entropy-member}
-\int\psi'\int_\T\rho_\alpha\log\rho_\alpha
+\int\psi\rho_\alpha d_\alpha
\leq\psi(0)\int_\T\rho_*\log\rho_*,
\qquad \psi\geq0.
\end{equation}
For the limiting continuity equation,
the initialized renormalized identity follows
from the spatial commutator argument
\cite{DL}, after extension to negative times.
Passing to the weak entropy limits gives
\begin{equation}\label{eq:construction-common-initial-entropy}
\begin{aligned}
-\int_0^T\psi'\int_\T
\overline{\rho\log\rho}
+\int_{Q_T}\psi z
&\leq\psi(0)\int_\T\rho_*\log\rho_*,\\
-\int_0^T\psi'\int_\T\rho\log\rho
+\int_{Q_T}\psi\rho d
&=\psi(0)\int_\T\rho_*\log\rho_*.
\end{aligned}
\end{equation}

Set
\[
\mathscr D(t)=
\int_\T(\overline{\rho\log\rho}
-\rho\log\rho)\,dx\geq0.
\]
Subtracting the two entropy relations and
taking
$\psi(t)=\int_t^T\zeta(s)\,ds$ with
$0\leq\zeta\in C_c^\infty(0,T)$ gives
\[
0\geq
\int_0^T\zeta(t)\mathscr D(t)\,dt
+\int_{Q_T}\psi(z-\rho d)\,dx\,dt
\geq0,
\]
where the last inequality follows from
\eqref{eq:inner-integrated-sign}.
Thus $\mathscr D=0$ almost everywhere.
Strict convexity of $s\log s$ eliminates
density oscillations, and the uniform
$L^{G+1}$ bound yields
\[
\rho_\alpha\to\rho
\quad\text{in }L^q(Q_T),\qquad
1\leq q<G+1,
\]
and
\[
\widehat P(\rho_\alpha)\to\widehat P(\rho)
\quad\text{in }L^s(Q_T),\qquad
1\leq s<\frac{G+1}{G}.
\]

Finally, the momentum pairing
\eqref{eq:construction-convection-pairing}
gives strong convergence of the weighted
velocities in $L^2(Q_T)$.
Convergence of the kinetic and internal energies,
together with lower semicontinuity of the
viscous dissipation, yields the initialized
energy inequality.

\medskip\noindent\emph{Step 3: Auxiliary pressure $\sigma\downarrow0$.} Fix $e,G$ and the smooth initial data, and let $\sigma_j\downarrow0$.
The energy inequality gives, uniformly in $j$,
\[
\rho_j\in L^\infty L^{b_*},\qquad
m_j\in L^\infty L^{3/2}\cap L^2L^2,\qquad
u_j\in L^2H^1,\qquad
\rho_ju_j\otimes u_j\in L^1L^{3/2}.
\]
Testing the momentum equation with $\eta A\rho_j$,
using $A:L^3\to W^{1,3}\hookrightarrow L^3$ and the
uniform energy bounds, gives
\begin{equation}\label{eq:sigma-gain}
\int_{Q_T}(\rho_j^{\gamma+1}+e\rho_j^4
+\sigma_j\rho_j^{G+1})\leq C,
\qquad
\|\sigma_j\rho_j^G\|_{(G+1)/G}
\leq C\sigma_j^{1/(G+1)}\longrightarrow0.
\end{equation}
The artificial-pressure norm in \eqref{eq:sigma-gain} follows directly
from
\[
\int_{Q_T}|\sigma_j\rho_j^G|^{(G+1)/G}
=\sigma_j^{1/G}\int_{Q_T}\sigma_j\rho_j^{G+1}
\leq C\sigma_j^{1/G}.
\]
Thus the pressure itself vanishes strongly, rather than merely acquiring
a small coefficient in front of a possibly unbounded density power.
In particular, $\rho_j$ is bounded in $L^{b_*+1}(Q_T)$ and
$P_e(\rho_j)$ in $L^{(b_*+1)/b_*}(Q_T)$.

The continuity equation bounds $\partial_t\rho_j$ in $L^2H^{-1}$.
The Fourier-mode estimates for momentum, now without the diffusion
correction, and
\[
\|(I-\mathbb P_L)f\|_{H^{-1}}
\leq C\lambda_{L+1}^{-1/4}\|f\|_{3/2}
\]
yield, after extraction,
$\rho_j\to\rho$ and $m_j\to m$ in $C([0,T];H^{-1})$,
with initial values $\rho_*,m_*$. The pairings in
\eqref{eq:construction-prepressure-identification} and
\eqref{eq:construction-convection-pairing} give $m=\rho u$ and
identify the convective limit. Denote the weak pressure limit by
$\overline P_e$.

For fixed $k$, set $g_j=T_k(\rho_j)$,
$h_j=a_k(\rho_j)d_j$, and $q_j=Ag_j$. The initialized
truncated equation satisfies
\begin{equation}\label{eq:construction-truncation-time-control}
\left\{
\begin{aligned}
\partial_tg_j+\Div(g_ju_j)+h_j&=0,\\
g_j(0)&=T_k(\rho_*),
\end{aligned}
\right.
\qquad
\|\partial_tg_j\|_{L^2H^{-1}}+\|h_j\|_2\leq C(k).
\end{equation}
Compactness gives $g_j\to\beta_k$ in $C([0,T];H^{-1})$.
Taking $h_j\rightharpoonup h_k$ in $L^2$ and passing to the limit,
we obtain
\[
\partial_t\beta_k+\Div(\beta_ku)+h_k=0,
\qquad \beta_k(0)=T_k(\rho_*).
\]
Since $A,A^*:H^{-1}\to L^2$ are bounded,
$q_j\to A\beta_k$ and $A^*m_j\to A^*m$ in $C_tL^2_x$.
The uniform $L^\infty W^{1,p}$ bounds for $q_j$, with $p>3$,
also give uniform convergence of $q_j$. Hence
\[
\int_{Q_T}\eta m_j\cdot Ah_j
=\int_{Q_T}\eta(A^*m_j)h_j
\longrightarrow\int_{Q_T}\eta(A^*m)h_k.
\]

For $B_{j,l}=g_jR_{il}m_{j,i}-m_{j,i}R_{il}g_j$,
Lemma~\ref{lem:construction-commutator-compactness}, with
$a=12$ and $b=3/2$, gives
\[
B_j\longrightarrow B\quad\text{in }L^2_tH^{-1}_x,
\]
since $1/12+2/3<5/6$. Testing the momentum equation with
$\eta q_j$ yields
\begin{equation}\label{eq:construction-truncated-flux}
\int\eta(\widehat P_j-\nu d_j)(g_j-\bar g_j)
=-\int\eta'm_j\cdot q_j
+\int\eta u_{j,l}B_{j,l}+\int\eta m_j\cdot Ah_j.
\end{equation}
Regularization justifies this test at fixed $j,k$. The limiting
test $\eta A\beta_k$ is admissible since its time derivative lies
in $L^2L^2$ and its gradient in $L^\infty L^p$ for finite $p$.
Passing to the limit and subtracting the two identities cancels the
right-hand sides; the mean terms cancel because
$\bar g_j\to\bar\beta_k$ uniformly in time. Moreover,
$\|\sigma_j\rho_j^GT_k(\rho_j)\|_1
\leq2k\|\sigma_j\rho_j^G\|_1\to0$. Thus
\begin{equation}\label{eq:construction-sigma-flux-limit}
\lim_j\int\eta[P_e(\rho_j)-\nu d_j]T_k(\rho_j)
=\int\eta[\overline P_e-\nu d]\beta_k.
\end{equation}

Let $Y$ be the density Young measure and write
$\Gamma_k=\av{P_eT_k}$ and $\beta_k=\av{T_k}$.
Monotonicity of $P_e$ and $T_k$ gives
\[
\Gamma_k-\overline P_e\beta_k
=\frac12\iint[P_e(s)-P_e(t)][T_k(s)-T_k(t)]
\dd Y(s)\dd Y(t)\geq0.
\]
By \eqref{eq:construction-sigma-flux-limit}, its spatial integral is
$\nu\int(z_k-\beta_kd)$, where
$z_k=\operatorname{wlim}_jT_k(\rho_j)d_j$.
The $L^4$ bound in \eqref{eq:sigma-gain} gives
\[
\|\rho_j-T_k(\rho_j)\|_2\leq C/k,\qquad
\|(\rho_j-T_k(\rho_j))d_j\|_1\leq C/k.
\]
Writing $z=\operatorname{wlim}_j\rho_jd_j$ in $L^{4/3}(Q_T)$,
we also have
\[
\|\rho-\beta_k\|_2\leq C/k,
\qquad \|z-z_k\|_1+\|(\beta_k-\rho)d\|_1\leq C/k.
\]
First letting $j\to\infty$ at fixed $k$ and then $k\to\infty$
recovers \eqref{eq:inner-integrated-sign}. The initialized
renormalized entropy equations for $\rho_j$ and $\rho$ (using
the negative-time extension for the limit), as in
\eqref{eq:construction-common-initial-entropy}, have the same
initial term. Subtracting them and testing in time with
$\psi(t)=\int_t^T\zeta(s)\dd s$, $\zeta\geq0$, the integrated sign
forces the entropy defect to vanish. Strict convexity yields
\[
\rho_j\longrightarrow\rho\quad\text{in }L^q(Q_T),
\qquad 1\leq q<b_*+1,
\]
and
\[
P_e(\rho_j)\longrightarrow P_e(\rho)
\quad\text{in }L^s(Q_T),
\qquad 1\leq s<\frac{b_*+1}{b_*}.
\]
The convective and weighted-velocity limits follow from
\eqref{eq:construction-strong-products}; the auxiliary internal
energy vanishes by \eqref{eq:sigma-gain}.

It remains to pass to the energy inequality and recover the initial
traces. At each of the three approximation limits, the weak limit of
$W_j=\sqrt{\rho_j}u_j$ in $L^2(Q_T)$ is
$W=\sqrt\rho\,u$. Strong momentum convergence in $C_tH^{-1}_x$
and weak convergence of $u_j$ in $L^2H^1$ give
\begin{equation}\label{eq:construction-kinetic-energy-recovery}
\begin{aligned}
\|W_j\|_{L^2(Q_T)}^2
&=\int_0^T\langle m_j,u_j\rangle\dd t
\longrightarrow\int_0^T\langle\rho u,u\rangle\dd t
=\|W\|_{L^2(Q_T)}^2,\\
W_j&\longrightarrow W\quad\text{in }L^2(Q_T),\\
|W_j|^2&\longrightarrow |W|^2\quad\text{in }L^1(Q_T).
\end{aligned}
\end{equation}
The internal energies converge in $L^1(Q_T)$ by
\[
\begin{array}{ll}
N\to\infty:& \rho_N\to \rho\text{ in }L^q,
\quad G<q<5G/3,\\[1mm]
\alpha\downarrow0:& \rho_\alpha\to\rho\text{ in }L^q,
\quad G<q<G+1,\\[1mm]
\sigma_j\downarrow0:& \rho_j\to\rho\text{ in }L^q,
\quad b_*<q<b_*+1,
\quad \sigma_j\rho_j^G\to0\text{ in }L^1.
\end{array}
\]
The initial kinetic energy converges under the weighted Galerkin
projection, the initial data remain fixed as $\alpha\downarrow0$, and
\[
\mathcal E_{\sigma,0}
=\mathcal E_{e,0}+\frac{\sigma}{G-1}\int_\T\rho_*^G
\longrightarrow\mathcal E_{e,0}.
\]
Thus, at each stage, $\mathcal E_j\to\mathcal E$ in $L^1(0,T)$
and $\mathcal E_{j,0}\to\mathcal E_0$. For
$0\leq\psi\in C_c^\infty([0,T))$,
\begin{equation}\label{eq:construction-initialized-energy-limit}
-\int\mathcal E\psi'+\int\mathcal D_u\psi
\leq\liminf_j\left[-\int\mathcal E_j\psi'
+\int\mathcal D_{u_j}\psi\right]
\leq\mathcal E_0\psi(0).
\end{equation}
The uniform energy and higher-integrability bounds control the endpoint
time layers. The initialized equations and negative-norm time
compactness give $\rho\in C_wL^{b_*}$ and
$\rho u\in C_wL^{2b_*/(b_*+1)}$, with initial values $\rho_*,m_*$;
the commutator argument provides the initial terms in the renormalized
equation. Finally, a diagonal subsequence on integer time intervals
gives the required global solution on $[0,\infty)$.
\end{proof}

\subsection{Auxiliary Periodic Commutator Lemma}
The following lemma was used in Steps 2 and 3 of the proof of Lemma~\ref{prop:supply}.
Its hypotheses are supplied by the negative-norm compactness of density
and momentum at those fixed-pressure stages. It passes a bilinear Riesz
commutator against the weak velocity by first obtaining a distributional
limit at each time in a common full-measure set, then using spatial
compactness and dominated convergence in time. The quantitative
high-frequency lemma used in Section~\ref{sec:scalar} has a different
role: it gives a remainder uniform in the final approximation index,
where these fixed-pressure spatial bounds are unavailable.

\begin{lemma}\label{lem:construction-commutator-compactness}
Let $1<a,b<\infty$ satisfy $1/a+1/b<5/6$. Suppose that
\[
 \begin{gathered}
 f_j\to f,\qquad g_j\to g
 \quad\text{strongly in }C([0,T];H^{-1}(\T)),\\
 \sup_j\bigl(\|f_j\|_{L^\infty_tL^a_x}
                 +\|g_j\|_{L^\infty_tL^b_x}\bigr)<\infty,
 \end{gathered}
\]
where $g_j=(g_{j,i})$ is vector-valued. For $R_{il}=\partial_i\partial_l\Delta^{-1}$ with its periodic zero-mode convention, set
\[
 B_{j,l}=f_jR_{il}g_{j,i}-g_{j,i}R_{il}f_j.
\]
Then $B_j\to B=fR_{il}g_i-g_iR_{il}f$ strongly in $L^2(0,T;H^{-1}(\T))$, with the limit initially understood in distributions.
\end{lemma}

\begin{proof}
Put $1/r=1/a+1/b<5/6$, so $r>6/5$ and $a>b'$, where $1/b+1/b'=1$. Boundedness of the periodic Riesz operators, obtained from the
Calder\'on--Zygmund estimates in \cite{Stein}, gives
\[
 \|B_j(t)\|_r\leq C\|f_j(t)\|_a\|g_j(t)\|_b,
 \qquad \sup_j\|B_j\|_{L^\infty_tL^r_x}<\infty.
\]
Choose the continuous $H^{-1}$ representatives. Remove the union of the null sets where either spatial bound fails, over the countable sequence, and denote the remaining full-measure set by $E$. At each time in $E$, weak compactness in $L^a$ and $L^b$, together with convergence in $H^{-1}$, identifies every cluster point. Thus
\[
 f_j(t)\rightharpoonup f(t)\text{ in }L^a,
 \qquad g_j(t)\rightharpoonup g(t)\text{ in }L^b,
 \qquad t\in E.
\]
The limits satisfy the same spatial bounds on $E$, which is independent of the spatial test function.

Fix $t\in E$ and $\phi\in C^\infty(\T)$. Self-adjointness of $R_{il}$ gives
\begin{equation}\label{eq:construction-commutator-pairing}
 \langle B_{j,l}(t),\phi\rangle
 =\int_\T g_{j,i}(t)\,[R_{il},\phi]f_j(t)\,dx.
\end{equation}
The contact term in the distribution kernel of $R_{il}$ cancels, because multiplication by $\phi$ commutes with that term. Near the diagonal the remaining kernel satisfies
\[
 |K_{il}(x-y)(\phi(y)-\phi(x))|
 \leq C_\phi|x-y|^{-2}.
\]
Decompose $[R_{il},\phi]=\mathcal C_h^{\rm near}+\mathcal C_h^{\rm far}$ with a smooth cutoff supported in $|x-y|<2h$. Young's inequality and $L^a(\T)\hookrightarrow L^{b'}(\T)$ give
\begin{equation}\label{eq:construction-commutator-near}
 \begin{aligned}
 \|\mathcal C_h^{\rm near}f\|_{b'}
 &\leq C\|\mathcal C_h^{\rm near}f\|_a
 \leq C_\phi h\|f\|_a,\\
 \left|\int_\T g\,\mathcal C_h^{\rm near}f\right|
 &\leq C_\phi h\|g\|_b\|f\|_a.
 \end{aligned}
\end{equation}
For fixed $h>0$, the far kernel is smooth on $\T\times\T$. It maps bounded subsets of $L^a$ to bounded subsets of $C^1(\T)$, and is therefore compact into $L^{b'}(\T)$. Thus
\[
 \mathcal C_h^{\rm far}f_j(t)\longrightarrow
 \mathcal C_h^{\rm far}f(t)\quad\text{in }L^{b'}.
\]
Weak convergence of $g_j(t)$ in $L^b$ then gives convergence of the far pairing. Apply \eqref{eq:construction-commutator-near} to the sequence and its limit and let $h\downarrow0$. This identifies the limit of \eqref{eq:construction-commutator-pairing} with the pairing for $f,g$, so $B_j(t)\to B(t)$ in distributions for every $t\in E$.

Since $r>6/5$, its conjugate exponent satisfies $r'<6$.
The compact Sobolev embedding $H^1(\T)\Subset L^{r'}(\T)$
\cite{Brezis}, followed by duality, gives
$L^r(\T)\Subset H^{-1}(\T)$. This strict inequality explains the
restriction $1/a+1/b<5/6$ in the lemma. The compact embedding and the
identified distributional limit imply strong $H^{-1}$ convergence at each such time: every subsequence has a strongly convergent further subsequence, and its limit must be $B(t)$. Finally,
\[
 \|B_j(t)-B(t)\|_{H^{-1}}^2\leq C\quad(t\in E),
 \qquad
 \int_0^T\|B_j(t)-B(t)\|_{H^{-1}}^2\dd t\longrightarrow0
\]
by dominated convergence.
\end{proof}

\begin{proof}[Proof of Lemma~\ref{prop:prepared-family}]
Apply the fixed-$e$ existence result to the data in \eqref{eq:periodic-preparation}, whose full initial energies converge to $E_0$. Coercivity of the periodic viscous term and Lemma~\ref{lem:mean} give the uniform bounds. For each approximate solution, the cubic density norm may grow as $e_n^{-1/3}$; it justifies the special test field, whereas the pressure bound is uniform in $n$.
\end{proof}

\section{Hilbert-Space Construction}
\label{app:barycentric-foundations}

We establish the Hilbert-space monotonicity and exact-marginal
representation used in Lemmas~\ref{lem:periodic-optimizer} and
\ref{lem:periodic-material}. Transport planning allows a prescribed
target law to be retained under relaxation \cite{ABC,GJ,GRST}; the
strong-convexity and subgradient arguments use the Hilbert-space
framework of \cite{BC}. The periodic lifts and the exact representation
needed by the restricted Jacobian argument are proved below.
The deterministic infimum first identifies the transport value.
Its convex envelope yields a unique residual optimizer, and an optimal
plan then realizes its subgradient as a conditional barycenter.
Section~\ref{sec:barycentric-spatial} converts that barycenter into a
periodic convex gradient.

\subsection{Hilbert-Space Envelope}\label{app:hilbert-envelope}\label{dp:app-hilbert}
Set $D=\operatorname{diam}(\T)=\sqrt3/2$. Let $(\mathcal A,\ell)$ be a standard atomless measure space of mass $c$, with positions $X(a)$, and set $\mathscr H=L^2(\mathcal A,\ell;\R^3)$ with inner product $\langle X,Y\rangle=\int_{\mathcal A}X\cdot Y\dd\ell$. Write $q:\R^3\to\T$ for the quotient map and define
\[
\mathcal K_{\rm law}
=\{Z\in\mathscr H:q_\#(Z_\#\ell)=c\,dx\}.
\]
For equal-mass measures, let $W_{2,\T}^2$ denote the unnormalized quadratic transport cost. Define
\begin{equation}\label{eq:periodic-hilbert-envelope}
\begin{aligned}
F(X)&=-\frac12\inf_{Z\in\mathcal K_{\rm law}}
\|X-Z\|_{\mathscr H}^2
=-\frac12W_{2,\T}^2(q_\#X_\#\ell,c\,dx),\\
F_K(X)&=\inf_{Y\in\mathscr H}
\left\{F(Y)+\frac K2\|X-Y\|_{\mathscr H}^2\right\},\\
\mathcal H_K(X)&=\frac12\|X\|_{\mathscr H}^2
+\frac{K-1}{K}F_K(X),
\qquad \mathscr T_K(X)=\nabla\mathcal H_K(X).
\end{aligned}
\end{equation}
Every deterministic target map induces an admissible plan, giving one inequality. Conversely, partition the source and target tori into cells of diameter at most $h$. Atomlessness allows the source labels to be split according to the mass assigned to each cell pair and mapped uniformly onto the corresponding target cell, preserving the marginal $c\,dx$. The costs converge as $h\downarrow0$, and nearest integer lifts give the Euclidean infimum. The infimum need not be attained by a deterministic map.

\begin{lemma}[Hilbert-space envelope]
	\label{lem:periodic-envelope}
	For $K>1$ and $X\in\mathscr H$, the minimizer $Y_K(X)$ in the definition of $F_K(X)$ is unique, with
	\[
	\mathscr T_K(X)=KX-(K-1)Y_K(X).
	\]
	The map $\mathscr T_K$ is $K$-Lipschitz and satisfies
	\begin{equation}\label{eq:periodic-firm}
	\langle\mathscr T_K(X)-\mathscr T_K(X'),X-X'\rangle
	\ge\frac1K\|\mathscr T_K(X)-\mathscr T_K(X')\|^2.
	\end{equation}
	The displacement $\mathscr T_K(X)-X$ is uniformly bounded and depends only on the source law and the current torus position. For every $H^1$ curve $X(t)$ in $\mathscr H$,
	\begin{equation}\label{eq:periodic-path}
	\langle\partial_t\mathscr T_K(X(t)),\dot X(t)\rangle
	\ge\frac1K\|\partial_t\mathscr T_K(X(t))\|^2
	\quad\text{for almost every }t.
	\end{equation}
\end{lemma}

\begin{proof}
	Expanding the square gives
	\[
	G(Y):=F(Y)+\tfrac12\|Y\|^2
	=\sup_{Z\in\mathcal K_{\rm law}}
	\{\langle Y,Z\rangle-\tfrac12\|Z\|^2\}.
	\]
	Thus $G$ is convex and weakly lower semicontinuous. Since $-cD^2/2\le F\le0$, the objective is coercive and equals
	\[
	G(Y)+\frac{K-1}{2}\|Y\|^2-K\langle X,Y\rangle
	+\frac K2\|X\|^2.
	\]
	A minimizing sequence is bounded, since its objective is at least
	$K\|X-Y\|^2/2-cD^2/2$. A weakly convergent subsequence exists in
	$\mathscr H$, and weak lower semicontinuity of the displayed expression
	gives a minimizer. Convexity of $G$ and $K-1>0$ imply uniqueness.
	More precisely, at the minimizer the objective satisfies
	\[
	\begin{aligned}
	F(Y)+\frac K2\|X-Y\|^2
	&\geq F(Y_K(X))+\frac K2\|X-Y_K(X)\|^2\\
	&\qquad+\frac{K-1}{2}\|Y-Y_K(X)\|^2.
	\end{aligned}
	\]
	The Hilbert-space subgradient rule \cite{BC} gives the optimality condition
	\begin{equation}\label{dp:subgradient}
	KX-(K-1)Y_K(X)\in\partial G(Y_K(X)).
	\end{equation}
	For two inputs, write $\Delta X=X-X'$, $\Delta Y=Y_K(X)-Y_K(X')$, and $\Delta\mathscr T=\mathscr T_K(X)-\mathscr T_K(X')$. Monotonicity of $\partial G$ gives
	\begin{align*}
	\langle\Delta\mathscr T,\Delta Y\rangle&\geq0,\\
	K\langle\Delta\mathscr T,\Delta X\rangle
	&=\|\Delta\mathscr T\|^2
	+(K-1)\langle\Delta\mathscr T,\Delta Y\rangle
	\geq\|\Delta\mathscr T\|^2.
	\end{align*}
	This gives \eqref{eq:periodic-firm} and, by Cauchy--Schwarz, the
	$K$-Lipschitz bound. Monotonicity also gives
	\[
	(K-1)\|\Delta Y\|^2
	\leq K\langle\Delta X,\Delta Y\rangle,
	\qquad
	\|\Delta Y\|\leq\frac{K}{K-1}\|\Delta X\|.
	\]
	For an increment $U$, use $Y_K(X)$ as competitor at $X+U$ and
	$Y_K(X+U)$ as competitor at $X$. This gives
	\[
	\begin{aligned}
	K\langle X-Y_K(X+U),U\rangle+\tfrac K2\|U\|^2
	&\leq F_K(X+U)-F_K(X),\\
	F_K(X+U)-F_K(X)
	&\leq K\langle X-Y_K(X),U\rangle+\tfrac K2\|U\|^2.
	\end{aligned}
	\]
	The Lipschitz bound for $Y_K$ makes the difference between the two
	linear terms $O(\|U\|^2)$. Thus differentiating the minimum value is
	justified and yields
	\[
	\nabla F_K(X)=K(X-Y_K(X)),\qquad
	\nabla\mathcal H_K(X)=KX-(K-1)Y_K(X),
	\]
	as required.
	
	The quadratic torus cost is $D$-Lipschitz in its first variable, so
	\[
	|F(X+U)-F(X)|\le D\|U\|_{L^1(\ell)}.
	\]
	Simultaneously translating the input and competitor shows that $F_K$ has the same bound. Localized variations then give $|K(X-Y_K(X))|\le D$ almost everywhere, and hence
	\begin{equation}\label{dp:bounded-displacement}
	|Y_K(X)-X|\le D/K,\qquad
	|\mathscr T_K(X)-X|\le (K-1)D/K\le D.
	\end{equation}
	
	For integer-valued $N\in\mathscr H$, uniqueness and translation invariance give
	\[
	F(X+N)=F(X),\qquad
	Y_K(X+N)=Y_K(X)+N,\qquad
	\mathscr T_K(X+N)=\mathscr T_K(X)+N.
	\]
	The displacement is therefore invariant under integer lifts and
measure-preserving relabelings. To prove that it depends only on position,
we use conditional Jensen's inequality for the convex law-invariant
functional $G$. Here is a direct justification. For $Y\in\mathscr H$,
realize conditionally independent copies $Y_1,Y_2,\ldots$ of $Y$ given
$X$ on a product label space. All $Y_j$ have the law of $Y$, and
\[
\begin{aligned}
G\left(\frac1N\sum_{j=1}^NY_j\right)
&\leq\frac1N\sum_{j=1}^NG(Y_j)=G(Y),\\
\left\|\frac1N\sum_{j=1}^NY_j-\mathbb E[Y\mid X]\right\|_2^2
&=\frac1N\|Y-\mathbb E[Y\mid X]\|_2^2\longrightarrow0.
\end{aligned}
\]
The same formula defining $G$ applies on every standard atomless label
space and depends only on the law of its argument. Since $F$ is
Lipschitz in $L^2$ and the squared norm is continuous, $G$ is continuous
in $L^2$. Passing to the limit proves
$G(\mathbb E[Y\mid X])\leq G(Y)$. This argument is an evaluation of the
law-invariant functional; it does not assume that minimizers persist
when the label space is enlarged.

Set $\bar Y=\mathbb E[Y\mid X]$. Orthogonality of conditional expectation
and $X$-measurability of $X$ give
\[
\begin{aligned}
\|Y\|^2&=\|\bar Y\|^2+\|Y-\bar Y\|^2,\\
\langle X,Y\rangle&=\langle X,\bar Y\rangle.
\end{aligned}
\]
In the objective written in terms of $G$, replacing $Y$ by $\bar Y$
therefore lowers the value by at least
$(K-1)\|Y-\bar Y\|^2/2$. At the minimizer this quantity must vanish,
so $Y_K(X)$ is a Borel function of $X$. Once competitors have this
form, their costs depend only on the source law. Uniqueness then gives
the same Borel function for any realization of that law, almost
everywhere for the law. Integer-lift invariance gives the asserted
dependence on the current torus position.
	
	Lipschitz continuity gives $\mathscr T_K(X(t))\in H^1$ whenever $X(t)\in H^1$. Applying \eqref{eq:periodic-firm} at neighboring times and passing to the difference-quotient limit proves \eqref{eq:periodic-path}.
\end{proof}

\subsection{Exact-Marginal Representation}
\label{app:exact-marginal-representation}

\begin{lemma}[Exact-marginal representation]\label{lem:exact-marginal-barycenter}
	Let $y\in\mathscr H$ be bounded and let $G$ be as in Lemma~\ref{lem:periodic-envelope}. Denote by $\mathfrak P_y$ the optimal plans for $\int|y(a)-z|^2\dd\pi(a,z)$, with source marginal $\ell$, target marginal $c\,dx$ modulo integers, and nearest integer lifts. Set
	\[
	\mathcal B_y=\{B_\pi(a)=\mathbb E_\pi[z\mid a]:\pi\in\mathfrak P_y\}
	\subset L^2(\mathcal A,\ell;\R^3).
	\]
	The set $\mathcal B_y$ is nonempty, convex, and weakly compact, and
	\begin{equation}\label{eq:exact-marginal-subgradient}
	\partial G(y)=\mathcal B_y.
	\end{equation}
	Thus every subgradient is represented by an optimal plan with the prescribed target marginal.
\end{lemma}

\begin{proof}
	Choose a compact cube $Q\subset\R^3$ containing all nearest lifts to
	the essential range of $y$. The compact set contains lifted targets,
	not just their barycenters: the exact target law must be preserved before
	taking a conditional expectation. In general the barycenter itself
	does not have the prescribed target law. Consider probability kernels $\kappa_a\in\mathcal P(Q)$ satisfying
	\begin{equation}\label{eq:kernel-exact-target}
	\begin{aligned}
	\pi(\dd a,\dd z)&=\ell(\dd a)\kappa_a(\dd z),\\
	\int_{\mathcal A}\int_Q\zeta(qz)\kappa_a(\dd z)\dd\ell(a)
	&=c\int_\T\zeta\,dx\quad(\zeta\in C(\T)).
	\end{aligned}
	\end{equation}
	Let $\mathfrak P_Q$ denote this set. It is nonempty: assign the uniform torus law to each label and choose nearest lifts measurably, with a fixed rule for ties.
	
	Set $\mathscr E_Q=L^1(\mathcal A,\ell;C(Q))$ and associate to a kernel the functional
	\[
	L_\pi(\Phi)=\int_{\mathcal A}\int_Q\Phi(a,z)\kappa_a(\dd z)\dd\ell(a),
	\qquad |L_\pi(\Phi)|\leq\|\Phi\|_{\mathscr E_Q}.
	\]
	Finite sums $f(a)g(z)$ are dense in $\mathscr E_Q$, and separability gives weak-star sequential compactness of the unit dual ball. A limit functional again defines a kernel: on a countable dense subspace of $C(Q)$, the maps $f\mapsto L(fg)$ yield measurable $k_g\in L^\infty(\ell)$ satisfying positivity, $k_1=1$, and $|k_g|\leq\|g\|_\infty$. Extension and the Riesz representation theorem give measurable probability measures $\kappa_a$. The marginal identities
	\[
	L_\pi(f\,1)=\int_{\mathcal A}f\dd\ell,
	\qquad
	L_\pi(\zeta\circ q)=c\int_\T\zeta\,dx
	\]
	are preserved in the limit. Hence $\mathfrak P_Q$ is convex and weak-star compact.
	
	To impose the nearest-lift constraint, define
	\begin{equation}\label{eq:kernel-nearest-lift-defect}
	d_y(a,z)=|y(a)-z|^2-d_\T(qy(a),qz)^2\geq0.
	\end{equation}
	Both $d_y$ and $c_y(a,z)=|y(a)-z|^2$ belong to $\mathscr E_Q$. The nearest-lift constraint $L_\pi(d_y)=0$ is closed in this weak-star
	topology, and the cost $L_\pi(c_y)$ attains its minimum on
	$\mathfrak P_Q$. Measurability in $a$ and continuity in $z$ suffice:
	the topology is defined against $L^1(\mathcal A;C(Q))$, so no continuity
	of $y$ on the label space is being assumed. The nonnegative defect
	turns the potentially discontinuous choice of a nearest lift at ties
	into a closed integral constraint.
	
	Replacing lifts by nearest ones preserves the marginals and decreases the cost by $L_\pi(d_y)$. Thus minimizers are supported on nearest lifts. Disintegration of a torus coupling and measurable lift selection give
	\[
	\min_{\pi\in\mathfrak P_Q}L_\pi(c_y)
	=W_{2,\T}^2(q_\#y_\#\ell,c\,dx).
	\]
	These minimizers form $\mathfrak P_y$. Since their lifts lie in a fixed bounded set, $\mathfrak P_y$ is nonempty, convex, and compact, independently of any further enlargement of $Q$.
	
	For $W\in\mathscr H$, we have $W(a)\cdot z\in\mathscr E_Q$ and
	\[
	\langle B_\pi,W\rangle=L_\pi(W\cdot z),\qquad
	\|B_\pi\|_{\mathscr H}\leq c^{1/2}\max_{z\in Q}|z|.
	\]
	Thus the barycenter map is weak-star-to-weak continuous, and $\mathcal B_y$ is convex and weakly compact.
	
	Writing
	\[
	\Lambda_w(\pi)=\int\left(w(a)\cdot z-\tfrac12|z|^2\right)\dd\pi,
	\]
	optimality gives $G(y)=\Lambda_y(\pi)$ for $\pi\in\mathfrak P_y$, while $G(w)\geq\Lambda_w(\pi)$ for every $w\in\mathscr H$. Hence
	\[
	G(y+W)\geq G(y)+\langle B_\pi,W\rangle,
	\qquad W\in\mathscr H,
	\]
	so $\mathcal B_y\subset\partial G(y)$.
	
	For the reverse inclusion, fix bounded $W\in\mathscr H$ and put $y_\varepsilon=y+\varepsilon W$, $0<\varepsilon\leq1$. Enlarging $Q$ to contain the corresponding nearest lifts gives
	\[
	\|c_{y_\varepsilon}-c_y\|_{\mathscr E_Q}
	+\|d_{y_\varepsilon}-d_y\|_{\mathscr E_Q}
	\leq C\varepsilon,
	\]
	with $C$ depending on $c,Q,\|y\|_\infty,\|W\|_\infty$. For $\pi_\varepsilon\in\mathfrak P_{y_\varepsilon}$ and any $\pi\in\mathfrak P_y$, optimality gives
	\[
	\begin{aligned}
	L_{\pi_\varepsilon}(c_y)
	&\leq L_{\pi_\varepsilon}(c_{y_\varepsilon})+C\varepsilon
	\leq L_\pi(c_{y_\varepsilon})+C\varepsilon
	\leq L_\pi(c_y)+2C\varepsilon,\\
	L_{\pi_\varepsilon}(d_y)&\leq C\varepsilon.
	\end{aligned}
	\]
	Every weak-star cluster point of $\pi_\varepsilon$ therefore lies in $\mathfrak P_y$, with both marginals and the nearest-lift constraint preserved.
	
	The affine maximization formula yields
	\begin{equation}\label{eq:exact-marginal-directional-squeeze}
	\langle B_\pi,W\rangle
	\leq\frac{G(y+\varepsilon W)-G(y)}{\varepsilon}
	\leq\langle B_{\pi_\varepsilon},W\rangle.
	\end{equation}
	Choose $\pi\in\mathfrak P_y$ maximizing
	$\langle B_\pi,W\rangle$; the maximum exists by weak compactness.
	For the upper bound, take a subsequence of $\pi_\varepsilon$ realizing
	the upper limit and use its weak-star cluster point in $\mathfrak P_y$.
	The two sides of \eqref{eq:exact-marginal-directional-squeeze} then
	coincide. Convexity of $G$ ensures existence of the one-sided derivative,
	and we obtain
	\begin{equation}\label{eq:exact-marginal-directional-derivative}
	G'(y;W)=\max_{B\in\mathcal B_y}\langle B,W\rangle.
	\end{equation}
	If $B_0\in\partial G(y)\setminus\mathcal B_y$, separation provides $W\in\mathscr H$ such that $\langle B_0,W\rangle>\sup_{B\in\mathcal B_y}\langle B,W\rangle$. Since $\mathcal B_y$ is bounded, truncation preserves this strict separation with bounded $W$. The subgradient inequality then contradicts \eqref{eq:exact-marginal-directional-derivative}, proving \eqref{eq:exact-marginal-subgradient}.
\end{proof}

For $x=X(a)$ and $y_*=Y_K(X)$, the optimal plan $\pi_*$ of \eqref{eq:periodic-residual} also minimizes the mixed quadratic cost
\begin{equation}\label{dp:barycentric-cost}
\mathcal J_K(\pi)=\frac12\int\left[
\frac{K-1}{K}|z-x|^2
+\frac1K|\mathbb E[z\mid x]-x|^2\right]\dd\pi(x,z),
\end{equation}
over plans with prescribed torus target marginal and variable integer lifts. Set
\[
L(y,\pi)=\frac{K-1}{2}\|y\|^2-K\langle x,y\rangle
+\int y\cdot z\dd\pi+\frac K2\|x\|^2
-\frac12\int|z|^2\dd\pi.
\]
Optimality of $\pi_*$ and $(K-1)y_*-Kx+\mathbb E_{\pi_*}[z\mid x]=0$ give a saddle point $(y_*,\pi_*)$. Completing the square yields
\[
\inf_y L(y,\pi)=-\frac{K}{K-1}\mathcal J_K(\pi),
\]
and the saddle-point inequalities yield
\[
\inf_yL(y,\pi)\leq L(y_*,\pi)
\leq L(y_*,\pi_*)=\inf_yL(y,\pi_*).
\]
Hence $\pi_*$ minimizes $\mathcal J_K$, with the prescribed target
marginal. The first term of $\mathcal J_K$ is linear in the plan and
the second is strictly convex in its conditional barycenter. If two
minimizers had barycenters $T_1\ne T_2$ in $L^2(\mu)$, their equally
weighted mixture would lower the cost by
\[
\frac{1}{8K}\int|T_1-T_2|^2\dd\mu>0.
\]
Thus the barycenter is unique even if the optimal plan is not.

\section{Auxiliary Entropy and High-Frequency Estimates}
\label{app:auxiliary-estimates}

Lemma~\ref{lem:vacuum-entropy} supplies the initialized identities used
in the entropy comparison of Section~\ref{sec:scalar}.
Lemma~\ref{lem:kernel} supplies the common kernel for the truncated
continuity source and convective commutator in
Lemma~\ref{lem:uniform-high-frequency}. These statements concern
different limiting operations: removal of the vacuum regularization
at fixed truncation level, and decay as the spatial smoothing scale
tends to zero. In the latter case, boundedness of the normalized
truncation and its renormalization coefficient makes the same kernel
estimate uniform in the truncation level. In the former, the truncation
level stays fixed until the initialized entropy identities have been
obtained and compared.

\begin{lemma}[Initialized truncated entropy identity]
	\label{lem:vacuum-entropy}
	Let $(\varrho,v)$ satisfy the initialized renormalized continuity equation
	on $Q_T$, with $\varrho\geq0$, $\varrho\in L^\infty(0,T;L^\gamma(\T))$,
	$v\in L^2(0,T;H^1(\T))$, and $\gamma>1$.
	For every $k\geq1$ and $\psi\in C_c^\infty([0,T))$,
	\[
	\begin{aligned}
	&-\int_0^T\psi'(t)\int_\T L_k(\varrho)\,dx\,dt
	+\int_0^T\psi(t)\int_\T T_k(\varrho)\Div v\,dx\,dt\\
	&\qquad=\psi(0)\int_\T L_k(\varrho_0)\,dx.
	\end{aligned}
	\]
	The identity applies to the approximating and limiting densities in
	Section~\ref{sec:scalar} and yields
	\eqref{eq:entropy-member-limit}--\eqref{eq:entropy-limit-companion}
	after passage to the weak limits.
\end{lemma}

\begin{proof}
	Fix $k\geq1$. Since $L_k$ is affine on $[3k,\infty)$, its slope is
	\[
	\alpha_k=\int_1^{3k}T_k(s)s^{-2}\dd s+\frac23.
	\]
	Choose a smooth nonnegative $w$ on $[0,\infty)$, supported in $[0,1]$,
	constant and positive near zero, with $\int_0^1w=1$. Set
	$\omega(s)=\int_0^s w(r)\dd r$ and, for $0<\varepsilon<1$, define
	\[
	\begin{gathered}
	a_{k,\varepsilon}(s)=T_k(s)\omega(s/\varepsilon),\\
	L_{k,\varepsilon}(z)=z\int_1^z a_{k,\varepsilon}(s)s^{-2}\dd s
	\quad(z>0),\qquad L_{k,\varepsilon}(0)=0.
	\end{gathered}
	\]
	The algebraic identities are
	\[
	zL_{k,\varepsilon}'(z)-L_{k,\varepsilon}(z)
	=a_{k,\varepsilon}(z),\qquad
	L_{k,\varepsilon}''(z)=\frac{a_{k,\varepsilon}'(z)}z\geq0
	\quad(z>0).
	\]
	Here $T_k$ and $\omega$ are nonnegative and nondecreasing. Near zero,
	$a_{k,\varepsilon}(z)=O(z^2/\varepsilon)$, so the first derivative
	extends continuously to zero. Thus the regularization is compatible
	with vacuum and preserves convexity. Then
	$L_{k,\varepsilon}-\alpha_k z$ is bounded and $C^1$, with derivative
	vanishing above $3k$. Moreover,
	\[
	\begin{gathered}
	0\leq L_{k,\varepsilon}(z)-L_k(z)
	\leq z\log(\varepsilon/z)\ind_{\{0<z<\varepsilon\}}
	\leq\frac{\varepsilon}{\exp(1)},\\
	|a_{k,\varepsilon}(z)-T_k(z)|\leq\varepsilon.
	\end{gathered}
	\]
	Writing $d_n=\Div u_n$ and
	$D_T=\sup_n\|d_n\|_{L^2(Q_T)}$, the source,
	time-derivative, and initial errors for a time test $\psi$ are bounded by
	\[
	\varepsilon\|\psi\|_\infty |Q_T|^{1/2}D_T,\qquad
	\frac{\varepsilon}{\exp(1)}\|\psi'\|_1,\qquad
	\frac{\varepsilon}{\exp(1)}|\psi(0)|.
	\]
	Apply initialized renormalization to
	$L_{k,\varepsilon}-\alpha_k z$ and add $\alpha_k$ times the continuity
	equation. For the approximating solution,
	\[
	\begin{aligned}
	&-\int_0^T\psi'(t)\int_\T L_{k,\varepsilon}(\rho_n)\,dx\,dt
	+\int_0^T\psi(t)\int_\T a_{k,\varepsilon}(\rho_n)d_n\,dx\,dt\\
	&\qquad=\psi(0)\int_\T L_{k,\varepsilon}(\rho_{0n})\,dx.
	\end{aligned}
	\]
	The same argument applies to $(\rho,u,\rho_0)$.
	Letting $\varepsilon\downarrow0$ proves the stated identity,
	including its initial term. The uniform error bounds and the weak limits of
	$L_k(\rho_n)$ and $T_k(\rho_n)d_n$ then give
	\eqref{eq:entropy-member-limit}; the limiting equation gives
	\eqref{eq:entropy-limit-companion}.
\end{proof}

Let $A=\nabla\Delta_\T^{-1}$ and $R_{ij}=\partial_j\partial_i\Delta_\T^{-1}$. Choose an even $\chi\in C_c^\infty(B_1)$ with $\chi\geq0$ and $\int\chi=1$, and use the Euclidean Fourier convention $\widehat\chi(\xi)=\int_{\R^3}e^{-2\pi\mathrm i x\cdot\xi}\chi(x)\dd x$. Define $S_a$ by convolution with the periodization of $\chi_a=a^{-3}\chi(\cdot/a)$, where $0<a<1/8$, and set $H_a=I-S_a$, $A_{a,i}=A_iH_a$, and $L_{ij,a}=R_{ij}H_a$. Since $H_a$ annihilates constants,
\begin{equation}\label{eq:zero-mode}
\overline{H_a f}=0,\qquad \sum_iL_{ii,a}=H_a.
\end{equation}
The unfiltered identity is $\sum_iR_{ii}f=f-\bar f$.

\begin{lemma}\label{lem:kernel}
	The periodic kernel $\mathcal A_a$ of $A_a$ satisfies $\|\mathcal A_a\|_1\leq Ca$. There is a nonnegative periodic kernel $\mathcal K_a$ such that
	\begin{equation}\label{eq:kernel-norm}
	\|\mathcal K_a\|_s\leq C_s(a^{3/s-2}+a^2),
	\qquad 1\leq s<3/2.
	\end{equation}
	If $g\in L^\infty$, $u\in H^1(\T)$, and $|h|\leq C_h|\nabla u|$, then
	\begin{equation}\label{eq:field-bound}
	|L_{ij,a}(gu_j)-u_jL_{ij,a}g+A_{a,i}h|
	\leq C(\|g\|_\infty+C_h)\mathcal K_a*_\T|\nabla u|
	\quad\hbox{almost everywhere}.
	\end{equation}
The source $h$ need not have zero spatial mean.
\end{lemma}

\begin{proof}
	The local potential and principal-value kernels are those of the
	classical singular-integral calculus \cite{Stein}. We retain the
	contact term and periodize the integrable high-frequency difference,
	so that the periodic zero mode is fixed by the construction itself.
	Begin in $\R^3$ with $A^\R(z)=z/(4\pi|z|^3)$ and form the integrable difference $\mathcal A_a^\R=A^\R-A^\R*\chi_a$ for periodization. Write $r=|z|$. For $r\leq2a$, the difference is bounded by $Cr^{-2}$. For $r>2a$, evenness of $\chi_a$ cancels the first-order Taylor term, giving
	\[
	\mathcal A_a^\R(z)
	=-\int\chi_a(y)\int_0^1(1-t)D^2A^\R(z-ty)[y,y]\dd t\dd y,
	\qquad |\mathcal A_a^\R(z)|\leq Ca^2r^{-4}.
	\]
	Thus
	\begin{equation}\label{eq:highpass-potential-cost}
	\|\mathcal A_a^\R\|_1
	\leq C\left(\int_0^{2a}1\dd r+a^2\int_{2a}^{\infty}r^{-2}\dd r\right)
	\leq Ca.
	\end{equation}
	Since $\mathcal A_a^\R$ is integrable and odd, its periodization has zero Fourier coefficient at frequency zero. Its nonzero coefficients are $-\mathrm i\ell[1-\widehat\chi(a\ell)]/(2\pi|\ell|^2)$, those of $A_a$. It is therefore the periodic kernel of $A_a$.
	
	Differentiation of the Euclidean kernel produces a Dirac term at the origin as well as a principal-value kernel:
	\[
	\left\{
	\begin{aligned}
	\partial_j A_i^{\mathbb R}
	&=\operatorname{pv}F_{ij}
	+\frac{\delta_{ij}}{3}\delta_0,\\[2mm]
	F_{ij}(z)
	&=\frac{\delta_{ij}|z|^2-3z_iz_j}{4\pi|z|^5}.
	\end{aligned}
	\right.
	\]
	The corresponding high-frequency part is $\operatorname{pv}F_{ij}-(\operatorname{pv}F_{ij})*\chi_a+\frac{\delta_{ij}}3(\delta_0-\chi_a)$. The Dirac term cancels in the commutator because the velocity difference vanishes at the origin, but the convolved term involving $\chi_a$ must remain. Away from the origin, the required kernel is
	\[
	W_{ij,a}=F_{ij}-(\operatorname{pv}F_{ij})*\chi_a
	-\frac{\delta_{ij}}3\chi_a,\qquad
	|W_a(z)|\leq
	\begin{cases}
	C(r^{-3}+a^{-3}),&r\leq2a,\\
	Ca^2r^{-5},&r>2a.
	\end{cases}
	\]
	For $r\leq2a$, rescale to $a=1$. Since $F_{ij}$ has zero integral on every spherical shell, the singular part of the convolution can be written as $\int_{|y|<1}F_{ij}(y)[\chi(z-y)-\chi(z)]\dd y$. The difference estimate for the smooth kernel bounds its absolute value by $C\int_0^1\dd r$. For $r>2a$, evenness of $\chi_a$ again cancels the first-order Taylor term, and the point lies outside the support of $\chi_a$.
	
	First let $u$ be smooth and periodic and $g$ bounded and periodic, and extend both periodically to $\R^3$. The commutator is
	\[
	[L_{ij,a},u_j]g(x)
	=\int_{\R^3}W_{ij,a}(z)[u_j(x-z)-u_j(x)]g(x-z)\dd z.
	\]
	Apply the line-segment formula to $u_j(x-z)-u_j(x)$ and change variables $w=tz$. The absolute value of the commutator is bounded by
	\[
	\|g\|_\infty\int_{\R^3}|\nabla u(x-w)|\,|w|
	\int_0^1|W_a(w/t)|t^{-4}\dd t\dd w.
	\]
	For $r=|w|\leq2a$, split the integral at $t_0=r/(2a)$ to obtain the bound
	\[
	Ca^2r^{-4}\int_0^{t_0}t\dd t+
	Cr^{-2}\int_{t_0}^1t^{-1}\dd t+
	Cra^{-3}\int_{t_0}^1t^{-4}\dd t
	\leq Cr^{-2}[1+\log(2a/r)].
	\]
	For $r>2a$, only the far-field estimate occurs, giving $Ca^2r^{-4}$. Adjusting the constant on $a<r<2a$, we may take
	\[
	K_a(w)=
	\begin{cases}
	|w|^{-2}[1+\log(a/|w|)],&0<|w|\leq a,\\
	a^2|w|^{-4},&|w|>a.
	\end{cases}
	\]
	The near- and far-field estimates also bound the potential kernel $\mathcal A_a^\R$ by $CK_a$, giving $|A_ah|\leq CC_hK_a*_{\R^3}|\nabla u|$. Thus $K_a$ controls both the continuity source and the convective commutator.
	
	Set $\mathcal K_a(x)=\sum_{\ell\in\mathbb Z^3}K_a(x+\ell)$. Radial integration gives
	\[
	\|K_a\|_{L^s(\R^3)}^s
	=4\pi a^{3-2s}
	\left[\int_0^1r^{2-2s}(1+|\log r|)^s\dd r
	+\int_1^\infty r^{2-4s}\dd r\right].
	\]
	The first integral is finite precisely when $2-2s>-1$, or $s<3/2$; the second is finite for $s\geq1$. On the centered fundamental cube, nonzero lattice points satisfy $|x+\ell|\geq c|\ell|$, so the remaining periodic terms are bounded by $Ca^2\sum_{\ell\ne0}|\ell|^{-4}\leq Ca^2$. This proves \eqref{eq:kernel-norm}, including the contribution from periodization.
	
	Finally, approximate $u\in H^1$ by smooth periodic functions. Young's inequality passes the convolution bound to the limit, and boundedness of $R_{ij}$ on finite-exponent $L^r$ spaces identifies the distributional limit. An almost-everywhere convergent subsequence gives \eqref{eq:field-bound}.
\end{proof}

\vspace{8mm}

\section*{Acknowledgements}

The research of G.Y. Hong was supported by the National Natural Science Foundation of China (Grant No. 12522111), the Guangdong Provincial Pearl River Talents Program (Grant No. 2023QN10X436), the Guangdong Basic and Applied Basic Research Foundation (Grant No. 2024A1515012306), and the Fundamental Research Funds for the
Central Universities (Grant No. 2025ZYGXZR040). The research of H.Y. Peng was supported by the National Natural Science Foundation of China $\#$12271112
and the Basic and Applied Basic Research Project of Guangzhou $\#$SL2023A04J01554. The research of C.J. Zhu was supported by the National Natural Science Foundation of China $\#$12571233 and the Basic and Applied Basic Research Foundation of Guangdong Province $\#$2026A1515012081.

\vspace{6mm}

\noindent {\bf Data Availability:} No datasets were generated or analysed for this study.

\noindent {\bf Conflicts of Interest:} The authors declare that they have no conflicts of interest.

\noindent {\bf AI declaration:} The authors acknowledge the use of generative AI tools in developing and checking mathematical
arguments and in drafting and revising the manuscript. The authors take full responsibility for the mathematical content.

\enlargethispage{4\baselineskip}
\end{document}